\documentclass{amsart}
\usepackage{amsfonts}
\usepackage{CJK,CJKnumb,CJKulem,times,dsfont,ifthen,mathrsfs,latexsym,amsfonts,color}
\usepackage{amsmath,amsthm,fontenc,amssymb,bm,graphicx,psfrag,listings,curves,extarrows}
\newtheorem{theorem}{Theorem}
\theoremstyle{plain}

\newtheorem{corollary}{Corollary}

\newtheorem{definition}{Definition}

\newtheorem{lemma}{Lemma}

\newtheorem{proposition}{Proposition}
\newtheorem{remark}{Remark}

\numberwithin{equation}{section}

\newcommand{\s}{\section}
\newcommand{\R}{\mathbb R}
\newcommand{\lab}{\label}
\newcommand{\bt}{\begin{theorem}}
\newcommand{\et}{\end{theorem}}
\newcommand{\bl}{\begin{lemma}}
\newcommand{\el}{\end{lemma}}
\newcommand{\bd}{\begin{definition}}
\newcommand{\ed}{\end{definition}}
\newcommand{\bc}{\begin{corollary}}
\newcommand{\ec}{\end{corollary}}
\newcommand{\bp}{\begin{proof}}
\newcommand{\ep}{\end{proof}}
\newcommand{\bo}{\begin{proposition}}
\newcommand{\eo}{\end{proposition}}
\newcommand{\br}{\begin{remark}}
\newcommand{\er}{\end{remark}}
\newcommand{\be}{\begin{equation}}
\newcommand{\ee}{\end{equation}}
\newcommand{\NN}{{\mathbb N}}

\begin{document}
\title[On a double-variable inequality and elliptic
systems]{On a double-variable inequality and elliptic
systems involving critical Hardy-Sobolev exponents}
\author{X. Zhong}
\address[X. Zhong and W. Zou]{Department of Mathematics, Sun Yat-Sen University Guangzhou  510275, China.\\
Department of Mathematical Sciences, Tsinghua University,Beijing 100084, China.}
\email[X. Zhong]{zhongxuexiu1989@163.com}
\author{W. Zou}
\email[W. Zou]{wzou@math.tsinghua.edu.cn}
\thanks{Supported by NSFC}
\date{March 3, 2016}
\subjclass[2010]{Primary 35B38, 35J10, 35J15, 35J20, 35J60; Secondary 49J40}
\keywords{Elliptic system, sharp constant, Hardy-Sobolev exponent, existence, nonexistence, ground state solution, infinitely  many sign-changing solutions}

\begin{abstract}
Let $\Omega\subset \R^N$ ($N\geq 3$)  be an  open domain which is not necessarily bounded.  The sharp constant and extremal functions to the following kind of double-variable inequalities
$$ S_{\alpha,\beta,\lambda,\mu}(\Omega)  \Big(\int_\Omega \big(\lambda \frac{|u|^{2^*(s)}}{|x|^s}+\mu \frac{|v|^{2^*(s)}}{|x|^s}+2^*(s)\kappa \frac{|u|^\alpha |v|^\beta}{|x|^s}\big)dx\Big)^{\frac{2}{2^*(s)}}$$
 $$\leq  \int_\Omega \big(|\nabla u|^2+|\nabla v|^2\big)dx$$
for $(u,v)\in {\mathscr{D}}:=D_{0}^{1,2}(\Omega)\times D_{0}^{1,2}(\Omega)$ will be explored. Further results about the sharp constant $S_{\alpha,\beta,\lambda,\mu}(\Omega)$ with its extremal functions when $\Omega$ is a general open domain will be involved. For this goal, we consider the following elliptic systems involving multiple Hardy-Sobolev critical exponents:
$$\begin{cases}
-\Delta u-\lambda \frac{|u|^{2^*(s_1)-2}u}{|x|^{s_1}}=\kappa\alpha \frac{1}{|x|^{s_2}}|u|^{\alpha-2}u|v|^\beta\quad &\hbox{in}\;\Omega,\\
-\Delta v-\mu \frac{|v|^{2^*(s_1)-2}v}{|x|^{s_1}}=\kappa\beta \frac{1}{|x|^{s_2}}|u|^{\alpha}|v|^{\beta-2}v\quad &\hbox{in}\;\Omega,\\
(u,v)\in \mathscr{D}:=D_{0}^{1,2}(\Omega)\times D_{0}^{1,2}(\Omega),
\end{cases}$$
where $s_1,s_2\in (0,2), \alpha>1,\beta>1, \lambda>0,\mu>0,\kappa\neq 0, \alpha+\beta\leq 2^*(s_2)$. Here,  $2^*(s):=\frac{2(N-s)}{N-2}$ is the critical Hardy-Sobolev exponent.  We mainly study the  critical case (i.e.,  $\alpha+\beta=2^*(s_2)$) when  $\Omega$ is a cone (in particular,  $\Omega=\R_+^N$ or $\Omega=\R^N$).
We will establish a sequence of fundamental results including regularity, symmetry, existence and multiplicity, uniqueness and nonexistence, {\it etc.}
\end{abstract}

\maketitle
\tableofcontents
\s{Introduction}
\renewcommand{\theequation}{1.\arabic{equation}}
\renewcommand{\theremark}{1.\arabic{remark}}
\renewcommand{\thedefinition}{1.\arabic{definition}}

Let   $\Omega\subset \R^N$ ($N\geq 3$)  be  an open domain which is  not necessarily bounded. We study the following nonlinear elliptic systems
\be\lab{P}
\begin{cases}
-\Delta u-\lambda \frac{|u|^{2^*(s_1)-2}u}{|x|^{s_1}}=\kappa\alpha \frac{1}{|x|^{s_2}}|u|^{\alpha-2}u|v|^\beta\quad &\hbox{in}\;\Omega,\\
-\Delta v-\mu \frac{|v|^{2^*(s_1)-2}v}{|x|^{s_1}}=\kappa\beta \frac{1}{|x|^{s_2}}|u|^{\alpha}|v|^{\beta-2}v\quad &\hbox{in}\;\Omega,\\
(u,v)\in \mathscr{D}:=D_{0}^{1,2}(\Omega)\times D_{0}^{1,2}(\Omega),
\end{cases}
\ee
where $s_1,s_2\in (0,2), \alpha>1,\beta>1, \lambda>0,\mu>0,\kappa\neq 0, \alpha+\beta\leq 2^*(s_2):=\frac{2(N-s_2)}{N-2}$.

The interest in studying the  nonlinear Schr\"odinger systems is motivated by  real problems
in  nonlinear optics,  plasma physics,  condensed matter physics,  etc.
For example, the coupled nonlinear Schr\"odinger systems arise in the description
of several physical phenomena such as the propagation of pulses in birefringent
optical fibers and Kerr-like photorefractive media, see \cite{AkhmedievAnkiewicz.1999,EvangelidesMollenauerGordonBergano.1992,Kaminow.1981,Menyuk.1987,Menyuk.1989,WaiMenyukChen.1991}, etc.
The problem comes from the physical phenomenon with a clear practical significance. The researches  on solutions under different situations not only corresponds  to different
physical interpretation, but also has a pure mathematical theoretical significance.
Hence, the coupled nonlinear Schr\"odinger systems are widely
studied in recently years, we refer the readers to \cite{AbdellaouiFelliPeral.2009,AmbrosettiColorado.2006,LinWei.2005,
MaiaMontefuscoPellacci.2006,Sirakov.2007}
and the
references therein.

For any $s\in [0,2]$, we define the measure $d\mu_s:=\frac{1}{|x|^s}dx$ and $\displaystyle|u|_{p,s}^{p}:=\int_{\Omega}|u|^pd\mu_s.$ We also use the notation $\displaystyle|u|_p:=|u|_{p,0}.$
The Hardy-Sobolev inequality\cite{CaffarelliKohnNirenberg.1984,CatrinaWang.2001,GhoussoubYuan.2000}
asserts that $D_{0}^{1,2}(\R^N)\hookrightarrow L^{2^*(s)}(\R^N,d\mu_s)$ is a continuous  embedding for $s\in [0,2]$. For a general open domain $\Omega$, there exists a positive constant $C(s,\Omega)$ such that
$$\int_\Omega |\nabla u|^2 dx\geq C(s, \Omega)|u|_{2*(s),s}^{2}, \;u\in D_{0}^{1,2}(\Omega).$$
Define
$\mu_{s_1}(\Omega)$ as
\be\lab{2014-2-27-e2}
\mu_{s_1}(\Omega):=\inf\left\{\frac{\int_\Omega |\nabla u|^2 dx}{\left(\int_\Omega \frac{|u|^{2^*(s_1)}}{|x|^{s_1}}\right)^{\frac{2}{2^*(s_1)}}}
\;:\;u\in D_{0}^{1,2}(\Omega)\backslash\{0\}\right\}
\ee
Consider the case of $\Omega=\R_+^N$, it is well known that
the extremal function of $\mu_{s_1}(\R_+^N)$ is  parallel to the ground state solution of the following problem:
\be\lab{BP-0}
\begin{cases}
-\Delta u=\frac{|u|^{2^*(s_1)-2}u}{|x|^{s_1}}\quad \hbox{in}\;\R_+^N,\\
u=0\;\hbox{on}\;\partial\R_+^N.
\end{cases}
\ee
We note that the existence of ground state solution of (\ref{BP-0}) for $0<s_1<2$ is solved by Ghoussoub and Robert \cite{GhoussoubRobert.2006}.
They also gave   some properties about the regularity, symmetry and decay estimates.
The instanton $U(x):=C\big(\kappa+|x|^{2-s_2}\big)^{-\frac{N-2}{2-s_2}}$  for  $0\leq s_2<2$  is a ground state solution to (\ref{BPP}) below (see \cite{Lieb.1983} and \cite{Talenti.1976a}):
\be\lab{BPP}
\begin{cases}
\Delta u+\frac{u^{2^*(s_2)-1}}{|x|^{s_2}}=0\quad \hbox{in}\;\R^N,\\
u>0\;\;\hbox{in}\;\R^N\;\;\;\hbox{and}\; u\rightarrow 0\;\hbox{as}\;|x|\rightarrow +\infty.
\end{cases}
\ee
The case that
$0\in \partial\Omega$ has become an  interesting  topic in recent years since the curvature of $\partial \Omega$ at $0$ plays an important role, see \cite{ChernLin.2010,GhoussoubKang.2004,GhoussoubRobert.2006,HsiaLinWadade.2010}, etc.
A lot of sufficient conditions are given in order  to ensure that
 $\mu_{s_1}(\Omega)<\mu_{s_1}(\R_+^N)$ in those papers. It is standard to apply the blow-up analysis to  show that $\mu_{s_1}(\Omega)$ can be achieved by some positive $u\in H_0^1(\Omega) $  (e.g., see \cite[Corollary 3.2]{GhoussoubKang.2004}), which is a ground state solution of $$\begin{cases}
-\Delta u=\frac{|u|^{2^*(s_1)-2}u}{|x|^{s_1}}\quad &\hbox{in}\;\Omega,\\
u=0\;&\hbox{on}\;\partial\Omega,
\end{cases}
$$
and the least energy  equals $(\frac{1}{2}-\frac{1}{2^*(s_1)})\mu_{s_1}(\Omega)^{\frac{N-s_1}{2-s_1}}$.

However, it seems there is   no article before involving  the system case  like   \eqref{P}  with  Hardy-Sobolev critical exponents,  which we are going to deal with in the current paper. It is well known that the main difficulty is the lack of compactness inherent in these problems involving Hardy-Sobolev critical exponents. The compactness concentration argument (see \cite{Lions.1985a}, etc.) is a powerful tool to handle with these critical problems. It is also well known that the compactness concentration argument depends heavily on the limit problem.
Consider a bounded domain $\Omega$, if $0\not\in \bar{\Omega}$, we see that $\frac{1}{|x|^{s_i}},i=1,2$ are regular. We are interested in the case of that $0\in \bar{\Omega}$.
It is easy to see that when $0\in \Omega$, the limit domain is $\R^N$, and when $0\in \partial\Omega$, the limit domain is usually a cone. Especially,  when $\partial\Omega$ possesses a suitable regularity (e.g.  $\partial\Omega\in C^2$ at $x=0$), the limit domain is $\R_+^N$ after a suitable rotation. Hence, in present paper, we mainly study the critical elliptic systems \eqref{P} with $\alpha+\beta=2^*(s_2)$ and $\Omega$ is a cone.

\bd\lab{cone}  A cone in $\R^N$ is an open  domain $\Omega$ with Lipschitz boundary and such that $tx\in \Omega$ for every $t>0$ and $x\in \Omega$.
\ed
We will establish a sequence of fundamental   results  to the system \eqref{P} including regularity, symmetry, existence and multiplicity, and nonexistence, {\it etc.}
Since there are a large number of conclusions in the current paper, we do not intend to list them here. This paper is organized as follows:

\noindent
{\bf In Section 2}, we will establish by a direct method  a type of interpolation inequalities,  which are essentially the variant  Caffarelli-Kohn-Nirenberg  (CKN)  inequalities, see \cite{CaffarelliKohnNirenberg.1984}.

\noindent
{\bf In Section 3}, we will study the regularity, symmetry and decay estimation about the nonnegative solutions of \eqref{P}. Taking $\R_+^N$ as a specific example, we will study the regularity based on the  technique  of Moser's iteration (see Proposition \ref{2014-11-21-prop1}). By the method of moving planes, we obtain the symmetry result (see Proposition \ref{2014-11-22-Prop1}). Due to the Kelvin transformation, we get the decay estimation (see Proposition \ref{2014-11-21-Prop2}).

\noindent
{\bf In Section 4}, we shall study the basic properties of the corresponding Nehari manifold.

\noindent
{\bf In Section 5}, we will give a nonexistence of nontrivial ground state solution of \eqref{P} for the case $s_2\geq s_1$, see Theorem \ref{2014-12-1-th1}.

\noindent
{\bf In Section 6}, we will give an existence result of positive solution result for a special case : $\lambda=\mu (\frac{\beta}{\alpha})^{\frac{2^*(s_1)-2}{2}}$, see Corollary \ref{2014-11-23-cro1}.  Further, we prepare  a sequence of preliminaries for the existence result which are not only useful for us to study the case of $s_1=s_2$  in Section 7, but also the case of $s_1\neq s_2$ in  Section 8.

\noindent
{\bf In Section 7}, we will focus on the case of $s_1=s_2=s\in (0,2)$ when $\Omega$ is a cone. In this case, the nonlinearities are homogeneous which enable us to define the following constant
\be\lab{2014-12-1-e1}
S_{\alpha,\beta,\lambda,\mu}(\Omega):=\inf_{(u,v)\in \widetilde{\mathscr{D}}} \frac{\int_\Omega \big(|\nabla u|^2+|\nabla v|^2\big)dx}{\Big(\int_\Omega \big(\lambda \frac{|u|^{2^*(s)}}{|x|^s}+\mu \frac{|v|^{2^*(s)}}{|x|^s}+2^*(s)\kappa \frac{|u|^\alpha |v|^\beta}{|x|^s}\big)dx\Big)^{\frac{2}{2^*(s)}}},
\ee
where
\be\lab{2014-12-1-e2}
\widetilde{\mathscr{D}}:=\left\{(u,v)\in \mathscr{D}:\;\int_\Omega \left(\lambda \frac{|u|^{2^*(s)}}{|x|^s}+\mu \frac{|v|^{2^*(s)}}{|x|^s}+2^*(s)\kappa \frac{|u|^\alpha |v|^\beta}{|x|^s}\right)dx>0\right\}.
\ee
In particular, we shall see that $\widetilde{\mathscr{D}}=\mathscr{D}\backslash \{(0,0)\}$   if and only if $$ \kappa>-\left(\frac{\lambda}{\alpha}\right)^{\frac{\alpha}{2^*(s)}}
\left(\frac{\mu}{\beta}\right)^{\frac{\beta}{2^*(s)}},$$
see Lemma  \ref{2014-11-26-l2}. When $\kappa<0$, we will prove that $S_{\alpha,\beta,\lambda,\mu}(\Omega)$ has no nontrivial extremals (see Lemma \ref{2014-11-26-l3}). Hence,   we will mainly   focus on the case of $\kappa>0$ and show that the system  \eqref{P} possesses a least energy solution and that $S_{\alpha,\beta,\lambda,\mu}(\Omega)$ is achieved. These conclusions  will produce the  sharp constant and extremal functions to the following kind of  inequalities with double-variable
 $$ S_{\alpha,\beta,\lambda,\mu}(\Omega)  \Big(\int_\Omega \big(\lambda \frac{|u|^{2^*(s)}}{|x|^s}+\mu \frac{|v|^{2^*(s)}}{|x|^s}+2^*(s)\kappa \frac{|u|^\alpha |v|^\beta}{|x|^s}\big)dx\Big)^{\frac{2}{2^*(s)}}$$
 $$\leq  \int_\Omega \big(|\nabla u|^2+|\nabla v|^2\big)dx$$
for $(u,v)\in {\mathscr{D}}.$  For this purpose,  a kind of Pohozaev identity will be established.  Then the existence,  regularity, uniqueness and nonexistence  results of the  positive ground state solution to the system  \eqref{P} can be seen in this section.
Under some proper hypotheses,  we will  show that the positive ground state solution  must be of the form   $\big(C(t_0)U, t_0C(t_0)U\big)$, where $t_0>0$ and $C(t_0)$ can be formulated  explicitly and $U$ is the ground state solution of
$$\begin{cases}
-\Delta u=\mu_{s}(\Omega) \frac{u^{2^*(s)-1}}{|x|^{s}}\;\quad&\hbox{in}\;\Omega,\\
u=0\;\;&\hbox{on}\;\partial \Omega.
\end{cases}
$$
Taking a special case $N=3,s=1,\alpha=\beta=2,\lambda=\mu=2\kappa$ in consideration, we will find  out all the positive ground state solutions to \eqref{P} .
Based on these conclusions, we may prove  the existence of infinitely many sign-changing solutions of the system  \eqref{P} on a cone $\Omega$ by gluing together suitable signed solutions corresponding to each sub-cone.  Further,   if  $\Omega$ is a general open domain,    the  sharp constant  $S_{\alpha,\beta,\lambda,\mu}(\Omega)$ and its extremal functions will be investigated.  We will find  a way to compute
   $S_{\alpha,\beta,\lambda,\mu}(\Omega)$ and     to judge
   when $S_{\alpha,\beta,\lambda,\mu}(\Omega)$  can be  achieved  if $\Omega$
   is a general open domain.

\noindent
{\bf In Section 8,}  the system  \eqref{P}  satisfying    $s_1\not=s_2\in (0,2)$    will be studied.  We shall consider  a new approximation to the original system  \eqref{P}.
 The estimation on the least energy  and the positive ground state along with its geometric structure to   the approximation   will be established.
 Finally, the existence of positive ground state solution to the original system will be given.

\vskip0.2in


\s{Interpolation inequalities}
\renewcommand{\theequation}{2.\arabic{equation}}
\renewcommand{\theremark}{2.\arabic{remark}}
\renewcommand{\thedefinition}{2.\arabic{definition}}
For $s_1\neq s_2$,  we note that there is no embedding relationship between $L^{2^*(s_1)}(\Omega,d\mu_{s_1})$ and $L^{2^*(s_2)}(\Omega,d\mu_{s_2})$ for any domain $\Omega$ with $0\in \bar{\Omega}$.    Hence, we are going to establish some interpolation inequalities in this section.

\bo\lab{2014-4-22-interpolation-wl1}  Let $\Omega\subset \R^N (N\geq 3)$ be an open set.
Assume that $0\leq s_1<s_2<s_3\leq 2$, then there exists  $\theta=\frac{(N-s_1)(s_3-s_2)}{(N-s_2)(s_3-s_1)}\in (0,1)$ such that
\be\lab{2014-4-22-interpolation-inequality-1}
|u|_{2^*(s_2), {s_2}}\leq |u|_{2^*(s_1), {s_1}}^{\theta}|u|_{2^*(s_3), {s_3}}^{1-\theta}
\ee
for all $u\in L^{2^*(s_1)}(\Omega, \frac{dx}{|x|^{s_1}})\cap L^{2^*(s_3)}(\Omega, \frac{dx}{|x|^{s_3}})$.
\eo
\bp
Define $\varrho=\frac{s_3-s_2}{s_3-s_1}$, then $1-\varrho=\frac{s_2-s_1}{s_3-s_1}$. A direct calculation shows that
\be\lab{2014-5-5-e1}
s_2=\varrho s_1+(1-\varrho)s_3
\ee
and
\be\lab{2014-5-5-e2}
2^*(s_2)=\varrho 2^*(s_1)+(1-\varrho)2^*(s_3).
\ee
It follows from the H\"{o}lder inequality that
\begin{align*}
\int_\Omega \frac{|u|^{2^*(s_2)}}{|x|^{s_2}}dx=&\int_\Omega \big(\frac{|u|^{2^*(s_1)}}{|x|^{s_1}}\big)^\varrho  \big(\frac{|u|^{2^*(s_3)}}{|x|^{s_3}}\big)^{1-\varrho}dx\\
\leq&\Big(\int_\Omega \frac{|u|^{2^*(s_1)}}{|x|^{s_1}}dx\Big)^\varrho \Big(\int_\Omega \frac{|u|^{2^*(s_3)}}{|x|^{s_3}}dx\Big)^{1-\varrho}.
\end{align*}
Let $\theta:=\frac{2^*(s_1)}{2^*(s_2)}\varrho$, then by (\ref{2014-5-5-e2}) again, $1-\theta=\frac{2^*(s_3)}{2^*(s_2)}(1-\varrho)$.
Then we obtain that
$$
|u|_{2^*(s_2), {s_2}}\leq |u|_{2^*(s_1), {s_1}}^{\theta}|u|_{2^*(s_3), {s_3}}^{1-\theta}
$$
for all $u\in L^{2^*(s_1)}(\Omega,\frac{dx}{|x|^{s_1}})\cap L^{2^*(s_3)}(\Omega,\frac{dx}{|x|^{s_3}})$, where
$$\theta=\frac{2^*(s_1)}{2^*(s_2)}\varrho =\frac{(N-s_1)(s_3-s_2)}{(N-s_2)(s_3-s_1)}\in (0,1)$$  has the following
properties.  Firstly, we note that $\theta>0$   since $s_1<s_2<s_3\leq 2<N$. Secondly,
$$ \theta<1 \Leftrightarrow  (N-s_1)(s_3-s_2)<(N-s_2)(s_3-s_1) \Leftrightarrow  (s_2-s_1)(N-s_3)>0.
$$
\ep

Define
\be\lab{2014-5-8-we1}
\vartheta(s_1,s_2):=\frac{N(s_2-s_1)}{s_2(N-s_1)}\;\quad \hbox{for \;$0\leq s_1\leq s_2\leq 2. $}
\ee
\bc\lab{2014-5-5-interpolation-corollary}  Let $\Omega\subset \R^N (N\geq 3)$ be an open set.
Assume  $0\leq s_1<2$.  Then for any $s_2\in [s_1,2]$ and $\theta\in [\vartheta(s_1,s_2), 1]$, there exists $C(\theta)>0$ such that
\be\lab{2014-5-8-e1}
|u|_{2^*(s_1), {s_1}}\leq C(\theta)\|u\|^\theta |u|_{2^*(s_2),  {s_2}}^{1-\theta}
\ee
for all $u\in D_{0}^{1,2}(\Omega)$, where $\|u\|:=\big(\int_\Omega |\nabla u|^2dx\big)^{\frac{1}{2}}$.
\ec
\bp
If $s_2=s_1=s,$ then $ \vartheta(s_1,s_2)=0$ and  (\ref{2014-5-8-e1}) is a direct conclusion of Hardy-Sobolev inequality and the best constant $C(\theta)=\mu_{s}(\Omega)^{-\frac{\theta}{2}},\forall\;\theta\in [0,1]$, where
$\mu_s(\Omega)$ is defined by (\ref{2014-2-27-e2}).
If $s_1=0$, then $\vartheta(s_1,s_2)=1, \theta=1$ and (\ref{2014-5-8-e1}) is just the well-known Sobolev inequality.

Next, we assume that $0<s_1<s_2\leq 2$. We also note that if $\theta=1$, (\ref{2014-5-8-e1}) is  just the well-known Sobolev inequality. Hence,  next we  always assume that $\theta<1$.  Define   $$\tilde{s}:=s_2-\frac{(N-s_2)(s_2-s_1)}{\theta(N-s_1)-(s_2-s_1)}.$$
Note  $\theta\in [\vartheta(s_1,s_2), 1)$, we have  that $0\leq \tilde{s}< s_1<s_2\leq 2$.
Then by Proposition \ref{2014-4-22-interpolation-wl1}, we have
$$ |u|_{2^*(s_1), {s_1}}\leq |u|_{2^*(\tilde{s}), {\tilde{s}}}^{\theta}|u|_{2^*(s_2), {s_2}}^{1-\theta}.$$
Recalling the Hardy-Sobolev inequality, we have
$$|u|_{2^*(\tilde{s}), {\tilde{s}}}\leq \mu_{\tilde{s}}(\Omega)^{-\frac{1}{2}}\|u\|.$$
Hence, there exists a $C(\theta)>0$ such that
$$|u|_{2^*(s_1),  {s_1}}\leq C(\theta)\|u\|^\theta|u|_{2^*(s_2), {s_2}}^{1-\theta}.$$
\ep

Define
\be\lab{2014-5-8-we2}
\varsigma(s_1,s_2):=\frac{(N-s_1)(2-s_2)}{(N-s_2)(2-s_1)}\;\quad \hbox{for   \;$0\leq s_1\leq s_2\leq 2.$}
\ee
\bc\lab{2014-4-22-interpolation-corollary}   Let $\Omega\subset \R^N (N\geq 3)$ be an open set.
Assume  $0<s_2\leq 2$.   Then for any $s_1\in [0,s_2]$ and $\sigma\in [0,\varsigma(s_1,s_2)]$, there exists a  $C(\sigma)>0$ such that
\be\lab{2014-5-8-we3}
|u|_{2^*(s_2),   {s_2}}\leq C(\sigma)\|u\|^{1-\sigma} |u|_{2^*(s_1),  {s_1}}^{\sigma}
\ee
for all $u\in D_{0}^{1,2}(\Omega)$.
\ec
\bp
We only need to consider that case of $s_1<s_2$ and $\sigma>0$.
Define
$$\bar{s}:=s_1+\frac{(N-s_1)(s_2-s_1)}{(N-s_1)-(N-s_2)\sigma}.$$
Recall that  $\sigma\in (0,\varsigma(s_1,s_2)]$,  we have $0\leq s_1<s_2<\bar{s}\leq 2$.
Then by Proposition \ref{2014-4-22-interpolation-wl1}, we have
$$|u|_{2^*(s_2), {s_2}}\leq |u|_{2^*(s_1), {s_1}}^{\sigma}|u|_{2^*(\bar{s}), {\bar{s}}}^{1-\sigma}.$$
Recalling the Hardy-Sobolev inequality, we have
$|u|_{2^*(\bar{s}), {\bar{s}}}\leq \mu_{\bar{s}}(\Omega)^{-\frac{1}{2}}\|u\|. $
Hence, there exists a  $C(\sigma)>0$ such that
$$|u|_{2^*(s_2), {s_2}}\leq C(\sigma)\|u\|^{1-\sigma}|u|_{2^*(s_1), {s_1}}^{\sigma}.$$
\ep

\br\lab{2015-2-15-r1}
 The above Corollary \ref{2014-5-5-interpolation-corollary} and Corollary \ref{2014-4-22-interpolation-corollary} are essentially the well known CKN inequality. However, based  on the Proposition \ref{2014-4-22-interpolation-wl1}, our proofs are very concise.  Moreover, the expressions of \eqref{2014-5-8-e1} and \eqref{2014-5-8-we3} are very convenient in our applications.
\er

\s{Regularity, symmetry and decay estimation}
\renewcommand{\theequation}{3.\arabic{equation}}
\renewcommand{\theremark}{3.\arabic{remark}}
\renewcommand{\thedefinition}{3.\arabic{definition}}
In this section, we will study the regularity, symmetry and decay estimation about the positive solutions.

\bl\lab{2014-7-20-l1}
Assume that $0<s_1\leq s_2<2, 0<u\in D_{0}^{1,2}(\R_+^N)$ and $|u|^{2^*(s_1)-1}/|x|^{s_1}\in L^q(B_1^+)$ for all $1\leq q<q_1$, where $B_1^+:=B_1(0)\cap \R_+^N$. Then
$$|u|^{2^*(s_2)-1}/|x|^{s_2}\in L^q(B_1^+)  \quad \hbox{  for all  }  1\leq q<\frac{N(N+2-2s_1)q_1}{N(N+2-2s_2)+(s_2-s_1)q_1}.$$
Further,  if $|u|^{2^*(s_1)-1}/|x|^{s_1}\in L^q(B_1^+)$ for all $1\leq q<\infty$, then we have $$|u|^{2^*(s_2)-1}/|x|^{s_2}\in L^q(B_1^+)\quad \hbox{ for all }  1\leq q<\frac{N(N+2-2s_1)}{(N+2)(s_2-s_1)}.$$
\el
\bp
When $q<\frac{N(N+2-2s_1)q_1}{N(N+2-2s_2)+(N+2)(s_2-s_1)q_1}$ and $0<s_1\leq s_2<2$, we see that
$$\frac{2^*(s_2)-1}{2^*(s_1)-1}\frac{q}{q_1}<1-\frac{s_2q-\frac{2^*(s_2)-1}{2^*(s_1)-1}s_1q}{N}\leq 1.$$
Then we can take some $\theta\in (0,1)$ such that
$$\frac{2^*(s_2)-1}{2^*(s_1)-1}\frac{q}{q_1}<\theta<1-\frac{s_2q-\frac{2^*(s_2)-1}{2^*(s_1)-1}s_1q}{N}.$$
Let
$$t=\frac{s_2q-\frac{2^*(s_2)-1}{2^*(s_1)-1}qs_1}{1-\theta},\;\quad
\tilde{q}=\frac{1}{\theta}\frac{2^*(s_2)-1}{2^*(s_1)-1}q.$$
 Then by the choice of $\theta$, we have $t<N$ and $\tilde{q}<q_1$. Hence, by the H\"older inequality, we have
\be\lab{2014-11-19-e1}
\int_{B_1^+} \frac{u^{(2^*(s_2)-1)q}}{|x|^{s_2q}}dx\leq \Big(\int_{B_1^+} \frac{u^{(2^*(s_1)-1)\tilde{q}}}{|x|^{s_1\tilde{q}}}dx\Big)^\theta \Big(\int_{B_1^+} \frac{1}{|x|^t}\Big)^{1-\theta}<+\infty.
\ee
It is easy to see that $\displaystyle \frac{N(N+2-2s_1)q_1}{N(N+2-2s_2)+(N+2)(s_2-s_1)q_1}$ is  increasing by $q_1$ and goes to $ \displaystyle \frac{N(N+2-2s_1)}{(N+2)(s_2-s_1)}$ as $q_1\rightarrow \infty$.
\ep

\bl\lab{2014-11-19-l1}
Assume that $0<s_2\leq s_1<2, 0<u\in D_{0}^{1,2}(\R_+^N)$ and $|u|^{2^*(s_1)-1}/|x|^{s_1}\in L^q(B_1^+)$ for all $1\leq q<q_1$, where $B_1^+:=B_1(0)\cap \R_+^N$. Then
$|u|^{2^*(s_2)-1}/|x|^{s_2}\in L^q(B_1^+)$ for all $1\leq q<\frac{2^*(s_1)-1}{2^*(s_2)-1}q_1$.
\el
\bp
For any $1\leq q<\frac{2^*(s_1)-1}{2^*(s_2)-1}q_1$, we set $t=\frac{2^*(s_2)-1}{2^*(s_1)-1}s_1q-s_2q$ and
$\tilde{q}=\frac{2^*(s_2)-1}{2^*(s_1)-1}q$. Then under the assumptions, it is easy to see that $t\geq0$ and $1<\tilde{q}<q_1$. Hence,
\begin{align}\lab{2014-11-19-xe1}
\int_{B_1^+} \frac{u^{(2^*(s_2)-1)q}}{|x|^{s_2q}}dx=\int_{B_1^+} \frac{u^{(2^*(s_1)-1)\tilde{q}}}{|x|^{s_1\tilde{q}}} |x|^tdx
\leq\int_{B_1^+} \frac{u^{(2^*(s_1)-1)\tilde{q}}}{|x|^{s_1\tilde{q}}} dx<+\infty.
\end{align}
\ep

We note that for some subset $\Omega_1$ and some $q\geq 1$ such that  $$|u|^{2^*(s_2)-1}/|x|^{s_2}, |v|^{2^*(s_2)-1}/|x|^{s_2}\in L^q(\Omega_1),$$ then by H\"older inequality, we also have
$\frac{|u|^{t_1}|v|^{t_2}}{|x|^{s_2}}\in L^q(\Omega_1)$ provided $0<t_1,t_2<2^*(s_2)-1$ and $t_1+t_2=2^*(s_2)-1$. Hence, we can obtain the following result:

\bo\lab{2014-11-21-prop1}
Assume $s_1,s_2\in(0,2), \kappa>0, \alpha>1, \beta>1,\alpha+\beta=2^*(s_2)$, then
any positive solution $(u,v)$ of
\be\lab{2014-7-17-e2}
\begin{cases}
-\Delta u-\lambda \frac{|u|^{2^*(s_1)-2}u}{|x|^{s_1}}=\kappa\alpha \frac{1}{|x|^{s_2}}|u|^{\alpha-2}u|v|^\beta\quad &\hbox{in}\;\R_+^N,\\
-\Delta v-\mu \frac{|v|^{2^*(s_1)-2}v}{|x|^{s_1}}=\kappa\beta \frac{1}{|x|^{s_2}}|u|^{\alpha}|v|^{\beta-2}v\quad &\hbox{in}\;\R_+^N,\\
(u,v)\in \mathscr{D}:=D_{0}^{1,2}(\R_+^N)\times D_{0}^{1,2}(\R_+^N),
\end{cases}
\ee
satisfying the following properties:
\begin{itemize}
\item[$(i)$] if $0<\max\{s_1,s_2\}<\frac{N+2}{N}$,  then $u,v\in C^2(\overline{\R_+^N})$;
\item[$(ii)$] if $\max\{s_1,s_2\}=\frac{N+2}{N}$,  then $u, v\in C^{1,\gamma}({\R_+^N})$ for all $0<\gamma<1$;
\item[$(iii)$]if $\max\{s_1,s_2\}>\frac{N+2}{N}$, then $u,v\in C^{1,\gamma}({\R_+^N})$ for all $0<\gamma<\frac{N(2-\max\{s_1,s_2\})}{N-2}$.
\end{itemize}
\eo
\bp

Indeed, it is enough to consider the regularity theorem at $0\in \partial\R_+^N$. By \cite[Lemma B.3]{Struwe.2008}, $u,v$ are locally bounded. Let $B_1^+:=B_1(0)\cap \R_+^N$.
We see that there exists some $C>0$ such that
\begin{align}\lab{2014-7-18-e1}
&|u(x)|^{2^*(s_1)-1}/|x|^{s_1}\leq C|x|^{-s_1}, |v(x)|^{2^*(s_1)-1}/|x|^{s_1}\leq C|x|^{-s_1},\nonumber\\
&|u(x)|^{2^*(s_2)-1}/|x|^{s_2}\leq C|x|^{-s_2}, |v(x)|^{2^*(s_2)-1}/|x|^{s_2}\leq C|x|^{-s_2}\;
\end{align}
for  $x\in B_1^+.$  Hence
$$|u|^{2^*(s_1)-1}/|x|^{s_1},|v|^{2^*(s_1)-1}/|x|^{s_1}\in L^q(B_1^+)\;\hbox{for all}\;1\leq q<\frac{N}{s_1}$$
and
$$\kappa\alpha \frac{1}{|x|^{s_2}}|u|^{\alpha-2}u|v|^\beta, \;\kappa\beta \frac{1}{|x|^{s_2}}|u|^{\alpha}|v|^{\beta-2}v\in  L^q(B_1^+)\;\hbox{for all}\;1\leq q<\frac{N}{s_2}.$$
Set $s_{max}:=\max\{s_1, s_2\}$ and $s_{min}:=\min\{s_1, s_2\}$.
Then we have that $u,v\in W^{2,q}(B_1^+)$ for all $1\leq q<\frac{N}{s_{max}}$.
Denote
$$\tau_u:=\sup\{\tau: \; \sup_{B_1^+}(|u(x)|/|x|^{\tau})<\infty, 0<\tau<1\},$$
$$\tau_v:=\sup\{\tau: \; \sup_{B_1^+}(|v(x)|/|x|^{\tau})<\infty, 0<\tau<1\}$$
and
$$\tau_0:=\min\{\tau_u,\tau_v\}.$$

\noindent{\bf Step 1:  We prove that $\tau_0=1$, i.e., $\tau_u=\tau_v=1$.}

\noindent{\bf Case 1: $s_{max}\leq 1$.}  By the Sobolev embedding,  we have $u, v\in C^\tau(\overline{B_1^+})$ for any $0<\tau<1$. Hence, $\tau_0=1$ in this case.

\noindent{\bf Case 2: $s_{max}>1$.}
For this case,  we have $u,v\in C^\tau(\overline{B_1^+})$ for all $0<\tau<\min\{2-s_{max}, 1\}$.
Then by the definition, we have $ 2-s_{max}\leq \tau_0\leq 1$.
For any $0<\tau<\tau_0$, we have $|u(x)|\leq C|x|^\tau$ and $|v(x)|\leq C|x|^\tau$ for $x\in\overline{B_1^+}$, then for any $x\in B_1^+$, there exists some $C>0$ such that
\begin{align}\lab{2014-7-18-e2}
|u(x)|^{2^*(s_i)-1}/|x|^{s_i}
,|v(x)|^{2^*(s_i)-1}/|x|^{s_i}\leq C|x|^{\big(2^*(s_i)-1\big)\tau-s_i}\;\quad (i=1,2).
\end{align}
Suppose $\tau_0<1$, then by (\ref{2014-7-18-e2}), there must hold $\big(2^*(s_{max})-1\big)\tau_0-s_{max}<0$. Otherwise,
$$ |u|^{2^*(s_{max})-1}/|x|^{s_{max}}\in L^q(B_1^+),\;|v|^{2^*(s_{max})-1}/|x|^{s_{max}}\in L^q(B_1^+)$$
for all $1\leq q<\infty.$
On the other hand,  by Lemma \ref{2014-11-19-l1} and H\"older inequality, it is easy to prove that
 $$\kappa\alpha \frac{1}{|x|^{s_2}}|u|^{\alpha-2}u|v|^\beta,\quad  \;\kappa\beta \frac{1}{|x|^{s_2}}|u|^{\alpha}|v|^{\beta-2}v\in  L^q(B_1^+)\;\hbox{for all}\;1\leq q<\infty.$$
 It follows that $u\in W^{2,q}(B_{1/2}^{+})$ for any $1\leq q<\infty$ and then by  the Sobolev embedding again we have $\tau_0=1$, a contradiction.
 Therefore,  $\big(2^*(s_{max})-1\big)\tau_0-s_{max}<0$ is proved and thus we have
 $$|u|^{2^*(s_{max})-1}/|x|^{s_{max}},  |v|^{2^*(s_{max})-1}/|x|^{s_{max}}\in L^q(B_1^+)$$ for all $\displaystyle 1\leq q<\frac{N}{s_{max}-\big(2^*(s_{max})-1\big)\tau_0}.$

 \noindent{\it Subcase 2.1:}
 If $s_{min}\leq 1$ or $\big(2^*(s_{min})-1\big)\tau_0-s_{min}\geq 0$, we have
 $$
 |u|^{2^*(s_{min})-1}/|x|^{s_{min}}, \;|v|^{2^*(s_{min})-1}/|x|^{s_{min}}\in L^q(B_1^+)$$
for all $1\leq q<N$.
  We claim that $s_{max}-\big(2^*(s_{max})-1\big)\tau_0>1$. If not, we see that $u,v\in W^{2,q}(B_1^+)$ for all $1\leq q<N$, and then by Sobolev embedding, we obtain that $\tau_0=1$, a contradiction.
 Hence, we have $u,v\in W^{2,q}(B_1^+)$ for all $1\leq q<\frac{N}{s_{max}-\big(2^*(s_{max})-1\big)\tau_0}$, and by  the Sobolev embedding again, we have $u,v\in C^\tau(\overline{B_{1/2}^{+}})$ for all $0<\tau<\min\{2-\big[s_{max}-\big(2^*(s_{max})-1\big)\tau_0\big], 1\}$.
 Then by the definition of $\tau_0$, we should have
 $$ 2-\big[s_{max}-\big(2^*(s_{max})-1\big)\tau_0]\leq \tau_0$$
 which implies that
 $$ 2-s_{max}+\big(2^*(s_{max})-2\big)\tau_0\leq 0.$$
 But $s_{max}<2, 2^*(s_{max})>2, \tau_0>0$, a contradiction again.

 \noindent{\it  Subcase 2.2:}
 If $1<s_{min}\leq s_{max}<2$ and $\big(2^*(s_{min})-1\big)\tau_0-s_{min}<0$, by Lemma \ref{2014-11-19-l1} again,
$$
 |u|^{2^*(s_{min})-1}/|x|^{s_{min}}, \;\;\; |v|^{2^*(s_{min})-1}/|x|^{s_{min}}\in L^q(B_1^+)$$
for all $$ 1\leq q<\frac{2^*(s_{max})-1}{2^*(s_{min})-1}\frac{N}{s_{max}-\big(2^*(s_{max})-1\big)\tau_0}.$$
 On the other hand, by the definition of $\tau_0$, we have that
 $$
 |u|^{2^*(s_{min})-1}/|x|^{s_{min}}, \;|v|^{2^*(s_{min})-1}/|x|^{s_{min}}\in L^q(B_1^+)$$
for all $\displaystyle 1\leq q<\frac{N}{s_{min}-\big(2^*(s_{min})-1\big)\tau_0}$.
Thus, $$
 |u|^{2^*(s_{min})-1}/|x|^{s_{min}}, \quad \;|v|^{2^*(s_{min})-1}/|x|^{s_{min}}\in L^q(B_1^+)$$
for all $$1\leq q<\max\Big\{\frac{N}{s_{min}-\big(2^*(s_{min})
-1\big)\tau_0},\frac{2^*(s_{max})-1}{2^*(s_{min})-1}\frac{N}{s_{max}
-\big(2^*(s_{max})-1\big)\tau_0}\Big\}.$$
Noting that
 \begin{align*}
 &\frac{2^*(s_{max})-1}{2^*(s_{min})-1}\frac{N}{s_{max}-\big(2^*(s_{max})-1\big)\tau_0}\\
 & \leq \frac{N}{s_{max}-\big(2^*(s_{max})-1\big)\tau_0}\\
 &\leq \frac{N}{s_{min}-\big(2^*(s_{min})-1\big)\tau_0},
 \end{align*}
we have that $\displaystyle
 |u|^{2^*(s_{min})-1}/|x|^{s_{min}}, \;|v|^{2^*(s_{min})-1}/|x|^{s_{min}}\in L^q(B_1^+)$
for all $1\leq q<\frac{N}{s_{min}-\big(2^*(s_{min})-1\big)\tau_0}$ and it follows that
$$u,v\in W^{2,q}(B_1^+) \;  \hbox{  for  } \;   1\leq q<\frac{N}{s_{max}-\big(2^*(s_{max})-1\big)\tau_0}.$$
Then apply the similar arguments as that in the subcase 2.1, we can deduce a contradiction.
Hence, $\tau_0=1$ is proved and then $\tau_u=\tau_v=1$, i.e., for any $0<\tau<1$,
 \begin{align}\lab{2014-11-21-we1}
|u(x)|^{2^*(s_i)-1}/|x|^{s_i},|v(x)|^{2^*(s_i)-1}/|x|^{s_i}\leq C|x|^{\big(2^*(s_i)-1\big)\tau-s_i},\quad (i=1,2)
\end{align}

\noindent{\bf Step 2: We prove that $u,v\in W^{2,q}(B_{1}^{+})$ for all $$1\leq q<\begin{cases}\infty\;&\hbox{if}\;2^*(s_{max})-1-s_{max}\geq 0\\
 \frac{N}{1+s_{max}-2^*(s_{max})}\;&\hbox{if}\;2^*(s_{max})-1-s_{max}<0\end{cases}.$$}
We divide the proof in two cases.

  \noindent{\it Case 1:\;$2^*(s_{max})-1-s_{max}\geq 0$, i.e., $s_{max}\leq \frac{N+2}{N}$.} By taking $\tau$ close to $1$, we see that
 $$|u|^{2^*(s_i)-1}/|x|^{s_i},|v|^{2^*(s_i)-1}/|x|^{s_i}
 \in L^q(B_1^+)\;\quad (i=1,2)$$
 for all  $\;1<q<\infty.$
 Meanwhile, by  the  H\"older inequality,
 $$\kappa\alpha \frac{1}{|x|^{s_2}}|u|^{\alpha-2}u|v|^\beta,\;\kappa\beta \frac{1}{|x|^{s_2}}|u|^{\alpha}|v|^{\beta-2}v\in  L^q(B_1^+)\;\hbox{for all}\;1\leq q<\infty.$$
 Hence, $u,v\in W^{2,q}(B_{\frac{1}{2}}^{+})$ for all $1\leq q<\infty$.

 \noindent{\it Case 2:\;$2^*(s_{max})-1-s_{max}<0$, i.e., $ \frac{N+2}{N}<s_{max}<2$.}
 In this case,  we have
 $$|u|^{2^*(s_{max})-1}/|x|^{s_{max}},|v|^{2^*(s_{max})-1}/|x|^{s_{max}}\in L^q(B_1^+)\;$$
 for all  $$ 1<q<\frac{N}{1+s_{max}-2^*(s_{max})}.$$
If  $2^*(s_{min})-1-s_{min}\geq 0$, then we see that
 $$|u|^{2^*(s_{min})-1}/|x|^{s_{min}},|v|^{2^*(s_{min})-1}/|x|^{s_{min}}\in L^q(B_1^+)\;\hbox{for all}\;1<q<\infty.$$
Hence, $u,v\in W^{2,q}(B_1^+)$ for all $1\leq q<\frac{N}{1+s_{max}-2^*(s_{max})}$.
 And if $2^*(s_{min})-1-s_{min}<0$,  we must have
 $$|u|^{2^*(s_{min})-1}/|x|^{s_{min}},|v|^{2^*(s_{min})-1}/|x|^{s_{min}}\in L^q(B_1^+)$$
 for all $\;1<q<\frac{N}{1+s_{min}-2^*(s_{min})}.$
 Since $$\frac{N}{1+s_{min}-2^*(s_{min})}\geq \frac{N}{1+s_{max}-2^*(s_{max})},$$ we also obtain that $u,v\in W^{2,q}(B_1^+)$ for all $1\leq q<\frac{N}{1+s_{max}-2^*(s_{max})}$.

\noindent{\bf Step 3:}
 By the Sobolev embedding  theorem,
 $$ u,v\in C^{1,\gamma}(\overline{B_{1/2}^{+}})\;\hbox{for all}\;0<\gamma<1\;\hbox{if}\;s_{max}\leq \frac{N+2}{N}.$$
 In particular, in the case $s_{max}<\frac{N+2}{N}$, there exists $q_0>N$ such that
 \begin{align*}
 &\|u\|_{W^{3,q_0}(B_{1/2}^{+})}\\
 \leq&C\Big(1+\|\frac{u^{2^*(s_1)-2}\nabla u}{|x|^{s_1}}\|_{L^{q_0}(B_1^+)}+\|\frac{u^{2^*(s_1)-1} }{|x|^{s_1+1}}\|_{L^{q_0}(B_1^+)}+\|\frac{u^{\alpha-2}v^\beta\nabla u}{|x|^{s_2}}\|_{L^{q_0}(B_1^+)}\\
 &+\|\frac{u^{\alpha-1}v^{\beta-1}\nabla v}{|x|^{s_2}}\|_{L^{q_0}(B_1^+)}+\|\frac{u^{\alpha-1}v^\beta}{|x|^{s_2+1}}\|_{L^{q_0}(B_1^+)}\Big)
 <\infty.
 \end{align*}
 Thus, we obtain that $u\in C^2(\overline{B_{1/2}^{+}})$. Similarly, we can also prove that $v\in C^2(\overline{B_{1/2}^{+}})$.
 If $s_{max}>\frac{N}{N+2}$, note that $\frac{N}{1+s_{max}-2^*(s_{max})}>N$, by taking $\tau$ close to $1$, we have $u,v\in C^{1,\gamma}(\overline{B_{1/2}^{+}})$ for all $0<\gamma<1-[1+s_{max}-2^*(s_{max})]=\frac{N(2-s_{max})}{N-2}$.
\ep

\bo\lab{2014-11-21-Prop2}
Assume that $s_1,s_2\in (0,2), \kappa>0, \alpha>1, \beta>1,\alpha+\beta=2^*(s_2)$. Let $(u,v)$ be a
positive solution of \eqref{2014-7-17-e2},
then there exists a constant $C$ such that $$|u(x)|, |v(x)|\leq C (1+|x|^{-(N-1)});\quad |\nabla u(x)|, |\nabla v(x)|\leq C|x|^{-N}.$$
\eo
\bp
Recalling the Kelvin transformation:
\be\lab{2014-11-19-xe2}
u^*(y):=|y|^{-(N-2)}u(\frac{y}{|y|^2}),\;\; v^*(y):=|y|^{-(N-2)}v(\frac{y}{|y|^2}).
\ee
It is well known that
\be\lab{2014-11-21-we2}
\Delta u^*(y)=\frac{1}{|y|^{N+2}}(\Delta u)(\frac{y}{|y|^2})\;\hbox{and}\;\Delta v^*(y)=\frac{1}{|y|^{N+2}}(\Delta v)(\frac{y}{|y|^2}).
\ee
Hence, a direct computation shows that $\big(u^*, v^*\big)$ is also a positive solution to the same equation.

By Proposition \ref{2014-11-21-prop1}, we see that $u^*, v^*\in C^{1,\gamma}(\overline{\R_+^N})$ for some $\gamma>0$. Then $|u^*(y)|, |v^*(y)|\leq C|y|$ for $y\in B_1^+$. Going back to $(u,v)$, we see that $|u(y)|, |v(y)|\leq C|y|^{-(N-1)}$ for $y\in \R_+^N$. Finally, it is standard to apply the gradient estimate, we obtain that $|\nabla u(y)|, |\nabla v(y)|\leq C|y|^{-N}$ for $y\in \R_+^N$.
\ep

\br\lab{zouwm-0611-1}   Checking the proofs of Lemmas  \ref{2014-7-20-l1}-\ref{2014-11-19-l1} and Propositions
\ref{2014-11-21-prop1}-\ref{2014-11-21-Prop2}, their conclusions  are valid for general cone $\Omega$. A little difference is that when $\Omega=\R^N$, the decay estimation is $$|u(x)|, |v(x)|\leq C (1+|x|^{-N});\quad |\nabla u(x)|, |\nabla v(x)|\leq C|x|^{-N-1}.$$ \er

\bo\lab{2014-11-22-Prop1}
Assume that $s_1,s_2\in (0,2), \kappa>0, \alpha>1, \beta>1,\alpha+\beta=2^*(s_2)$. Let $(u,v)$ be a
positive solution of \eqref{2014-7-17-e2}.
Then we have that $u\circ \sigma=u, v\circ \sigma=v$ for all isometry of $\R^N$ such that $\sigma(\R_+^N)=\R_+^N$. In particular, $\big(u(x',x_N),v(x',x_N)\big)$ is axially symmetric with respect to the $x_N-$axis, i.e., $u(x',x_N)=u(|x'|,x_N)$ and $v(x',x_N)=v(|x'|,x_N)$.
\eo
\bp
We prove the result by the well-known method of moving planes. Denote by $\overrightarrow{e_N}$ the $N^{th}$ vector of the canonical basis of $\R^N$ and consider the open ball $D:=B_{\frac{1}{2}}(\frac{1}{2}\overrightarrow{e_N})$. Set
\be\lab{2014-11-22-xe1}
\varphi(x):=|x|^{2-N}u(-\overrightarrow{e_N}+\frac{x}{|x|^2}),
\psi(x):=|x|^{2-N}v(-\overrightarrow{e_N}+\frac{x}{|x|^2})
\ee
for all $x\in \bar{D}\backslash\{0\}$ and $\varphi(0)=\psi(0)=0$.
By Proposition \ref{2014-11-21-prop1}, $\varphi(x), \psi(x)\in C^2(D)\cap C^1(\bar{D}\backslash\{0\})$. We note that it is easy to see that $\varphi(x)>0, \psi(x)>0$ in $D$ and $\varphi(x)=\psi(x)=0$ on $\partial D\backslash \{0\}$. On the other hand, by Proposition \ref{2014-11-21-Prop2}, there exists $C>0$ such that
\be\lab{2014-11-22-xe2}
\varphi(x)\leq C|x|, \psi(x)\leq C|x|\;\hbox{for all}\;x\in \bar{D}\backslash \{0\}.
\ee
Since $\varphi(0)=\psi(0)=0$, we have that $\varphi(x),\psi(x)\in C^0(\bar{D})$. By a direct computation, $\big(\varphi(x),\psi(x)\big)$ satisfies the following eqaution
\be\lab{2014-11-22-xe3}
\begin{cases}
-\Delta \varphi-\lambda \frac{\varphi^{2^*(s_1)-1}}{\big|x-|x|^2\overrightarrow{e_N}\big|^{s_1}}=\kappa \alpha \frac{\varphi^{\alpha-1}\psi^{\beta}}{\big|x-|x|^2\overrightarrow{e_N}\big|^{s_2}}\\
-\Delta \psi-\mu \frac{\psi^{2^*(s_1)-1}}{\big|x-|x|^2\overrightarrow{e_N}\big|^{s_1}}=\kappa \beta \frac{\varphi^{\alpha}\psi^{\beta-1}}{\big|x-|x|^2\overrightarrow{e_N}\big|^{s_2}}
\end{cases}\;\hbox{in}\;D.
\ee
Noting that
\be\lab{2014-11-22-xe4}
\big|x-|x|^2\overrightarrow{e_N}\big|=\big|x\big| \big|x-\overrightarrow{e_N}\big|,
\ee
we have
\be\lab{2014-11-22-xe5}
\begin{cases}
-\Delta \varphi-\lambda \frac{\varphi^{2^*(s_1)-1}}{\big|x\big|^{s_1} \big|x-\overrightarrow{e_N}\big|^{s_1}}=\kappa \alpha \frac{\varphi^{\alpha-1}\psi^{\beta}}{\big|x\big|^{s_2} \big|x-\overrightarrow{e_N}\big|^{s_2}}\\
-\Delta \psi-\mu \frac{\psi^{2^*(s_1)-1}}{\big|x\big|^{s_1} \big|x-\overrightarrow{e_N}\big|^{s_1}}=\kappa \beta \frac{\varphi^{\alpha}\psi^{\beta-1}}{\big|x\big|^{s_1} \big|x-\overrightarrow{e_N}\big|^{s_1}}
\end{cases}\;\hbox{in}\;D.
\ee
Since $\overrightarrow{e_N}\in \partial D\backslash \{0\}$ and $\varphi(x),\psi(x)\in C^1(\bar{D}\backslash \{0\})\cap C^0(\bar{D})$, there exists $C>0$ such that
\be\lab{2014-11-22-xe6}
\varphi(x)\leq C|x-\overrightarrow{e_N}|, \psi(x)\leq C|x-\overrightarrow{e_N}|\;\hbox{for all}\;x\in \bar{D}.
\ee
Noting that $2^*(s_i)-1-s_i>-N$ for $i=1,2$, then by \eqref{2014-11-22-xe2}, \eqref{2014-11-22-xe5}, \eqref{2014-11-22-xe6} and the standard elliptic theory, we obtain that $\varphi(x),\psi(x)\in C^1(\bar{D})$. By $\varphi(x)>0,\psi(x)>0$ in $D$, we obtain that $\frac{\partial \varphi}{\partial \nu}<0, \frac{\partial \psi}{\partial \nu}<0$ on $\partial D$, where $\nu$ denotes the outward unit normal to $D$ at $x\in \partial D$.

For any $\eta\geq 0$ and any $x=(x_1,x')\in \R^N$, where $x'=(x_2,\cdots,x_N)\in \R^{N-1}$, we let
\be\lab{2014-11-22-xe7}
x_\eta=(2\eta-x_1, x')\;\hbox{and}\;D_\eta:=\{x\in D\big| x_\eta\in D\}.
\ee
We say that $(P_\eta)$ holds iff
$$D_\eta\neq \emptyset\;\hbox{and}\;\varphi(x)\geq \varphi(x_\eta), \psi(x)\geq \psi(x_\eta)\;\hbox{for all}\;x\in D_\eta\;\hbox{such that}\;x_1\leq \eta.$$

\noindent{\bf Step 1:} We shall prove that $(P_\eta)$ holds if $\eta<\frac{1}{2}$ and close to $\frac{1}{2}$ sufficiently.

\noindent
Indeed, it is easily to follow the Hopf's Lemma (see the arguments above) that there exists $\varepsilon_0>0$ such that $(P_\eta)$ holds for $\eta\in (\frac{1}{2}-\varepsilon_0, \frac{1}{2})$.
Now, we let
\be\lab{2014-11-22-xe8}
\sigma:=\min\{\eta\geq 0; (P_\delta)\;\hbox{holds for all}\;\delta\in (\eta,\frac{1}{2})\}.
\ee

\noindent{\bf Step 2:} We shall prove that $\sigma=0$.

\noindent
We prove it by way of negation. Assume that $\sigma>0$, then we see that $D_\sigma\neq\emptyset$ and that $(P_\sigma)$ holds. Now, we set
\be\lab{2014-11-22-xe9}
\hat{\varphi}(x):=\varphi(x)-\varphi(x_\sigma)\;\hbox{and}\;\hat{\psi}(x):=\psi(x)-\psi(x_\sigma).
\ee
Then we have
\begin{align}\lab{2014-11-22-xe10}
-\Delta \hat{\varphi}=&[-\Delta \varphi(x)]-[-\Delta \varphi(x_\sigma)]\\
=&\lambda \frac{\big(\varphi(x)\big)^{2^*(s_1)-1}}{|x-|x|^2\overrightarrow{e_N}|^{s_1}}
+\kappa \alpha \frac{\big(\varphi(x)\big)^{\alpha-1}\big(\psi(x)\big)^\beta}{|x-|x|^2\overrightarrow{e_N}|^{s_2}}\nonumber\\
&-\lambda \frac{\big(\varphi(x_\sigma)\big)^{2^*(s_1)-1}}{|x_\sigma-|x_\sigma|^2\overrightarrow{e_N}|^{s_1}}
-\kappa \alpha \frac{\big(\varphi(x_\sigma)\big)^{\alpha-1}\big(\psi(x_\sigma)\big)^\beta}{|x_\sigma-|x_\sigma|^2\overrightarrow{e_N}|^{s_2}}\nonumber\\
\geq&\lambda \big(\varphi(x_\sigma)\big)^{2^*(s_1)-1}\big[\frac{1}{|x-|x|^2\overrightarrow{e_N}|^{s_1}}-\frac{1}{|x_\sigma-|x_\sigma|^2\overrightarrow{e_N}|^{s_1}}\big]\nonumber\\
&+\kappa \alpha \big(\varphi(x_\sigma)\big)^{\alpha-1}\big(\psi(x)\big)^\beta\big[\frac{1}{|x-|x|^2\overrightarrow{e_N}|^{s_2}}-\frac{1}{|x_\sigma-|x_\sigma|^2\overrightarrow{e_N}|^{s_2}}\big]\nonumber
\end{align}
for all $x\in D_\sigma\cap \{x_1<\sigma\}$.
Noting that
\begin{align}\lab{2014-11-22-xe11}
\big|x_\sigma+|x_\sigma|^2\overrightarrow{e_N}\big|^2-\big|x+|x|^2\overrightarrow{e_N}\big|^2
=4\sigma(\sigma-x_1)(1+|x_\sigma|^2+|x|^2+2x_N),
\end{align}
we obtain that
\be\lab{2014-11-22-xe12}
-\Delta \hat{\varphi}(x)>0\;\hbox{for all}\;x\in D_\sigma\cap \{x_1<\sigma\}.
\ee
Similarly, we also have
\be\lab{2014-11-22-xe13}
-\Delta \hat{\psi}(x)>0\;\hbox{for all}\;x\in D_\sigma\cap \{x_1<\sigma\}.
\ee
Then by the Hopf's Lemma and the strong comparison principle, we have
\be\lab{2014-11-22-xe14}
\hat{\varphi},\hat{\psi}>0\;\hbox{in}\;D_\sigma\cap \{x_1<\sigma\}\;\hbox{and}\;\frac{\partial \hat{\varphi}}{\partial \nu}, \frac{\partial \hat{\psi}}{\partial \nu}<0\;\hbox{on}\;D_\sigma\cap \{x_1=\sigma\}.
\ee
Here we use the assumption $\sigma>0$. By definition, there exists a subsequence $\{\sigma_i\}_{i\in \NN}\subset \R^+$ and a sequence $\{x^i\}_{i\in \NN}\subset D$ such that $\displaystyle \sigma_i<\sigma, x^i\in D_{\sigma_i}, (x^i)_1<\sigma_i, \lim_{i\rightarrow \infty}\sigma_i=\sigma$ and
\be\lab{2014-11-22-xe14}
\varphi(x^i)<\varphi((x^i)_{\sigma_i})\;\hbox{or}\;\psi(x^i)<\psi((x^i)_{\sigma_i}).
\ee
Up to a subsequence, we may assume that $\varphi(x^i)<\varphi((x^i)_{\sigma_i})$ without loss of generality. Since $\{x^i\}_{i\in \NN}$ is bounded, going to a subsequence again, we assume that $\displaystyle\lim_{i\rightarrow \infty}x_i=x\in \overline{D_\sigma}\cap \{x_1\leq \sigma\}$ due to the choice of $\{x^i\}$. Then we have $\varphi(x)\leq \varphi(x_\sigma)$, i.e., $\hat{\varphi}(x)\leq 0$. Combining with \eqref{2014-11-22-xe14}, we obtain that $\hat{\varphi}(x)=0$ and then $x\in \partial\big(D_\sigma\cap \{x_1<\sigma\}\big)$.

\noindent
{\bf Case 1:} If $x\in \partial D$, then $\varphi(x)=0$. It follows that $\varphi(x_\sigma)=0$. Since $x_\sigma\in D$ and $\varphi>0$ in $D$, we also have $x_\sigma\in\partial D$. We say that $x=x_\sigma$. If not, $x$ and $x_\sigma$ are symmetric with respect to the hyperplane $x_1=\sigma$. This is impossible since that $D$ is a ball, $\sigma>0$ and $x, x_\sigma\in \partial D$. Now recalling that $\varphi\in C^1$, by the mean value theorem, there exists a sequence $\tau_i\in \big((x^i)_1, 2\sigma_i-(x^i)_1\big)$ such that
\be\lab{2014-11-22-xe15}
\varphi(x^i)-\varphi((x^i)_{\sigma_i})=2\partial_1\varphi(\tau_i, (x')^i) \;\big((x^i)_1-\sigma_i\big).
\ee
Using the facts $(x^i)_1<\sigma_i$ and $\varphi(x^i)<\varphi\big((x^i)_{\sigma_i}\big)$, we let $i$ go to infinity and then obtain that
\be\lab{2014-11-22-xe16}
\partial_1\varphi(x)\geq 0.
\ee
On the other hand,
\begin{align}\lab{2014-11-22-xe17}
\partial_1\varphi(x)=\frac{\partial \varphi}{\partial \nu} \langle \nu(x), \overrightarrow{e_1}\rangle
=\frac{\sigma}{|x-\frac{1}{2}\overrightarrow{e_N}|}\frac{\partial\varphi}{\partial\nu}<0,
\end{align}
a contradiction. Here we use the assumption $\sigma>0$ again and the fact $x_1=\sigma$ since $x=x_\sigma$.

\noindent
{\bf Case 2:} If $x\in D$, since $\varphi(x)=\varphi(x_\sigma)>0$, we have $x_\sigma\in D$. Since $x\in \partial\big(D_\sigma\cap \{x_1<\sigma\}\big)$, we have $x\in D\cap \{x_1=\sigma\}$. Then apply the similar arguments in Case 1, we can also obtain that $\partial_1\varphi(x)\geq 0$. On the other hand, by \eqref{2014-11-22-xe14}, we have $2\partial_1\varphi(x)=\partial_1\hat{\varphi}(x)<0$, also a contradiction.
Hence, $\sigma=0$ is proved.

\noindent{\bf Step 3:} By Step 2, $\sigma=0$. Hence, we have
\be\lab{2014-11-22-xe18}
\varphi(x_1,x')\geq \varphi(-x_1,x'), \psi(x_1,x')\geq \psi(-x_1,x')\;\hbox{for all}\;x\in D_0=D.
\ee
Apply the same argument one can obtain the reverse inequality. Thus,
\be\lab{2014-11-22-xe18}
\varphi(x_1,x')= \varphi(-x_1,x'), \psi(x_1,x')= \psi(-x_1,x')\;\hbox{for all}\;x\in D.
\ee
Hence, $\varphi$ and $\psi$ are symmetric with respect to the hyperplane $\{x_1=0\}$. Noting that the arguments above are valid for any hyperplane containing $\overrightarrow{e_N}$. By going back to $(u,v)$, we obtain the conclusions of this proposition.
\ep

\br\lab{2014-11-27-open}{\bf (Open problem)}  The Proposition \ref{2014-11-22-Prop1} is established  under the  assumption
that $\Omega=\R_+^N$. But so far  we do not  know whether the  Proposition \ref{2014-11-22-Prop1} is true for general cone
$\Omega$.  It remains an open problem.
\er

\s{Nehari manifold}
\renewcommand{\theequation}{4.\arabic{equation}}
\renewcommand{\theremark}{4.\arabic{remark}}
\renewcommand{\thedefinition}{4.\arabic{definition}}
In this section, we  study the  Nehari manifold  corresponding to the following equation:
\be\lab{2014-7-17-e2+0}
\begin{cases}
-\Delta u-\lambda \frac{|u|^{2^*(s_1)-2}u}{|x|^{s_1}}=\kappa\alpha \frac{1}{|x|^{s_2}}|u|^{\alpha-2}u|v|^\beta\quad &\hbox{in}\;\R_+^N,\\
-\Delta v-\mu \frac{|v|^{2^*(s_1)-2}v}{|x|^{s_1}}=\kappa\beta \frac{1}{|x|^{s_2}}|u|^{\alpha}|v|^{\beta-2}v\quad &\hbox{in}\;\R_+^N,\\
(u,v)\in \mathscr{D}:=D_{0}^{1,2}(\R_+^N)\times D_{0}^{1,2}(\R_+^N),
\end{cases}
\ee
For $(u,v)\in \mathscr{D}$, we define the norm
$$\|(u,v)\|_{\mathscr{D}}=\big(\|u\|^{2}+\|v\|^{2}\big)^{\frac{1}{2}},$$
where
$\|u\|:=\big(\int_{\R_+^N} |\nabla u|^2dx\big)^{\frac{1}{2}}\;\hbox{for}\;u\in D_{0}^{1,2}(\R_+^N)$.
A pair of function $(u,v)$ is said to be a weak solution of \eqref{2014-7-17-e2+0} if and only if
\begin{align*}
&\int_{\R_+^N}\nabla u\nabla \varphi_1+\nabla v\nabla\varphi_2 dx
-\lambda \int_{\R_+^N}\frac{|u|^{2^*(s_1)-2}u\varphi_1}{|x|^{s_1}}dx
-\mu \int_{\R_+^N}\frac{|v|^{2^*(s_1)-2}v\varphi_2}{|x|^{s_1}}dx\\
&-\kappa\alpha\int_{\R_+^N} \frac{|u|^{\alpha-2}u|v|^\beta\varphi_1}{|x|^{s_2}}dx
-\kappa\beta\int_{\R_+^N}\frac{|u|^{\alpha}|v|^{\beta-2}v\varphi_2}{|x|^{s_2}}dx
=0\;\hbox{for all}\;(\varphi_1,\varphi_2)\in \mathscr{D}.
\end{align*}
The corresponding energy functional of problem \eqref{2014-7-17-e2+0} is defined as
$\displaystyle
\Phi(u,v)=\frac{1}{2}a(u,v)-\frac{1}{2^*(s_1)}b(u,v)
-\kappa c(u,v)
$
for all $(u,v)\in\mathscr{D}$, where
\be\lab{2014-11-23-xe1}
\begin{cases}
a(u,v):=\|(u,v)\|_{\mathscr{D}}^{2},\\
b(u,v):=\lambda \int_{\R_+^N}\frac{|u|^{2^*(s_1)}}{|x|^{s_1}}dx+\mu \int_{\R_+^N}\frac{|v|^{2^*(s_1)}}{|x|^{s_1}}dx,\\
c(u,v):=\int_{\R_+N}\frac{|u|^\alpha |v|^\beta}{|x|^{s_2}}dx.
\end{cases}
\ee
We consider the corresponding Nehari manifold
$$ \mathcal{N}:=\left\{(u,v)\in \mathscr{D}\backslash (0,0)| J(u,v)=0\right\}$$
where
$$
J(u,v):=\langle \Phi'(u,v), (u,v)\rangle
=a(u,v)-b(u,v)-\kappa(\alpha+\beta)c(u,v)
$$
and $\Phi'(u,v)$ denotes the Fr\'{e}chet derivative of $\Phi$ at $(u, v)$ and $\langle \cdot,\cdot \rangle$ is the duality product between $\mathscr{D}$ and its dual space $\mathscr{D}^*$.

\bl\lab{2014-11-23-xl1}
Assume $s_1,s_2\in (0,2), \lambda, \mu\in (0,+\infty), \alpha>1,\beta>1$ and $ \alpha+\beta= 2^*(s_2)$.
Then for any $(u,v)\in \mathscr{D}\backslash \{(0,0)\}$, there exists a unique $t=t_{(u,v)}>0$ such that $t(u,v)=(tu,tv)\in \mathcal{N}$ if one of the following assumptions is satisfied:
\begin{itemize}
\item[(i)]$\kappa>0$.
\item[(ii)] $\kappa<0$ and $s_2>s_1$.
\item[(iii)] $s_2=s_1$  and $ \kappa<0$ with $ |\kappa|$ small enough.
\end{itemize}
Moreover, $\mathcal{N}$ is closed and bounded away from $0$.
\el
\bp
For any $(u,v)\in \mathscr{D}$, we use the notations defined by \eqref{2014-11-23-xe1} and we will write them as $a, b, c$ for simplicity. It is easy to see that for any $(u,v)\neq (0,0)$, we have $a>0,b>0,c\geq 0$.
For any $(u,v)\in \mathscr{D}\backslash \{(0,0)\}$ and $t>0,$ we have
\begin{align}\lab{2014-1-4-e1}
\Phi(tu,tv)=&\frac{1}{2}at^2
-\frac{1}{2^*(s_1)}bt^{2^*(s_1)}-\kappa ct^{2^*(s_2)}.
\end{align}
Denote $\frac{d\Phi(tu,tv)}{dt}:=-tg(t)$, where
$$
g(t)=bt^{2^*(s_1)-2}+\kappa 2^*(s_2)c t^{2^*(s_2)-2}
-a.
$$
For the cases of $(i)$ and $(ii)$, it is easy to see that $g(+\infty)=+\infty$ and $g(0)=-a<0$. Also for the case $(iii)$, by the Young inequality, one can prove that there exists some $C>0$ such that
$c(u,v)\leq  C b(u,v)$ for all $(u,v)\in \mathscr{D}$. Thus, for the case of $s_2=s_1$, if $\kappa<0$ with $|\kappa|$ small enough, we obtain that
\be\lab{2014-5-9-e1}
b(u,v)+2^*(s_2)\kappa c(u,v)>0\;\hbox{ for all}\; (u,v)\in\mathscr{D}\backslash \{(0,0)\}.
 \ee
 Hence, we also have $g(+\infty)=+\infty$ and $g(0)=-a<0$.

 Thus, we obtain that there exists some $t>0$ such that $g(t)=0$ due to the continuity of $g(t)$. It follows that $tu\in \mathcal{N}$.
By the Hardy-Sobolev inequality and the Young inequality, there exists some $C>0$ such that
$$b(u,v)\leq C \|(u,v)\|_{\mathscr{D}}^{2^*(s_1)},\quad  c(u,v)\leq C\|(u,v)\|_{\mathscr{D}}^{2^*(s_2)}.$$
Let $(u,v)\in \mathcal{N}$, by $2^*(s_i)>2, i=1,2$, we have
$$
a=b+\kappa(2^*(s_2))c
\leq C\Big(a^{\frac{2^*(s_1)}{2}}+a^{\frac{2^*(s_2)}{2}}\Big),
$$
which implies that there exists some $\delta_0>0$ such that
\be\lab{2014-4-9-xe2}
\|(u,v)\|_{\mathscr{D}}=a^{\frac{1}{2}}\geq \delta_0\;\hbox{for all}\;(u,v)\in \mathcal{N}.
\ee
Thus, $\mathcal{N}$ is bounded away from $(0,0)$ and obviously,  $\mathcal{N}$ is closed.

For any $(u,v)\neq (0,0)$, let $t_0:=\inf\{t|g(t)=0, t>0\}$. Then we see that $t_0>0$ and $g(t_0)=0$. Without loss of generality, we may assume that $t_0=1$, that is, $g(t)<0$ for $0<t<1$ and $g(1)=0=b+\kappa 2^*(s_2)c-a$.
We note that
$$
g'(t)=\big(2^*(s_1)-2\big)bt^{2^*(s_1)-3}+\kappa 2^*(s_2)\big(2^*(s_2)-2\big)ct^{2^*(s_2)-3}.
$$

\noindent{$(i)$}  If $\kappa>0$, then $g'(t)>0$ for all $t>0$.

\noindent{$(ii)$} If $\kappa<0, s_2>s_1$, recalling that $0=b+\kappa 2^*(s_2)c-a$, we have
\begin{align*}
g'(t)\equiv&\big(2^*(s_1)-2\big)bt^{2^*(s_1)-3}+\kappa 2^*(s_2)\big(2^*(s_2)-2\big)ct^{2^*(s_2)-3}\\
=&:h(t)t^{2^*(s_2)-3},
\end{align*}
with
$\displaystyle h(t):=\big(2^*(s_1)-2\big)bt^{2^*(s_1)-2^*(s_2)}+\kappa 2^*(s_2)\big(2^*(s_2)-2\big)c.$
When $t>1$, we have
\begin{align*}
h(t)>&\big(2^*(s_1)-2\big)b+\kappa 2^*(s_2)\big(2^*(s_2)-2\big)c\\
=&\big(2^*(s_1)-2\big)a-\kappa 2^*(s_2)(2^*(s_1)-2^*(s_2))c>0.
\end{align*}
Hence,  $g'(t)>0$ for all $t>1$.

\noindent{$(iii)$} If $\kappa<0, s_2=s_1$, similar to the arguments as case $(ii)$ above, we  know $$ h(t)=\big(2^*(s_1)-2\big)\big(b+2^*(s_1)\kappa c\big)>0$$
when $\kappa$ is small enough by (\ref{2014-5-9-e1}).

The arguments above imply that $g(t)>0$ for $t>1$. Hence, $t=1$ is the unique solution of $g(t)=0$. It follows  that for any $(u,v)\neq (0,0)$, there exists a unique $t_{(u,v)}>0$ such that
$t_{(u,v)}(u,v)\in \mathcal{N}$
and
$$ \Phi(t_{(u,v)} u, t_{(u,v)}v)=\underset{t>0}{\max}\;\Phi(tu, tv).$$
\ep

\bl\lab{2014-11-23-xl2}
Under the assumptions of Lemma \ref{2014-11-23-xl1},   any $(PS)_m$ sequence of $\Phi(u,v)$ i.e.,
$$\begin{cases}\Phi(u_n, v_n)\rightarrow m\\ \Phi'(u_n, v_n)\rightarrow 0\;\hbox{in}\;\mathscr{D}^* \end{cases}$$
is bounded in $\mathscr{D}$.
\el
\bp
Let $\{(u_n,v_n)\}\subset \mathscr{D}$ be a $(PS)_m$ sequence of $\Phi(u,v)$. We  tend to use the previous marks $a,b,c$ and denote $a(u_n,v_n),b(u_n,v_n),c(u_n,v_n)$  by $a_n,b_n,c_n$ for the simplicity. Then we have
\be\lab{2014-4-10-e2}
\Phi(u_n,v_n)=\frac{1}{2}a_n-\frac{1}{2^*(s_1)}b_n
-\kappa c_n=m+o(1)
\ee
and
\be\lab{2014-4-10-e3}
J(u_n,v_n)=a_n-b_n-\kappa(\alpha+\beta)c_n=o(1)\|(u_n,v_n)\|_{\mathscr{D}}.
\ee

\noindent{(1)} If $\kappa>0$,
for the case of $s_2\leq s_1$, we have
\begin{align*}
&m+o(1)\big(1+\|(u_n,v_n)\|_{\mathscr{D}}\big)
 = \Phi(u_n,v_n)-\frac{1}{2^*(s_1)}J(u_n,v_n)\\
& = \big(\frac{1}{2}-\frac{1}{2^*(s_1)}\big)a_n+\big(\frac{2^*(s_2)}{2^*(s_1)}-1\big)c_n \geq\big(\frac{1}{2}-\frac{1}{2^*(s_1)}\big)\|(u_n,v_n)\|_{\mathscr{D}}^{2}
\end{align*}
and for the case of $s_2>s_1$, we have
\begin{align*}
&  m+o(1)\big(1+\|(u_n,v_n)\|_{\mathscr{D}}\big) = \Phi(u_n,v_n)-\frac{1}{2^*(s_2)}J(u_n,v_n)\\
& = \big(\frac{1}{2}-\frac{1}{2^*(s_2)}\big)a_n+\big(\frac{1}{2^*(s_2)}-\frac{1}{2^*(s_1)}\big)b_n
\geq \big(\frac{1}{2}-\frac{1}{2^*(s_2)}\big)\|(u_n,v_n)\|_{\mathscr{D}}^{2}.
\end{align*}

\noindent{(2)} If $\kappa<0, s_2\geq s_1$, similarly  we obtain that
$$ m+o(1)\big(1+\|(u_n,v_n)\|_{\mathscr{D}}\big)\geq \big(\frac{1}{2}-\frac{1}{2^*(s_2)}\big)\|(u_n,v_n)\|_{\mathscr{D}}^{2}.$$
Based on the above arguments, we can see that $\{(u_n,v_n)\}$ is bounded in $\mathscr{D}$.
\ep

\bl\lab{2014-11-23-xl3}
Under the assumptions of Lemma \ref{2014-11-23-xl1}, let $\{(u_n,v_n)\}\subset\mathcal{N}$ be a $(PS)_c$  sequence for $\Phi\big|_{\mathcal{N}}$, i.e., $\Phi(u_n,v_n)\rightarrow c$ and $\Phi'\big|_{\mathcal{N}}(u_n,v_n)\rightarrow 0$ in $\mathscr{D}^{*}$. Then $\{(u_n,v_n)\}$ is also a $(PS)_c$ sequence for $\Phi$.
\el
\bp
For any $(u,v)\in \mathcal{N}$,
we will follow the previous marks $a,b,c$ defined by \eqref{2014-11-23-xe1}.  Then we have
$\displaystyle a-b-\kappa 2^*(s_2)c=0$
and
$\displaystyle\langle J'(u,v), (u,v)\rangle =2a-2^*(s_1)b-\kappa\big(2^*(s_2)\big)^2c.$

\noindent{(1)} If $\kappa>0$,
\begin{align*}
\langle J'(u,v), (u,v)\rangle=&2\big[b+\kappa 2^*(s_2)c\big]-2^*(s_1)b-k\big(2^*(s_2)\big)^2c\\
=&\big[2-2^*(s_1)\big]b+\big(2-2^*(s_2)\big) 2^*(s_2)\kappa c\\
<&\max\{2-2^*(s_1),2-2^*(s_2)\}\big[b+ 2^*(s_2)\kappa c\big]\\
=&\max\{2-2^*(s_1),2-2^*(s_2)\}a.
\end{align*}

\noindent{(2)} If $\kappa<0, s_2\geq s_1$,
\begin{align*}
\langle J'(u,v), (u,v)\rangle=&2a-2^*(s_1)\big[a-\kappa(\alpha+\beta)c\big]-\kappa(\alpha+\beta)^2c\\
=&\big[2-2^*(s_1)\big]a+\kappa\big[2^*(s_1)-\alpha-\beta\big](\alpha+\beta)c\\
\leq&\big[2-2^*(s_1)\big]a.
\end{align*}
Hence, by (\ref{2014-4-9-xe2}), we obtain that
\be\lab{2014-4-22-xe1}
\langle J'(u,v), (u,v)\rangle\leq \max\{2-2^*(s_1),2-\alpha-\beta\}\delta_0^2<0\;\hbox{for all}\;(u,v)\in \mathcal{N},
\ee
where $\delta_0$ is given by (\ref{2014-4-9-xe2}).
By the similar arguments as in Lemma \ref{2014-11-23-xl2}, we can prove that $\{(u_n,v_n)\}$ is bounded in $\mathscr{D}$. Let $\{t_n\}\subset \R$ be a sequence of multipliers satisfying
$$\Phi'(u_n,v_n)=\Phi'|_{\mathcal{N}}(u_n,v_n)+t_n J'(u_n,v_n).$$
Testing by $(u_n,v_n)$, we obtain that
$$ t_n \langle J'(u_n,v_n), (u_n,v_n)\rangle\rightarrow 0.$$
Recalling (\ref{2014-4-22-xe1}), we obtain $t_n\rightarrow 0$.
We can also have that $J'(u_n,v_n)$ is bounded due to the boundedness of $(u_n,v_n)$.
Hence, it follows that
$\Phi'(u_n,v_n)\rightarrow 0\;\hbox{in}\;\mathscr{D}^*.$
\ep

\vskip 0.26in
Define
\be\lab{2014-4-11-e1}
c_0:=\inf_{(u,v)\in \mathcal{N}}\Phi(u,v)
\ee
and
$\displaystyle \eta:= \frac{1}{2}-\frac{1}{2^*(s_{max})},$
where $s_{max}:=\max\{s_1,s_2\}$.
From the arguments in the proof of Lemma \ref{2014-11-23-xl2}, we obtain that
\be\lab{buzze1}
c_0\geq \eta \|(u,v)\|_{\mathscr{D}}^{2}.
\ee
Combined with Lemma \ref{2014-11-23-xl1}, we have
\be\lab{bane-1}
c_0\geq \eta \delta_0^2,
\ee
where $\delta_0$ is given by (\ref{2014-4-9-xe2}).
If $m_0$ is achieved by some $(u,v)\in \mathcal{N}$, then $(u,v)$ is a ground state solution of \eqref{2014-7-17-e2}.

\s{Nonexistence of nontrivial ground state solution}
\renewcommand{\theequation}{5.\arabic{equation}}
\renewcommand{\theremark}{5.\arabic{remark}}
\renewcommand{\thedefinition}{5.\arabic{definition}}
In this section, we continue to study the equation (\ref{2014-7-17-e2+0}).

\bd In the sequel, we call $(u, v)$   nontrivial iff $u\not=0$ and $v\not=0$,  and call $(u, v)$ semi-trivial iff
either $u=0$ or $v=0$ but not all zero.
\ed
 We obtain the nonexistence of nontrivial ground state solution
of  (\ref{2014-7-17-e2+0}), i.e.,  the least energy  $c_0:=\inf_{(u,v)\in \mathcal{N}}\Phi(u,v)$ defined in (\ref{2014-4-11-e1}) can only be
attained by semi-trivial pairs.   Denote
\be\lab{2014-11-23-we1}
\mu_s(\R_+^N):=\inf\left\{\frac{\int_{\R_+^N} |\nabla u|^2 dx}{(\int_{\R_+^N}\frac{|u|^{2^*(s)}}{|x|^s}dx)^{\frac{2}{2^*(s)}}}: \;u\in D_{0}^{1,2}(\R_+^N)\backslash\{0\}\right\}.
\ee
By the result of  Egnell \cite{Egnell.1992}, $\mu_{s_1}(\R_+^N)$ is achieved  and the extremals are parallel to $U(x)$, a ground state solution of the following problem:
\be\lab{2014-11-23-we2}
\begin{cases}
-\Delta u=\mu_{s_1}(\R_+^N) \frac{u^{2^*(s_1)-1}}{|x|^{s_1}}\;&\hbox{in}\;\R_+^N,\\
u=0\;\;&\hbox{on}\;\partial \R_+^N.
\end{cases}
\ee
Define the functional
\be\lab{2014-11-23-we3}
\Psi_\lambda(u)=\frac{1}{2}\int_{\R_+^N}|\nabla u|^2 dx-\frac{\lambda}{2^*(s_1)}\int_{\R_+^N}\frac{|u|^{2^*(s_1)}}{|x|^{s_1}}dx.
\ee
Then a direct computation shows  that $u$ is a least energy critical point of $\Psi_\lambda$ if and only if
\be\lab{2014-11-23-we6}
u=U_\lambda:=\left(\frac{\mu_{s_1}(\R_+^N)}{\lambda}\right)^{\frac{1}{2^*(s_1)-2}} U,
\ee
where $U$ is a ground state solution of \eqref{2014-11-23-we2}. And the corresponding ground state value is denoted by
\be\lab{2014-11-23-we4}
m_\lambda=\Psi_\lambda(U_\lambda)=\left[\frac{1}{2}-\frac{1}{2^*(s_1)}\right]\left(\mu_{s_1}(\R_+^N)\right)^{\frac{2^*(s_1)}{2^*(s_1)-2}} \lambda^{-\frac{2}{2^*(s_1)-2}}.
\ee
Then we see that $m_\lambda$ is decreasing by $\lambda$ and
\be\lab{2014-11-23-we5}
c_0\leq \min\{m_\lambda, m_\mu\},
\ee
where $c_0$ is defined by \eqref{2014-4-11-e1}.

\bt\lab{2014-12-1-th1}
Assume that $\alpha+\beta=2^*(s_2)$.
If one of the following conditions is satisfied:
\begin{itemize}
\item[(i)]  $\kappa<0$ and $s_2\geq s_1$;
\item[$(ii)$] $\min\{\alpha,\beta\} \frac{(N-s_1)(2-s_2)}{(N-s_2)(2-s_1)}>2, s_2\geq s_1$ and $\kappa>0$ small enough,
\end{itemize}
then we have
$\displaystyle c_0=\min\{m_\lambda, m_\mu\}$.
Moreover,  $c_0$ is achieved by and only by semitrivial solution
$$
\begin{cases}
(U_\lambda, 0) \quad &\hbox{if}\;\lambda>\mu,\\
(0,U_\mu)\quad &\hbox{if}\;\lambda<\mu,\\ (U_\lambda,0)\;\hbox{or}\;(0,U_\lambda)&\hbox{if}\;\lambda=\mu,
\end{cases}
$$
 where $U_\lambda,U_\mu$ are defined by \eqref{2014-11-23-we6}.
\et

\br   Theorem \ref{2014-12-1-th1} means that the system   (\ref{2014-7-17-e2+0}) has only semi-trivial ground state
under the hypotheses of the theorem.\er

\br\lab{2014-11-29-r1}
If $s_1=s_2=s\in (0,2)$, $\min\{\alpha,\beta\}>2$, we must have  $N=3$. In this case,  the assumption  that ``$\kappa$ is small enough" can be removed (see Theorem \ref{2014-11-30-th2}).
\er
\bp
Without loss of generality, we only prove the case of $\lambda>\mu$.
By \eqref{2014-11-23-we5}, we see that $c_0\leq m_\lambda$. By (\ref{bane-1}), we also have $c_0>0$. Now,we proceed by contradiction. Assume that $c_0$ is achieved by some $(u,v)\in \mathscr{D}$ such that $ u\neq 0, v\neq 0$. Without loss of generality, we may assume that $u\geq 0,v\geq 0$ since $c_0$ is the least energy.

$(i)$ If $\kappa<0$, we obtain that
$$\int_{\R_+^N} |\nabla u|^2-\lambda\int_{\R_+^N} \frac{|u|^{2^*(s_1)}}{|x|^{s_1}}dx=\kappa\alpha\int_{\R_+^N} \frac{|u|^\alpha|v|^\beta}{|x|^{s_2}}dx\leq 0.$$
Recalling that
$$ \mu_{s_1}(\R_+^N)\Big(\int_{\R_+^N} \frac{|u|^{2^*(s_1)}}{|x|^{s_1}}dx\Big)^{\frac{2}{2^*(s_1)}}\leq \|u\|^2,$$
if $u\neq 0$,
we obtain that
\be\lab{2014-4-14-we1}
\int_{\R_+^N} \frac{|u|^{2^*(s_1)}}{|x|^{s_1}}dx\geq \big(\frac{\mu_{s_1}(\R_+^N)}{\lambda}\big)^{\frac{2^*(s_1)}{2^*(s_1)-2}}.
\ee
Then by $s_2\geq s_1$,
{\allowdisplaybreaks
\begin{align*} &\big(\frac{1}{2}-\frac{1}{2^*(s_2)}\big)\|u\|^2+\big(\frac{1}{2^*(s_2)}-\frac{1}{2^*(s_1)}\big)\lambda \int_{\R_+^N} \frac{|u|^{2^*(s_1)}}{|x|^{s_1}}dx\\
&\geq \big(\frac{1}{2}-\frac{1}{2^*(s_2)}\big)\mu_{s_1}(\R_+^N)\big(\int_{\R_+^N} \frac{|u|^{2^*(s_1)}}{|x|^{s_1}}dx\big)^{\frac{2}{2^*(s_1)}}+\\
&\quad\quad \big(\frac{1}{2^*(s_2)}-\frac{1}{2^*(s_1)}\big)\lambda \int_{\R_+^N} \frac{|u|^{2^*(s_1)}}{|x|^{s_1}}dx\\
&\geq \big(\frac{1}{2}-\frac{1}{2^*(s_2)}\big)\lambda^{-\frac{2}{2^*(s_1)-2}}\big(\mu_{s_1}(\R_+^N)\big)^{\frac{2^*(s_1)}{2^*(s_1)-2}}+\\
&\quad\quad \big(\frac{1}{2^*(s_2)}-\frac{1}{2^*(s_1)}\big)\lambda^{-\frac{2}{2^*(s_1)-2}}\big(\mu_{s_1}(\R_+^N)\big)^{\frac{2^*(s_1)}{2^*(s_1)-2}}\\
&=\big(\frac{1}{2}-\frac{1}{2^*(s_1)}\big)\lambda^{-\frac{2}{2^*(s_1)-2}}\big(\mu_{s_1}(\R_+^N)\big)^{\frac{2^*(s_1)}{2^*(s_1)-2}}\\
&=m_\lambda.
\end{align*}
}
Similarly, if $v\neq 0$, we have
\be\lab{2014-4-15-e1}
\int_{\R_+^N} \frac{|v|^{2^*(s_1)}}{|x|^{s_1}}dx\geq \big(\frac{\mu_{s_1}({\R_+^N})}{\mu}\big)^{\frac{2^*(s_1)}{2^*(s_1)-2}}
\ee
and
$$\big(\frac{1}{2}-\frac{1}{2^*(s_2)}\big)\|v\|^2+\big(\frac{1}{2^*(s_2)}-\frac{1}{2^*(s_1)}\big)\mu \int_{\R_+^N} \frac{|v|^{2^*(s_1)}}{|x|^{s_1}}dx\geq m_\mu>m_\lambda.$$
Then,
\begin{align*}
c_0=&\Phi(u,v)=\big(\frac{1}{2}-\frac{1}{2^*(s_2)}\big)a(u,v)+\big(\frac{1}{2^*(s_2)}-\frac{1}{2^*(s_1)}\big)b(u,v)\\
=&\big(\frac{1}{2}-\frac{1}{2^*(s_2)}\big)\|u\|^2+\big(\frac{1}{2^*(s_2)}-\frac{1}{2^*(s_1)}\big)\lambda \int_{\R_+^N} \frac{|u|^{2^*(s_1)}}{|x|^{s_1}}dx\\
&+\big(\frac{1}{2}-\frac{1}{2^*(s_2)}\big)\|v\|^2+\big(\frac{1}{2^*(s_2)}-\frac{1}{2^*(s_1)}\big)\mu \int_{\R_+^N} \frac{|v|^{2^*(s_1)}}{|x|^{s_1}}dx\\
\geq&\begin{cases}m_\lambda\quad &\hbox{if}\;v=0\\
m_\lambda+m_\mu &\hbox{if}\;v\neq 0.
\end{cases}
\end{align*}
Hence, $c_0=m_\lambda$ is proved and we see that $v=0$, i.e., $(u,v)=(U_\lambda, 0)$.

$(ii)$ If $\kappa>0$, we denote
$\displaystyle\sigma:=\int_{\R_+^N} \frac{|u|^{2^*(s_1)}}{|x|^{s_1}}dx,\;\delta:=\int_{\R_+^N} \frac{|v|^{2^*(s_1)}}{|x|^{s_1}}dx.$
Then we have
$$\sigma\leq \big(\frac{\mu_{s_1}(\R_+^N)}{\lambda}\big)^{\frac{2^*(s_1)}{2^*(s_1)-2}}.$$
If not, apply the above similar arguments, we have
$$\Phi(u,v)\geq \big(\frac{1}{2}-\frac{1}{2^*(s_2)}\big)\|u\|^2+\big(\frac{1}{2^*(s_2)}-\frac{1}{2^*(s_1)}\big)\lambda \int_{\R_+^N} \frac{|u|^{2^*(s_1)}}{|x|^{s_1}}dx>m_\lambda,$$
a contradiction.
Similarly, we also have
$$\delta\leq \big(\frac{\mu_{s_1}(\R_+^N)}{\mu}\big)^{\frac{2^*(s_1)}{2^*(s_1)-2}}.$$
Similar to the arguments of Lemma \ref{2014-11-23-xl1}, we have
$\Phi(u,v)\geq (\frac{1}{2}-\frac{1}{2^*(s_2)}) \|(u,v)\|_{\mathscr{D}}^{2}$.
Hence, $\|u\|^2, \|v\|^2\leq (\frac{1}{2}-\frac{1}{2^*(s_2)})^{-1}c_0$.
By Corollary \ref{2014-4-22-interpolation-corollary}, under the assumption of $\min\{\alpha,\beta\} \frac{(N-s_1)(2-s_2)}{(N-s_2)(2-s_1)}>2$, we can choose some proper  $\eta_1\geq 2, \eta_2\geq 2$ and $C>0$ such that
\be\lab{2014-5-9-we1}
\int_{\R_+^N} \frac{|u|^\alpha |v|^\beta}{|x|^{s_2}}dx\leq C |u|_{2^*(s_1), s_1}^{\eta_1}=C\sigma^{\frac{\eta_1}{2^*(s_1)}}
\ee
and
\be\lab{2014-5-9-we2}
\int_{\R_+^N} \frac{|u|^\alpha |v|^\beta}{|x|^{s_2}}dx\leq C |v|_{2^*(s_1), s_1}^{\eta_2}=C\delta^{\frac{\eta_2}{2^*(s_1)}}.
\ee
It follows that there exists some $C>0$ such that
\be\lab{2014-4-15-e2}
\mu_{s_1}(\R_+^N)\sigma^{\frac{2}{2^*(s_1)}}-\lambda \sigma\leq \kappa C \sigma^{\frac{\eta_1}{2^*(s_1)}}
\ee
and that
\be\lab{2014-4-15-e3}
\mu_{s_1}(\R_+^N)\delta^{\frac{2}{2^*(s_1)}}-\mu \delta\leq \kappa C \delta^{\frac{\eta_2}{2^*(s_1)}}.
\ee
Define $g_i:\R_+\mapsto \R_+, i=1,2$ with
$$ g_1(t):=\lambda t^{\frac{2^*(s_1)-2}{2^*(s_1)}}+\kappa C t^{\frac{\eta_1-2}{2^*(s_1)}}$$
and
$$ g_2(t):=\mu t^{\frac{2^*(s_1)-2}{2^*(s_1)}}+\kappa C t^{\frac{\eta_2-2}{2^*(s_1)}}.$$
Since $\eta_1, \eta_2\geq 2$, $g_i(t)$ is strictly increasing in terms of  $t$. It is easy to check that there exists some $\kappa_0>0$ such that when $\kappa<\kappa_0$, a direct calculation shows   that
$$ g_1\Big(\frac{1}{2}\big(\frac{\mu_{s_1}(\R_+^N)}{\lambda}\big)^{\frac{2^*(s_1)}{2^*(s_1)-2}}\Big)<\mu_{s_1}(\R_+^N)$$
and
$$ g_2\Big(\frac{1}{2}\big(\frac{\mu_{s_1}(\R_+^N)}{\mu}\big)^{\frac{2^*(s_1)}{2^*(s_1)-2}}\Big)<\mu_{s_1}(\R_+^N).$$
Hence, if $\sigma\neq 0, \delta\neq 0$, by (\ref{2014-4-15-e2}) and (\ref{2014-4-15-e3}), we obtain that
$$\sigma>\frac{1}{2}\big(\frac{\mu_{s_1}(\R_+^N)}{\lambda}\big)^{\frac{2^*(s_1)}{2^*(s_1)-2}}$$
and
$$ \delta>\frac{1}{2}\big(\frac{\mu_{s_1}(\R_+^N)}{\mu}\big)^{\frac{2^*(s_1)}{2^*(s_1)-2}}.$$
Then{\allowdisplaybreaks
\begin{align*} &\big(\frac{1}{2}-\frac{1}{2^*(s_2)}\big)\|u\|^2+\big(\frac{1}{2^*(s_2)}-\frac{1}{2^*(s_1)}\big)\lambda \int_{\R_+^N} \frac{|u|^{2^*(s_1)}}{|x|^{s_1}}dx\\
\geq&\big(\frac{1}{2}-\frac{1}{2^*(s_2)}\big)\mu_{s_1}(\R_+^N)\sigma^{\frac{2}{2^*(s_1)}}+\big(\frac{1}{2^*(s_2)}-\frac{1}{2^*(s_1)}\big)\lambda\sigma\\
\geq&\big(\frac{1}{2}-\frac{1}{2^*(s_2)}\big)\mu_{s_1}(\R_+^N)\big[\frac{1}{2}\big(\frac{\mu_{s_1}(\R_+^N)}{\lambda}\big)^{\frac{2^*(s_1)}{2^*(s_1)-2}}\big]^{\frac{2}{2^*(s_1)}}\\
&+\big(\frac{1}{2^*(s_2)}-\frac{1}{2^*(s_1)}\big)\lambda \frac{1}{2}\big(\frac{\mu_{s_1}(\R_+^N)}{\lambda}\big)^{\frac{2^*(s_1)}{2^*(s_1)-2}}\\
>&\frac{1}{2}\big(\frac{1}{2}-\frac{1}{2^*(s_1)}\big)\lambda^{-\frac{2}{2^*(s_1)-2}}\mu_{s_1}(\R_+^N)^{\frac{2^*(s_1)}{2^*(s_1)-2}}=\frac{1}{2}m_\lambda.
\end{align*}}
Similarly, we have
$$\big(\frac{1}{2}-\frac{1}{2^*(s_2)}\big)\|v\|^2+\big(\frac{1}{2^*(s_2)}-\frac{1}{2^*(s_1)}\big)\mu \int_{\R_+^N} \frac{|v|^{2^*(s_1)}}{|x|^{s_1}}dx>\frac{1}{2}m_\mu>\frac{1}{2}m_\lambda.$$
Thus,
\begin{align*}
\Phi(u,v)
=&\big(\frac{1}{2}-\frac{1}{2^*(s_2)}\big)\|u\|^2+\big(\frac{1}{2^*(s_2)}-\frac{1}{2^*(s_1)}\big)\lambda \int_{\R_+^N} \frac{|u|^{2^*(s_1)}}{|x|^{s_1}}dx\\
&+\big(\frac{1}{2}-\frac{1}{2^*(s_2)}\big)\|v\|^2+\big(\frac{1}{2^*(s_2)}-\frac{1}{2^*(s_1)}\big)\mu \int_{\R_+^N} \frac{|v|^{2^*(s_1)}}{|x|^{s_1}}dx>m_\lambda,
\end{align*}
a contradiction.

The arguments above imply that $\sigma=0$ or $\delta=0$, i.e., $u=0$ or $v=0$. If $u=0$, then $v\neq 0$ is a critical point of $\Psi_\mu$ and then $\Phi(u,v)=\Psi_\mu(v)\geq m_\mu>m_\lambda$, also a contradiction.
Thus, we obtain that $u\neq 0, v=0$. Hence $u$ is a critical point of $\Psi_\lambda$, and $c_0=\Phi(u,v)=\Psi_\lambda(u)\geq m_\lambda$.
Then, we have $c_0=m_\lambda$ and $u=U_\lambda$.
\ep

\br   We remark that  the  Theorem \ref{2014-12-1-th1} of this section is valid for any cone  $\Omega$.
\er

 \vskip0.4in

\s{Preliminaries for the existence results}
\renewcommand{\theequation}{6.\arabic{equation}}
\renewcommand{\theremark}{6.\arabic{remark}}
\renewcommand{\thedefinition}{6.\arabic{definition}}
\br\lab{2014-11-27-zr1}
Without loss of generality, we only consider the case of $\Omega=\R_+^N$. We remark that  the results of this section are still valid for any domain $\Omega$ as long as  $\mu_{s_1}(\Omega)$ is attained (e.g. $\Omega$ is a cone).
\er
Since  the system  \eqref{2014-7-17-e2+0} possesses semitrivial solution  $(u, v)$, we are interested in the nontrivial solutions. Firstly, we recall the following result   due to Ghoussoub and Robert \cite[Theorem 1.2]{GhoussoubRobert.2006} (see also \cite[Lemma 2.1]{HsiaLinWadade.2010}, \cite[Lemma 2.6]{LinWadadeothers.2012}) for the  scalar  equation.

\bl\lab{CL-2} (\cite[Theorem 1.2]{GhoussoubRobert.2006})  Let $u\in D_{0}^{1,2}(\R_+^N)$ be an entire solution to the problem
\be\lab{BP}
\begin{cases}
\Delta u+\frac{u^{2^*(s_1)-1}}{|x|^{s_1}}=0\quad &\hbox{in}\;\R_+^N,\\
u>0\;\;\hbox{in}\;\R_+^N\;\;\;\hbox{and}& u=0\;\hbox{on}\;\partial\R_+^N.
\end{cases}
\ee
Then, the following hold:
\begin{itemize}
\item[(i)]$$\begin{cases}u\in C^2(\overline{\R_+^N})\quad &\hbox{if}\;s_1<1+\frac{2}{N},\\
u\in C^{1,\beta}(\overline{\R_+^N})\quad &\hbox{for all}\;0<\beta<1\;\hbox{if}\;s_1=1+\frac{2}{N},\\
u\in C^{1,\beta}(\overline{\R_+^N})\quad &\hbox{for all}\;0<\beta<\frac{N(2-s_1)}{N-2}\;\hbox{if}\;s_1>1+\frac{2}{N}.
 \end{cases}$$
\item[(ii)] There is a constant $C$ such that $|u(x)|\leq C (1+|x|)^{1-N}$ and $|\nabla u(x)|\leq C(1+|x|)^{-N}$.
\item[(iii)]$u(x', x_N)$ is axially symmetric with respect to the $x_N$-axis, i.e., $u(x', x_N)=u(|x'|, x_N)$, where $x'=(x_1,\cdots,x_{N-1})$.
\end{itemize}
\el

\subsection{Existence of positive solution for  the  case : $\lambda=\mu \left(\frac{\beta}{\alpha}\right)^{\frac{2^*(s_1)-2}{2}}$}
The following result is essentially due to \cite[Theorem 1.2]{LiLin.2012}:
\bl\lab{2014-11-23-xbul1}
Let $N\geq 3$, $s_1,s_2\in (0,2),\lambda>0$, $\alpha>1,\beta>1$ and $\alpha+\beta=2^*(s_2)$. Then the following problem
\be\lab{2014-11-23-xbuel1}
\begin{cases}
-\Delta w-\lambda \frac{w^{2^*(s_1)-1}}{|x|^{s_1}}-\kappa \alpha (\frac{\beta}{\alpha})^{\frac{\beta}{2}}\frac{w^{2^*(s_2)-1}}{|x|^{s_2}}=0\;\hbox{in}\;\R_+^N,\\
w(x)\in D_{0}^{1,2}(\R_+^N),\;w(x)>0\;\hbox{in}\;\R_+^N,
\end{cases}
\ee
has a least-energy solution provided further one of the following holds:
\begin{itemize}
\item[$(i)$] $0<s_1<s_2<2$ and $\kappa\in \R$.
\item[$(ii)$] $s_1>s_2$ and $\kappa\geq0$.
\item[$(iii)$] $s_1=s_2$ and $\kappa>-\lambda \frac{1}{\alpha}(\frac{\alpha}{\beta})^{\frac{\beta}{2}}$.
\end{itemize}
\el
\br\lab{2014-11-23-xbur1}
The case of $\kappa=0$ or $s_1=s_2$ with $\kappa>-\lambda \frac{1}{\alpha}(\frac{\alpha}{\beta})^{\frac{\beta}{2}}$, \eqref{2014-11-23-xbuel1} is essentially the problem \eqref{BP}. And the existence result was firstly given by Egnell \cite{Egnell.1992}.
\er

\bc\lab{2014-11-23-cro1}
Let $N\geq 3$, $s_1,s_2\in (0,2),\mu>0,\kappa\neq 0$, $\alpha>1,\beta>1$ and $\alpha+\beta=2^*(s_2)$. If $\lambda=\mu (\frac{\beta}{\alpha})^{\frac{2^*(s_1)}{2}}$, then $(w, \sqrt{\frac{\beta}{\alpha}}w)$ is a positive solution of \eqref{2014-7-17-e2} provided further one of the following holds:
\begin{itemize}
\item[$(i)$] $0<s_1<s_2<2$ and $\kappa\in \R\backslash\{0\}$.
\item[$(ii)$] $s_1>s_2$ and $\kappa>0$.
\item[$(iii)$] $s_1=s_2$ and $\kappa>-\lambda \frac{1}{\alpha}(\frac{\alpha}{\beta})^{\frac{\beta}{2}}$.
\end{itemize}
Here $w$ is a least-energy solution of \eqref{2014-11-23-xbuel1}.
\ec
\bp
This proof can be got through via  a direct computation. We omit the details.
\ep

\bc\lab{2014-11-23-wcor1}
 Assume that $N\geq 3, \alpha>1,\beta>1, \alpha+\beta=2^*(s_2), \lambda=\mu (\frac{\beta}{\alpha})^{\frac{2^*(s_1)}{2}}>0$ and one of the following holds
\begin{itemize}
\item[$(i)$] $0<s_1<s_2<2,\kappa<0$,
\item[$(ii)$]$s_1=s_2\in (0,2), -\lambda \frac{1}{\alpha}(\frac{\alpha}{\beta})^{\frac{\beta}{2}}<\kappa <0$,
\item[$(iii)$]$\min\{\alpha,\beta\} \frac{(N-s_1)(2-s_2)}{(N-s_2)(2-s_1)}>2, s_2\geq s_1$ and $\kappa>0$ small enough.
\end{itemize}
Then $(w, \sqrt{\frac{\beta}{\alpha}}w)$ is a positive solution to equation \eqref{2014-7-17-e2} but problem \eqref{2014-7-17-e2} has no nontrivial ground state solution.
\ec
\bp
It is a  straightforward consequence of  Theorem \ref{2014-12-1-th1} and Corollary \ref{2014-11-23-cro1}.
\ep

\subsection{Estimation on the upper bound  of  $\displaystyle c_0:=\inf_{(u,v)\in \mathcal{N}}\Phi(u,v)$}

In order  to prove the existence of positive ground state solution to  the equation \eqref{2014-7-17-e2+0}, we have to give an estimation on the
upper bound  of $c_0$,  including  the cases of $s_1=s_2$ and $s_1\neq s_2$.

Let $1<\alpha, 1<\beta, \alpha+\beta= 2^*(s_2)$. Let $u:=U_\lambda$ be a function defined by \eqref{2014-11-23-we6}.
 Then we have $u>0$ in $\R_+^N$ and
$$c(u,v):=\int_{\R_+^N}\frac{|u|^\alpha |v|^\beta}{|x|^{s_2}}dx>0, \quad \;\forall\;v\in D_{0}^{1,2}(\R_+^N)\backslash\{0\}.$$
Assume that the assumptions of Lemma \ref{2014-11-23-xl1} are satisfied, i.e., one of the following holds:
\begin{itemize}
\item[(i)]$\kappa>0$,
\item[(ii)] $\kappa<0$ and $s_2>s_1$,
\item[(iii)] $s_2=s_1$  and $ \kappa<0$ small enough,
\end{itemize}
then we see that for any $\varepsilon\in \R$,  there exists a  unique positive number $t(\varepsilon)>0$ such that $\big(t(\varepsilon)u, t(\varepsilon)\varepsilon v\big)\in \mathcal{N}$. The function $t(\varepsilon):\R\mapsto \R_+$ is implicitly defined by the equation
\begin{align}\lab{2014-4-12-xe1}
\|u\|^2+\varepsilon^2\|v\|^2=&\Big[\lambda|u|_{2^*(s_1),s_1}^{2^*(s_1)}+\mu|v|_{2^*(s_1),s_1}^{2^*(s_1)}|\varepsilon|^{2^*(s_1)}\Big]\big[t(\varepsilon)\big]^{2^*(s_1)-2}\\
&+\kappa 2^*(s_2)c(u,v)\big[t(\varepsilon)\big]^{2^*(s_2)-2}|\varepsilon|^\beta.\nonumber
\end{align}
We notice that $t(0)=1$. Moreover, from the Implicit Function Theorem, it follows that $t(\varepsilon)\in C^1(\R)$ and $t'(\varepsilon)=\frac{P_v(\varepsilon)}{Q_v(\varepsilon)}$, where
\begin{align*}
Q_v(\varepsilon):=&\big[2^*(s_1)-2\big]\Big[\lambda|u|_{2^*(s_1),s_1}^{2^*(s_1)}+\mu|v|_{2^*(s_1),s_1}^{2^*(s_1)}|s|^{2^*(s_1)}\Big]\big[t(\varepsilon)\big]^{2^*(s_1)-3}\nonumber\\
&+\kappa 2^*(s_2)[2^*(s_2)-2]c(u,v)\big[t(\varepsilon)\big]^{2^*(s_2)-3}|\varepsilon|^\beta
\end{align*}
and
\begin{align*}
P_v(\varepsilon):=&2\|v\|^2 \varepsilon-2^*(s_1)\mu|v|_{2^*(s_1),s_1}^{2^*(s_1)}\big[t(\varepsilon)\big]^{2^*(s_1)-2}|\varepsilon|^{2^*(s_1)-2}\varepsilon\nonumber\\
&-\kappa 2^*(s_2)\beta c(u,v)\big[t(\varepsilon)\big]^{2^*(s_2)-2}|\varepsilon|^{\beta-2}\varepsilon.
\end{align*}

\bl\lab{2014-11-24-l1}({\bf  The case of  $\beta<2$})
Assume that  $1<\alpha, 1<\beta<2, \alpha+\beta=2^*(s_2)$ and one of the following holds:
\begin{itemize}
\item[(i)]$\kappa>0$.
\item[(ii)] $\kappa<0$ and $s_2>s_1$.
\item[(iii)] $s_2=s_1$  and $\kappa<0$  with $|\kappa|$ small enough.
\end{itemize}
Let $U_\lambda:=\big(\frac{\mu_{s_1}(\R_+^N)}{\lambda}\big)^{\frac{1}{2^*(s_1)-2}}U$, where $U$ is a ground state solution of \eqref{2014-11-23-we2}.
Then
\begin{itemize}
\item[(a)]if $\kappa<0$, $(u,0)$ is a local minimum point of $\Phi$ in $\mathcal{N}$.
\item[(b)] if $\kappa>0$, then
$$ c_0:=\inf_{(\phi,\varphi)\in \mathcal{N}}\Phi(\phi,\varphi)<\Phi(U_\lambda,0)=\Psi_{\lambda}(U_\lambda)=m_{\lambda}.$$
\end{itemize}
\el
\bp  Let $u:=U_\lambda$ and take  any $v\in D_{0}^{1,2}(\R_+^N)\backslash\{0\}.$
When $\beta<2$, we have
$$P_v(\varepsilon)=-\kappa2^*(s_2)\beta c(u,v)|\varepsilon|^{\beta-2}\varepsilon\big(1+o(1)\big)\;\hbox{as}\;\varepsilon\rightarrow 0$$
and
$$Q_v(\varepsilon)=\Big(\big[2^*(s_1)-2\big]\lambda|u|_{2^*(s_1),s_1}^{2^*(s_1)}
\Big)\big(1+o(1)\big)\;\hbox{as}\;\varepsilon\rightarrow 0.$$
Hence,
$$ t'(\varepsilon)=-M(v)\beta |\varepsilon|^{\beta-2}\varepsilon\big(1+o(1)\big),$$
where
$$ M(v):=\frac{\kappa2^*(s_2)c(u,v)}{\big[2^*(s_1)-2\big]\lambda|u|_{2^*(s_1),s_1}^{2^*(s_1)}
}.$$
By the Taylor formula, we obtain that
$$t(\varepsilon)=1-M(v)|\varepsilon|^\beta\big(1+o(1)\big),$$
and
$$\big[t(\varepsilon)\big]^{2^*(s_i)}=1-2^*(s_i)M(v)|\varepsilon|^\beta\big(1+o(1)\big)\quad (i=1,2).$$
Noting that   for any $(\phi,\varphi)\in \mathcal{N}$, we have
$$
\Phi(\phi,\varphi)=\big(\frac{1}{2}-\frac{1}{2^*(s_1)}\big)b(\phi,\varphi)+\frac{2^*(s_2)-2}{2}\kappa c(\phi,\varphi),
$$
where $b(\phi,\varphi),c(\phi,\varphi)$ are defined by \eqref{2014-11-23-xe1}.
Thus,{\allowdisplaybreaks
\begin{align*}
&\Phi\big(t(\varepsilon)u, t(\varepsilon)\varepsilon v\big)-\Phi(u,0)\\
=&\big(\frac{1}{2}-\frac{1}{2^*(s_1)}\big)\Big[\lambda|u|_{2^*(s_1),s_1}^{2^*(s_1)}\big[t(\varepsilon)\big]^{2^*(s_1)}
+\mu|v|_{2^*(s_1),s_1}^{2^*(s_1)}\big[t(\varepsilon)\big]^{2^*(s_1)}|\varepsilon|^{2^*(s_1)}\\
&-\lambda|u|_{2^*(s_1),s_1}^{2^*(s_1)}\Big]
+\frac{2^*(s_2)-2}{2}\kappa c(u,v)\big[t(\varepsilon)\big]^{2^*(s_2)}|\varepsilon|^\beta\\
=&\big(\frac{1}{2}-\frac{1}{2^*(s_1)}\big)\lambda|u|_{2^*(s_1),s_1}^{2^*(s_1)}\big[-2^*(s_1)M(v)\big]|\varepsilon|^\beta \big(1+o(1)\big)\\
&+\frac{2^*(s_2)-2}{2}\kappa c(u,v)\big[t(\varepsilon)\big]^{2^*(s_2)}|\varepsilon|^\beta\\
=&-\frac{\kappa2^*(s_2)}{2}c(u,v)|\varepsilon|^\beta\big(1+o(1)\big)+\frac{2^*(s_2)-2}{2}\kappa c(u,v)\big(1+o(1)\big)|\varepsilon|^\beta\\
=&-\kappa|\varepsilon|^\beta c(u,v)\big(1+o(1)\big),
\end{align*}
}
which implies the results  since  $c(u,v)>0$ for any $0\neq v\in D_{0}^{1,2}(\R_+^N)$.
\ep

\bl\lab{2014-11-24-l2}({\bf the case of $\beta>2$})
Assume that  $1<\alpha, 2<\beta, \alpha+\beta=2^*(s_2)$ and one of the following holds:
\begin{itemize}
\item[(i)]$\kappa>0$.
\item[(ii)] $\kappa<0$ and $s_2>s_1$.
\item[(iii)] $s_2=s_1$  and $ \kappa<0$ with $|\kappa|$ small enough.
\end{itemize}
Let $U_\lambda:=\big(\frac{\mu_{s_1}(\R_+^N)}{\lambda}\big)^{\frac{1}{2^*(s_1)-2}}U$, where $U$ is a ground state solution of \eqref{2014-11-23-we2}.
 Then $(U_\lambda,0)$ is a local minimum point of $\Phi$ in $\mathcal{N}$.
\el
\bp   Let $u:=U_\lambda$ and take  any $v\in D_{0}^{1,2}(\R_+^N)\backslash\{0\}.$
When $\beta>2$, we have
$$P_v(\varepsilon)=2\|v\|^2\varepsilon\big(1+o(1)\big),
Q_v(\varepsilon)=\Big(\big[2^*(s_1)-2\big]\lambda|u|_{2^*(s_1),s_1}^{2^*(s_1)}
\Big)\big(1+o(1)\big).$$
Hence,
$$ t'(\varepsilon)=2\tilde{M}(v) \varepsilon\big(1+o(1)\big),$$
where
$$ \tilde{M}(v):=\frac{\|v\|^2}{\big[2^*(s_1)-2\big]\lambda|u|_{2^*(s_1),s_1}^{2^*(s_1)}
}.$$
By the Taylor formula, we obtain that
$$t(\varepsilon)=1+\tilde{M}(v)|\varepsilon|^2\big(1+o(1)\big)$$
and
$$\big[t(\varepsilon)\big]^{2^*(s_i)}=1+2^*(s_i)\tilde{M}(v)|\varepsilon|^2\big(1+o(1)\big)\quad (i=1,2).$$
Hence a direct computation shows that
\begin{align*}
&\Phi\big(t(\varepsilon)u, t(\varepsilon)\varepsilon v\big)-\Phi(u,0)\\
=&\big(\frac{1}{2}-\frac{1}{2^*(s_1)}\big)\lambda|u|_{2^*(s_1),s_1}^{2^*(s_1)}\big[2^*(s_1)\tilde{M}(v)\big]|\varepsilon|^2 \big(1+o(1)\big)+o(|\varepsilon|^2)\\
=&\frac{1}{2}\|v\|^2|\varepsilon|^2\big(1+o(1)\big)
>0\;\hbox{when $\varepsilon$ is small enough}.
\end{align*}
\ep

Define $\displaystyle \eta_1:=\inf_{v\in \Xi}\|v\|^2$, where
$$ \Xi:=\{v\in D_{0}^{1,2}(\R_+^N): \quad \int_{\R_+^N} \frac{|U_\lambda|^\alpha |v|^2}{|x|^{s_2}}dx=1\}.$$
We note that by the  Hardy-Sobolev inequality, $U_\lambda\in L^{2^*(s_2)}(\R_+^N, \frac{dx}{|x|^{s_2}})$, then by the  H\"older inequality, $\int_{\R_+^N} \frac{|U_\lambda|^\alpha |v|^2}{|x|^{s_2}}dx$ is well defined  for all $v\in D_{0}^{1,2}(\R_+^N)$ when $\alpha=2^*(s_2)-2$.

Define $\langle \phi, \psi\rangle :=\int_{\Omega}\frac{|U_\lambda|^\alpha \phi\psi}{|x|^{s_2}}dx$, then it is easy to check that $\langle \cdot, \cdot\rangle$ is an inner product. We say that $\phi$ and $\psi$ are  orthogonal if and only if $\langle \phi, \psi\rangle=0$.
Then we have the following result:

\bl\lab{2014-11-24-l3}
Assume that $1<\alpha=2^*(s_2)-2, \beta=2$,
then there exists
$\eta_1>0$ and
 some $0<v\in D_{0}^{1,2}(\R_+^N)$ such that
\be\lab{2014-4-14-gze1}
-\Delta v=\eta_1 \frac{|U_\lambda|^\alpha}{|x|^{s_2}}v\;\hbox{in}\;\R_+^N,\;
v=0\;\hbox{on}\;\partial\R_+^N.
\ee
Furthermore, the eigenvalue $\eta_1$ is simple and  satisfying
\be\lab{2014-4-14-ze1}
\int_{\R_+^N}\frac{|U_\lambda|^\alpha|v|^2}{|x|^{s_2}}dx\leq \frac{1}{\eta_1}\|v\|^2\;\hbox{for all}\;v\in D_{0}^{1,2}(\R_+^N).
\ee
In particular, if $s_1=s_2=s$, we have
\be\lab{2014-11-27-lae2}
\eta_1=\lambda,
\ee
and it is only attained by $v=U_\lambda$.
\el
\bp
It is easy to see that $\eta_1\geq 0$ under our assumptions. Let $\{v_n\}\subset \Xi$ be  such that $\|v_n\|^2\rightarrow \eta_1$. Then $\{v_n\}$ is bounded in $D_{0}^{1,2}(\R_+^N)$. Going to a subsequence if necessary, we may assume that $v_n\rightharpoonup v_0$ in $D_{0}^{1,2}(\R_+^N)$ and $v_n\rightarrow v_0$ a.e. in $\R_+^N$.
By H\"{o}lder inequality, we have
$$
\Big|\int_\Lambda \frac{u^\alpha v_n^2- u^\alpha v_0^2}{|x|^{s_2}}dx\Big|
\leq\Big(\int_{\Lambda} \frac{u^{\alpha+2}}{|x|^{s_2}}dx\Big)^{\frac{\alpha}{\alpha+2}}\Big(\int_{\Lambda} \frac{|v_n^2-v^2|^{\frac{\alpha+2}{2}}}{|x|^{s_2}}\Big)^{\frac{2}{\alpha+2}}\rightarrow 0
$$
as $|\Lambda|\rightarrow 0$ due to the absolute continuity of the integral and the boundness of $v_n$.
Similarly, we also have that $\{\frac{u^\alpha v_n^2}{|x|^{s_2}}\}$ is a tight sequence, i.e., for $\varepsilon>0$, there exists some $R_\varepsilon>0$ such that
\be\lab{2014-11-24-xe1}
\Big|\int_{\R_+^N\cap B_{R_\varepsilon}^{c}} \frac{u^\alpha v_n^2}{|x|^{s_2}}dx\Big|\leq \varepsilon\;\hbox{uniformly for all}\; n\in\NN.
\ee
Combine with  the Egoroff Theorem, it is easy to prove that
\be\lab{2014-6-25-xe1}
\int_{\R_+^N} \frac{u^\alpha v_n^2}{|x|^{s_2}}dx\rightarrow \int_{\R_+^N} \frac{u^\alpha v_0^2}{|x|^{s_2}}dx.
\ee
Hence, we prove that
$$D_{0}^{1,2}(\R_+^N)\mapsto \R\;\hbox{with}\;\chi(v)=\int_{\R_+^N} \frac{|u|^\alpha|v|^2}{|x|^{s_2}}dx$$
is weak continuous, which implies that $\Xi$ is weak closed.
Hence,  $v_0\in \Xi$ and we have
$$\|v_0\|^2\leq \liminf_{n\rightarrow \infty}\|v_n\|^2=\eta_1.$$
On the other hand, by the definition of $\eta_1$ and $v_0\in \Xi$, we have
$\displaystyle\|v_0\|^2\geq \eta_1$.
Thus, $v_0$ is a minimizer of $\|v\|^2$ constraint on $\Xi$.
It is easy to see that $|v_0|$ is also a minimizer. Hence, we may assume that $v_0\geq 0$ without loss of generality.
We see that there exists some Lagrange multiplier $\eta\in \R$ such that
$$ -\Delta v_0=\eta \frac{|u|^\alpha v_0}{|x|^{s_2}}.$$
It follows that $\eta=\eta_1$. Since $v_0\in \Xi$, we get that $v_0\neq 0$ and $\eta_1>0$.

Let $a(x):=\eta_1\frac{|u|^\alpha}{|x|^{s_2}}$, it is easy to
see that $a(x)\in L_{loc}^{\frac{N}{2}}(\R_+^N)$.
Then the Br\'{e}zis-Kato theorem in \cite{BrezisKato.1978} implies that $v\in L_{loc}^{r}(\R_+^N)$ for all $1\leq r<\infty$. Then $v_0\in W_{loc}^{2,r}(\R_+^N)$ for all $1\leq r<\infty$. By  the elliptic regularity theory, $v_0\in C^2(\R_+^N)$.
Finally, by the maximum principle, we obtain that $v_0$ is positive. Finally, $(\ref{2014-4-14-ze1})$ is an easy conclusion from the definition of $\eta_1$.

It is standard to prove that $\eta_1$ is simple, we omit the details. Next, we will compute the  value of $\eta_1$ when $s_1=s_2=s$.
A direct computation shows that
\be\lab{2014-11-27-glae3}
-\Delta U_\lambda=\lambda\frac{U_{\lambda}^{2^*(s)-1}}{|x|^s}.
\ee
Testing \eqref{2014-4-14-gze1} by $U_\lambda$, we have
\be\lab{2014-11-27-lae3}
\int_{\R_+^N}(\nabla v\cdot \nabla U_\lambda)dx=\eta_1\int_{\R_+^N}\frac{U_\lambda^\alpha}{|x|^s}vU_\lambda dx.
\ee
Testing \eqref{2014-11-27-glae3} by $v$, we also have
\be\lab{2014-11-27-lae4}
\int_{\R_+^N} (\nabla U_\lambda\cdot\nabla v)dx=\lambda\int_{\R_+^N}\frac{U_\lambda^\alpha}{|x|^s}U_\lambda v dx.
\ee
Hence,
\be\lab{2014-11-27-lae5}
(\eta_1-\lambda)\int_{\R_+^N}\frac{U_\lambda^\alpha}{|x|^s}U_\lambda vdx=0.
\ee
Since $v$ and $U_\lambda$ are positive, we obtain that $\eta_1=\lambda$.
\ep

\bl\lab{2014-11-24-l4}({\bf the case of $\beta=2$})
Assume that $1<\alpha=2^*(s_2)-2, \beta=2$ and one of the following holds:
\begin{itemize}
\item[(i)]$\kappa>0$.
\item[(ii)] $\kappa<0$ and $s_2>s_1$.
\item[(iii)] $s_2=s_1$  and $ \kappa<0$ with $|\kappa|$ small enough.
\end{itemize}
Let $U_\lambda:=\big(\frac{\mu_{s_1}(\R_+^N)}{\lambda}\big)^{\frac{1}{2^*(s_1)-2}}U$, where $U$ is a ground state solution of \eqref{2014-11-23-we2}. Then there exists a positive $k_0=\frac{\eta_1}{2^*(s_2)}$ such that
\begin{itemize}
\item[$(a)$]if $\kappa<0$, $(U_\lambda,0)$ is a local minimum point of $\Phi$ in $\mathcal{N}$.
\item[$(b)$]if $0<\kappa<k_0$, then $(U_\lambda,0)$ is a local minimum point of $\Phi$ in $\mathcal{N}$.
\item[$(c)$] if $\kappa>k_0$, then
$$ c_0:=\inf_{(\phi,\varphi)\in \mathcal{N}}\Phi(\phi,\varphi)<\Phi(U_\lambda,0)=\Psi_{\lambda}(U_\lambda)=m_\lambda,$$
\end{itemize}
where  $\displaystyle \eta_1:=\inf_{v\in \Xi}\|v\|^2$, with
$$\Xi:=\{v\in D_{0}^{1,2}(\R_+^N):   \; \int_{\R_+^N}\frac{|U_\lambda|^\alpha |v|^2}{|x|^{s_2}}dx=1\}.$$
\el
\bp

 Let $u:=U_\lambda$ and take  any $v\in D_{0}^{1,2}(\R_+^N)\backslash\{0\}.$
In this case, we have
$$P_v(\varepsilon)=2\Big(\|v\|^2-\kappa 2^*(s_2) c(u,v)\Big)\varepsilon\big(1+o(1)\big), Q_v(\varepsilon)=\Big(\big[2^*(s_1)-2\big]\lambda|u|_{2^*(s_1),s_1}^{2^*(s_1)}
\Big)\big(1+o(1)\big).$$
Hence,
$\displaystyle t'(\varepsilon)=2\bar{M}(v) \varepsilon\big(1+o(1)\big),$
where
$\displaystyle \bar{M}(v):=\frac{\|v\|^2-\kappa 2^*(s_2) c(u,v)}{\big[2^*(s_1)-2\big]\lambda|u|_{2^*(s_1),s_1}^{2^*(s_1)}
}.$
Similar to the arguments above, we obtain that
$$\Phi\big(t(s)u, t(s)sv\big)-\Phi(u,0)=\frac{1}{2}\big(\|v\|^2-\kappa 2^*(s_2) c(u,v)\big)|s|^2\big(1+o(1)\big).$$
If $\kappa<0$, we obtain the result of $(a)$.

Since $u>0$ is given, by Lemma \ref{2014-11-24-l3}, we have
\be\lab{2014-4-22-bue1}
c(u,v)=\int_{\R_+^N} \frac{|u|^\alpha|v|^2}{|x|^{s_2}}dx\leq \frac{1}{\eta_1}\|v\|^2\;\hbox{for all}\;v\in D_{0}^{1,2}(\R_+^N).
\ee
Define
$\displaystyle k_0:=\frac{\eta_1}{2^*(s_2)}.$
Hence, when $\kappa<k_0$, we have
$$\|v\|^2-2^*(s_2)\kappa  c(u,v)>0\;\hbox{for all}\; v\in D_{0}^{1,2}(\R_+^N)\backslash\{0\},$$
which implies $(b)$.
For $\kappa>k_0=\frac{\eta_1}{2^*(s_2)}$, by the definition of $\eta_1$, there exists some $v\in D_{0}^{1,2}(\R_+^N)\backslash\{0\}$ such that
$\displaystyle \|v\|^2-\kappa 2^*(s_2) c(u,v)<0,$
and then it follows $(c)$.
\ep

Summarize the above conclusions, we obtain the following result:
\bc\lab{2014-11-24-xcro1}
Assume that $ 1<\alpha, 1<\beta, \alpha+\beta\leq 2^*(s_2)$ and one of the following holds:
\begin{itemize}
\item[(i)]$\kappa>0$.
\item[(ii)] $\kappa<0$ and $s_2>s_1$.
\item[(iii)] $s_2=s_1$  and $ \kappa<0$ with $|\kappa|$ small enough.
\end{itemize}
Let $U_\lambda:=\big(\frac{\mu_{s_1}(\R_+^N)}{\lambda}\big)^{\frac{1}{2^*(s_1)-2}}U$, where $U$ is a ground state solution of \eqref{2014-11-23-we2}.
\begin{itemize}
\item[(1)] Assume that either  $\kappa<0$ or $\beta>2$ or $\beta=2$ but with $\kappa<\frac{\eta_1}{2^*(s_2)}$, then $(U_\lambda,0)$ is a local minimum point of $\Phi$ in $\mathcal{N}$.
\item[(2)] Assume  that either $\beta<2$ and $\kappa>0$ or $\beta=2$ but with $\kappa>\frac{\eta_1}{2^*(s_2)}$, then
$$ c_0:=\inf_{(\phi,\varphi)\in \mathcal{N}}\Phi(\phi,\varphi)<\Phi(U_\lambda,0)=\Psi_{\lambda}(U_\lambda)=m_\lambda,$$
\end{itemize}
where $\eta_1$ is defined as that in Lemma \ref{2014-11-24-l3}. In particular, $\eta_1=\lambda$ if $s_1=s_2=s$.
\ec
Apply the similar arguments, we can obtain the following result:
\bc\lab{2014-4-12-Cro4}
Assume that $ 1<\alpha, 1<\beta, \alpha+\beta\leq 2^*(s_2)$ and one of the following holds:
\begin{itemize}
\item[(i)]$\kappa>0$.
\item[(ii)] $\kappa<0$ and $s_2>s_1$.
\item[(iii)] $s_2=s_1$  and $ \kappa<0$ with $|\kappa|$ small enough.
\end{itemize}
Let $U_\mu:=\big(\frac{\mu_{s_1}(\R_+^N)}{\mu}\big)^{\frac{1}{2^*(s_1)-2}}U$, where $U$ is a ground state solution of \eqref{2014-11-23-we2}.  Then
\begin{itemize}
\item[(1)] if  either $\kappa<0$ or $\alpha>2$ or $\alpha=2$ but with $\kappa<\frac{\eta_2}{2^*(s_2)}$,  then $(0,U_\mu)$ is a local minimum point of $\Phi$ in $\mathcal{N}$;
\item[(2)] if $\alpha<2$ and $\kappa>0$ or $\alpha=2$ but with $\kappa>\frac{\eta_2}{2^*(s_2)}$, then
$$ c_0=\inf_{(\phi,\varphi)\in \mathcal{N}}\Phi(\phi,\varphi)<\Phi(0,U_\mu)=\Psi_{\mu}(U_\mu)=m_\mu,$$
\end{itemize}
where
$\displaystyle \eta_2:=\inf_{u\in \Theta}\|u\|^2$, and
$\displaystyle\Theta:=\Big\{u\in D_{0}^{1,2}(\R_+^N): \;\; \int_{\R_+^N}\frac{|u|^2 |U_\mu|^\beta}{|x|^{s_2}}dx=1\Big\}.$ In particular $\eta_2=\mu$ whenever  $s_1=s_2=s$.
\ec

\s{The case of $s_1=s_2=s\in (0,2)$: nontrivial ground state and uniqueness;   sharp constant  $S_{\alpha,\beta,\lambda,\mu}(\Omega)$;   existence of infinitely many sign-changing solutions}
\renewcommand{\theequation}{7.\arabic{equation}}
\renewcommand{\theremark}{7.\arabic{remark}}
\renewcommand{\thedefinition}{7.\arabic{definition}}

  In this section, we focus on the case of $s_1=s_2:=s\in (0,2)$;  the case $s_1\neq s_2$ will be studied in the forthcoming paper (Part II). That is, we study the following problem
\be\lab{2015-1-5-e1}
\begin{cases}
-\Delta u-\lambda \frac{|u|^{2^*(s)-2}u}{|x|^{s}}=\kappa\alpha \frac{1}{|x|^{s}}|u|^{\alpha-2}u|v|^\beta\quad &\hbox{in}\;\Omega,\\
-\Delta v-\mu \frac{|v|^{2^*(s)-2}v}{|x|^{s}}=\kappa\beta \frac{1}{|x|^{s}}|u|^{\alpha}|v|^{\beta-2}v\quad &\hbox{in}\;\Omega,\\
(u,v)\in \mathscr{D}:=D_{0}^{1,2}(\Omega)\times D_{0}^{1,2}(\Omega),
\end{cases}
\ee
where  $\Omega$ is  a cone in $\R^N$(especially, $\Omega=\R^N$ or $\Omega=\R_+^N$) or $\Omega\subset \R^N$ is  an open   domain but $0\not\in \bar{\Omega}$.  In this section, we are aim to study the existence of nontrivial ground state solution to the system \ref{2015-1-5-e1}. Thanks to the fact of that $s_1=s_2$, we shall  obtain further results: the uniqueness of the nontrivial ground state solution; the relationship between the sharp constant  $S_{\alpha,\beta,\lambda,\mu}(\Omega)$ and the domain $\Omega$; the existence of infinitely many sign-changing solutions to  system \ref{2015-1-5-e1}.   Finally, we will explore some approaches for studying  the sharp constant $S_{\alpha,\beta,\lambda,\mu}(\Omega)$ when
$\Omega$ is not necessarily a cone.

 Noting that for  the special case $s_1=s_2=s\in (0,2)$ and $\alpha+\beta=2^*(s)$, the nonlinearities are homogeneous which enable us to define the following constant
\be\lab{2014-11-26-e0}
S_{\alpha,\beta,\lambda,\mu}(\Omega):=\inf_{(u,v)\in \widetilde{\mathscr{D}}} \frac{\int_\Omega \big(|\nabla u|^2+|\nabla v|^2\big)dx}{\Big(\int_\Omega \big(\lambda \frac{|u|^{2^*(s)}}{|x|^s}+\mu \frac{|v|^{2^*(s)}}{|x|^s}+2^*(s)\kappa \frac{|u|^\alpha |v|^\beta}{|x|^s}\big)dx\Big)^{\frac{2}{2^*(s)}}},
\ee
where
\be\lab{2014-11-26-xe1}
\widetilde{\mathscr{D}}:=\left\{(u,v)\in \mathscr{D}:\;\int_\Omega \big(\lambda \frac{|u|^{2^*(s)}}{|x|^s}+\mu \frac{|v|^{2^*(s)}}{|x|^s}+2^*(s)\kappa \frac{|u|^\alpha |v|^\beta}{|x|^s}\big)dx>0\right\}.
\ee
The above constant determines the following kind of  inequalities:

\begin{align}\lab{zz-06-13}
  & S_{\alpha,\beta,\lambda,\mu}(\Omega)  \Big(\int_\Omega \big(\lambda \frac{|u|^{2^*(s)}}{|x|^s}+\mu \frac{|v|^{2^*(s)}}{|x|^s}+2^*(s)\kappa \frac{|u|^\alpha |v|^\beta}{|x|^s}\big)dx\Big)^{\frac{2}{2^*(s)}}\\
&\quad\quad\quad\quad\quad\quad\quad\quad\quad    \leq  \int_\Omega \big(|\nabla u|^2+|\nabla v|^2\big)dx\nonumber
\end{align}
for $(u,v)\in  {\mathscr{D}}$ or $\widetilde{\mathscr{D}}$, which can be viewed as the double-variable CKN inequality. To the best of our knowledge,
such kind of inequality and its sharp constant with extremal functions  have  not been studied before.

\vskip0.23in

Denote
\be\lab{2014-11-26-xe2}
\mu_s(\Omega):=\inf\left\{\int_{\Omega} |\nabla u|^2 dx\;:\;\int_{\Omega}\frac{|u|^{2^*(s)}}{|x|^s}dx=1\right\},
\ee
then $\mu_s({\Omega})$  can be  attained  when  $\Omega$ is  a cone in $\R^N$  and $0<s<2$  (see  \cite{Egnell.1992}).
Noting that $\big(D_{0}^{1,2}(\Omega)\backslash\{0\}\big)\times\{0\}\subset \widetilde{\mathscr{D}}$ and $\{0\}\times \big(D_{0}^{1,2}(\Omega)\backslash\{0\}\big)\subset \widetilde{\mathscr{D}}$, by the definition, we  have that
\be\lab{2014-11-26-xe3}
S_{\alpha,\beta,\lambda,\mu}(\Omega)\leq \big(\max\{\lambda, \mu\}\big)^{-\frac{2}{2^*(s)}}\mu_s(\Omega).
\ee
By Young's inequality, it is easy to see that $S_{\alpha,\beta,\lambda,\mu}(\Omega)>0$.  The following
statement is obvious.

\bo\lab{zzz=111}
Assume that $(u,v)$ is an  extremal of $S_{\alpha,\beta,\lambda,\mu}(\Omega)$ such that $$\int_\Omega \frac{1}{|x|^s}\big[\lambda |u|^{2^*(s)}+\mu|v|^{2^*(s)}+2^*(s)\kappa |u|^\alpha|v|^\beta\big]dx=1$$
and
$$ \int_\Omega \big(|\nabla u|^2+|\nabla v|^2\big)dx=S_{\alpha,\beta,\lambda,\mu}(\Omega).$$
Then    $$(\phi,\psi):=\Big((S_{\alpha,\beta,\lambda,\mu}(\Omega))^{\frac{1}{2^*(s)-2}}u, \;\; (S_{\alpha,\beta,\lambda,\mu}(\Omega))^{\frac{1}{2^*(s)-2}}v\Big)$$ is a ground state solution to problem \eqref{2015-1-5-e1} and the corresponding energy
$$ c_0=\Phi(\phi,\psi)=\big(\frac{1}{2}-\frac{1}{2^*(s)}\big)\big[S_{\alpha,\beta,\lambda,\mu}(\Omega)\big]^{\frac{2^*(s)}{2^*(s)-2}}.$$
\eo

\bl\lab{2014-11-26-l1}
Assume that $\alpha>0,\beta>0,\lambda>0,\mu>0$, then there exists a best constant
\be\lab{2014-11-26-e1}
\kappa(\alpha,\beta,\lambda,\mu)=(\alpha+\beta)\left(\frac{\lambda}{\alpha}\right)^{\frac{\alpha}{\alpha+\beta}}
(\frac{\mu}{\beta})^{\frac{\beta}{\alpha+\beta}}
\ee
such that
$$\kappa(\alpha,\beta,\lambda,\mu)\int_\Omega |u|^\alpha |v|^\beta d\nu\leq \lambda\int_\Omega |u|^{\alpha+\beta}d\nu+\mu \int_\Omega |v|^{\alpha+\beta}d\nu$$ for all $(u,v)\in L^{\alpha+\beta}(\Omega,d\nu)\times L^{\alpha+\beta}(\Omega,d\nu)$.
\el
\bp
By Young's inequality with $\varepsilon$,
\be\lab{2014-11-26-e2}
xy\leq \varepsilon x^p+C(\varepsilon) y^q\quad (x, y>0, \varepsilon>0)
\ee
where $\frac{1}{p}+\frac{1}{q}=1, C(\varepsilon)=(\varepsilon p)^{-q/p}q^{-1}$. Take $p=\frac{\alpha+\beta}{\alpha}, q=\frac{\alpha+\beta}{\beta}$ and $\displaystyle x=|u|^\alpha, y=|v|^\beta, \displaystyle\varepsilon=\frac{1}{\alpha+\beta}
(\frac{\lambda}{\mu})^{\frac{\beta}{\alpha+\beta}}
\alpha^{\frac{\alpha}{\alpha+\beta}}\beta^{\frac{\beta}{\alpha+\beta}},$
then we obtain that
\be\lab{2014-11-26-e2}
|u|^\alpha|v|^\beta\leq \frac{1}{\kappa(\alpha,\beta,\lambda,\mu)}\big(\lambda |u|^{\alpha+\beta}+\mu|v|^{\alpha+\beta}\big).
\ee
Hence, for all $(u,v)\in  L^{\alpha+\beta}(\Omega,d\nu)\times L^{\alpha+\beta}(\Omega,d\nu)$, we have
\be\lab{2014-11-26-e3}
\kappa(\alpha,\beta,\lambda,\mu)\int_\Omega |u|^\alpha |v|^\beta d\nu\leq \lambda\int_\Omega |u|^{\alpha+\beta}d\nu+\mu \int_\Omega |v|^{\alpha+\beta}d\nu.
\ee
And we note that when $(u, tu)$ with
$\displaystyle t=(\frac{\lambda \beta}{\mu\alpha})^{\frac{1}{\alpha+\beta}},$
the constant $\kappa(\alpha,\beta,\lambda,\mu)$ is attained.
Hence, $\kappa(\alpha,\beta,\lambda,\mu)$ is the best constant.
\ep

\bl\lab{2014-11-26-l3} Assume  $\kappa\leq 0$. Then
$$S_{\alpha,\beta,\lambda,\mu}(\Omega)=\big(\max\{\lambda, \mu\}\big)^{-\frac{2}{2^*(s)}}\mu_s(\Omega).$$
In particular, in this case $S_{\alpha,\beta,\lambda,\mu}(\Omega)$   can only  be attained by the
semi-trivial pairs: $$
\begin{cases}
(U,0)\quad &\hbox{if}\;\lambda>\mu;\\
(0,U)\quad &\hbox{if}\;\lambda<\mu;\\
(U,0)\;\hbox{or}\;(0,U)\;&\hbox{if}\;\lambda=\mu;
\end{cases}
$$
where $U$ is an extremal function of $\mu_s(\Omega)$.  Hence, the  ground state   to problem \eqref{2015-1-5-e1}
can only be attained by semi-trivial pairs.

\el
\bp
By \eqref{2014-11-26-xe3}, we only need to prove the reverse inequality.
Indeed, for any $(u,v)\in \tilde{\mathscr{D}}$, we have
\begin{align}\lab{2014-11-26-xe4}
&\frac{\int_\Omega \big(|\nabla u|^2+|\nabla v|^2\big)dx}{\Big(\int_\Omega \big(\lambda \frac{|u|^{2^*(s)}}{|x|^s}+\mu \frac{|v|^{2^*(s)}}{|x|^s}+2^*(s)\kappa \frac{|u|^\alpha |v|^\beta}{|x|^s}\big)dx\Big)^{\frac{2}{2^*(s)}}}\\
\geq& \frac{\int_\Omega \big(|\nabla u|^2+|\nabla v|^2\big)dx}{\Big(\int_\Omega \big(\lambda \frac{|u|^{2^*(s)}}{|x|^s}+\mu \frac{|v|^{2^*(s)}}{|x|^s}\Big)^{\frac{2}{2^*(s)}}}\nonumber\\
\geq&\big(\max\{\lambda, \mu\}\big)^{-\frac{2}{2^*(s)}}\frac{\int_\Omega \big(|\nabla u|^2+|\nabla v|^2\big)dx}{\Big(\int_\Omega \big( \frac{|u|^{2^*(s)}}{|x|^s}+ \frac{|v|^{2^*(s)}}{|x|^s}\big)dx\Big)^{\frac{2}{2^*(s)}}}\nonumber\\
\geq&\big(\max\{\lambda, \mu\}\big)^{-\frac{2}{2^*(s)}}\mu_s(\Omega)
\frac{\big(\int_\Omega \frac{|u|^{2^*(s)}}{|x|^s}dx\big)^{\frac{2}{2^*(s)}}+\big(\int_\Omega \frac{|v|^{2^*(s)}}{|x|^s}dx\big)^{\frac{2}{2^*(s)}}}{\big(\int_\Omega \frac{|u|^{2^*(s)}}{|x|^s}+\frac{|v|^{2^*(s)}}{|x|^s}dx\big)^{\frac{2}{2^*(s)}} }\nonumber\\
\geq&\big(\max\{\lambda, \mu\}\big)^{-\frac{2}{2^*(s)}}\mu_s(\Omega).\nonumber
\end{align}
By taking the infimum over $\tilde{\mathscr{D}}$, we obtain that
$\displaystyle S_{\alpha,\beta,\lambda,\mu}(\Omega)\geq \big(\max\{\lambda, \mu\}\big)^{-\frac{2}{2^*(s)}}\mu_s(\Omega)$.
Moreover, by the processes of \eqref{2014-11-26-xe4}, we see that
$S_{\alpha,\beta,\lambda,\mu}(\Omega)$ is only achieved by
$(U,0)$ if $\lambda>\mu$ (resp. $(0,U)$ if $\lambda<\mu$; $(U,0)$ or $(0,U)$ if $\lambda=\mu$),
where $U$ is an extremal function of $\mu_s(\Omega)$.
\ep

\bo\lab{2014-11-26-l2}
$\displaystyle \widetilde{\mathscr{D}}=\mathscr{D}\backslash \{(0,0)\}$  when $\displaystyle \kappa>-(\frac{\lambda}{\alpha})^{\frac{\alpha}{2^*(s)}}
(\frac{\mu}{\beta})^{\frac{\beta}{2^*(s)}}$.
\eo
\bp
By Lemma \ref{2014-11-26-l1}, if $\displaystyle \kappa>-(\frac{\lambda}{\alpha})^{\frac{\alpha}{2^*(s)}}
(\frac{\mu}{\beta})^{\frac{\beta}{2^*(s)}}$, then there exists some $C>0$ such that
\be\lab{2014-11-26-e4}
\int_\Omega\big(\lambda \frac{|u|^{2^*(s)}}{|x|^s}+\mu \frac{|v|^{2^*(s)}}{|x|^s}+2^*(s)\kappa \frac{|u|^\alpha |v|^\beta}{|x|^s}\big)dx\geq C\Big(\int_\Omega(\lambda \frac{|u|^{2^*(s)}}{|x|^s}+\mu \frac{|v|^{2^*(s)}}{|x|^s}\big)dx\Big).
\ee
Thereby this proposition  is proved.
\ep

\br\lab{zz-0613=1}  By the Lemma \ref{2014-11-26-l3} and Proposition \ref{2014-11-26-l2}, we know that,  when $0>\kappa>-(\frac{\lambda}{\alpha})^{\frac{\alpha}{2^*(s)}} (\frac{\mu}{\beta})^{\frac{\beta}{2^*(s)}}$, then
the inequality (\ref{zz-06-13}) is meaningful but the sharp constant $S_{\alpha,\beta,\lambda,\mu}(\Omega)$ can be
reached only by semi-trivial extremals.
\er

\noindent{\bf Note:} In view of   Lemma \ref{2014-11-26-l3}, Proposition \ref{2014-11-26-l2} and Remark  \ref{zz-0613=1},
we have  to consider the case  that $\kappa>0$.     Therefore, throughout the remaining part of the current paper, we assume
$\kappa>0$.

We obtain the following result:

\bt\lab{2014-11-26-th1}
Assume  $\Omega$ is a cone in $\R^N$(especially, $\Omega=\R^N$ or $\Omega=\R_+^N$) or $\Omega$ is an open domain with $0\not\in \bar{\Omega}$. If  $0<s<2,\alpha>1,\beta>1,\alpha+\beta=2^*(s)$ and
 $\kappa>0$, then $S_{\alpha,\beta,\lambda,\mu}(\Omega)$ is always achieved by some nonnegative function $(u,v)$.
\et

\br   Theorem \ref{2014-11-26-th1} asserts that the constant $S_{\alpha,\beta,\lambda,\mu}(\Omega)$ can be attained by a
 nonnegative extremal function. But at this moment, we can not exclude the possibility of the semi-triviality of the extremal
 function.   Fortunately, we will see further results in Theorem \ref{2014-11-27-wth1} below for the nontrivial  extremal
 function.

\er


When we  study the critical problem of  a scalar  equation, say $(P)$, blow-up analysis is one of the classical method.   Its ideas  usually are
as follows:  Consider a modified subcritical problem define on a convex domain, say $(P_\varepsilon)$, which approximates the problem $(P)$. Then study the existence of positive solution $u_\varepsilon$  to the modified problem $(P_\varepsilon)$ and the regularity of the positive solutions. By the standard Pohozaev identity, $u_\varepsilon$ must blow-up as $\varepsilon\rightarrow 0$.  Apply the standard blow-up arguments to obtain the existence of positive solution of $(P)$.  It follows that the set of solutions  is  nonempty. Therefore, one may  take  a minimizing sequence to
 approach the ground state. In this process,   the approximation involves the domains and the format of the nonlinearities.  However, for the
 system \ref{2015-1-5-e1}, the customary skills can not be applied directly. Since we want to get rid of the semi-trivial solution, the main obstacles lies in that we can not get the precise estimate about  the limit of the least energy of the approximate problem.  Therefore,
  in the current paper, we introduce a new approximation system  where we just modify the singularity.

\subsection{Approximating the  problems}
  For any $\varepsilon\geq 0$, set
\be\lab{2014-10-8-e1}
a_\varepsilon(x):=\begin{cases}
\frac{1}{|x|^{s-\varepsilon}}\quad &\hbox{for}\;|x|<1,\\
\frac{1}{|x|^{s+\varepsilon}}&\hbox{for}\;|x|\geq 1.
\end{cases}
\ee
\bl\lab{2014-11-26-xl1}
Let $\varepsilon>0$. Then
for any $u\in D_{0}^{1,2}(\Omega)$, $\int_\Omega a_\varepsilon(x)|u|^{2^*(s)}dx$ is well defined and decreasing by $\varepsilon$.
\el
\bp
Let $\varepsilon_1>\varepsilon_2\geq 0$. By the definition of $a_\varepsilon(x)$, it is easy to obtain the result by noting that $a_{\varepsilon_1}(x)<a_{\varepsilon_2}(x)\leq a_0(x)$.
\ep

We also note that for any compact set $\Omega_1\subset \Omega$ such that $0\not\in\bar{\Omega}_1$, $a_\varepsilon(x)\rightarrow a_0(x)$ uniformly on $\Omega_1$ as $\varepsilon\rightarrow 0$. For any fixed $\varepsilon>0$, we consider the ground state solution to the following problem:
\be\lab{2014-11-26-xbue2}
\begin{cases}
-\Delta u-\lambda a_\varepsilon(x)|u|^{2^*(s)-2}u=\kappa\alpha a_\varepsilon(x)|u|^{\alpha-2}u|v|^\beta\quad &\hbox{in}\;\Omega,\\
-\Delta v-\mu a_\varepsilon(x)|v|^{2^*(s)-2}v=\kappa\beta a_\varepsilon(x)|u|^{\alpha}|v|^{\beta-2}v\quad &\hbox{in}\;\Omega,\\
\kappa>0,(u,v)\in \mathscr{D}:=D_{0}^{1,2}(\Omega)\times D_{0}^{1,2}(\Omega),
\end{cases}
\ee
Consider the following variational problem
\be\lab{2014-11-26-xbue3}
\min \int_\Omega \big(|\nabla u|^2+|\nabla v|^2\big)dx
\ee
\be\lab{2014-11-26-xbue3+0}
\hbox{s.t. }\;\int_\Omega a_\varepsilon(x)\big[\lambda |u|^{2^*(s)}+\mu|v|^{2^*(s)}+2^*(s)\kappa |u|^\alpha|v|^\beta\big]dx=1.\ee
Let $S_{\alpha,\beta,\lambda,\mu}^{\varepsilon}(\Omega)$ be the minimize value of \eqref{2014-11-26-xbue3}, then we have:

\bl\lab{2014-11-26-xl2} The constant
$S_{\alpha,\beta,\lambda,\mu}^{\varepsilon}(\Omega)$ is increasing with  respect to $\varepsilon$ and $$\displaystyle \lim_{\varepsilon\rightarrow 0^+}S_{\alpha,\beta,\lambda,\mu}^{\varepsilon}(\Omega)=S_{\alpha,\beta,\lambda,\mu}(\Omega).$$
\el
\bp
By the definition of $S_{\alpha,\beta,\lambda,\mu}^{\varepsilon}(\Omega)$ and $a_\varepsilon(x)$, it is easy to see that
$$S_{\alpha,\beta,\lambda,\mu}^{\varepsilon_1}(\Omega)\geq S_{\alpha,\beta,\lambda,\mu}^{\varepsilon_2}(\Omega)\;\hbox{for any}\;\varepsilon_1>\varepsilon_2\geq 0.$$
Hence, we have
\be\lab{2014-11-26-xbue4}
\liminf_{\varepsilon\rightarrow 0^+}S_{\alpha,\beta,\lambda,\mu}^{\varepsilon}(\Omega)\geq S_{\alpha,\beta,\lambda,\mu}(\Omega).
\ee
On the other hand, for any $\eta>0$, there exists $(u,v)\in C_c^\infty(\Omega)\times C_c^\infty(\Omega)$ such that
$$\int_\Omega a_0(x)\big[\lambda |u|^{2^*(s)}+\mu|v|^{2^*(s)}+2^*(s)\kappa |u|^\alpha|v|^\beta\big]dx=1$$
and
$$ \int_\Omega \big(|\nabla u|^2+|\nabla v|^2\big)dx<S_{\alpha,\beta,\lambda,\mu}(\Omega)+\eta.$$
Since $a_\varepsilon(x)\rightarrow a_0(x)$ in $L^\infty(supp(u))$ as $\varepsilon\rightarrow 0$, we obtain that
\begin{align*}
&\limsup_{\varepsilon\rightarrow 0^+} S_{\alpha,\beta,\lambda,\mu}^{\varepsilon}(\Omega)\\
=&\limsup_{\varepsilon\rightarrow 0^+}\frac{\int_\Omega \big(|\nabla u|^2+|\nabla v|^2\big)dx}{\Big(\int_\Omega a_\varepsilon(x)\big[\lambda |u|^{2^*(s)}+\mu|v|^{2^*(s)}+2^*(s)\kappa |u|^\alpha|v|^\beta\big]dx\Big)^{\frac{2}{2^*(s)}}}\\
=&\frac{\int_\Omega \big(|\nabla u|^2+|\nabla v|^2\big)dx}{\Big(\int_\Omega a_0(x)\big[\lambda |u|^{2^*(s)}+\mu|v|^{2^*(s)}+2^*(s)\kappa |u|^\alpha|v|^\beta\big]dx\Big)^{\frac{2}{2^*(s)}}}\\
=&\int_\Omega \big(|\nabla u|^2+|\nabla v|^2\big)dx
<S_{\alpha,\beta,\lambda,\mu}(\Omega)+\eta.
\end{align*}
By the arbitrariness of $\eta$, we have
\be\lab{2014-11-26-xbue5}
\limsup_{\varepsilon\rightarrow 0^+}S_{\alpha,\beta,\lambda,\mu}^{\varepsilon}(\Omega)\leq S_{\alpha,\beta,\lambda,\mu}(\Omega).
\ee
Thus the proof is completed by \eqref{2014-11-26-xbue4} and \eqref{2014-11-26-xbue5}.
\ep

Let $L^{p}(\Omega,a_\varepsilon(x)dx)$ denote  the space of $L^{p}$-integrable functions with respect to the measure $a_\varepsilon(x)dx$ and the corresponding norm is indicated by $$|u|_{p, \varepsilon}:=\ \big(\int_\Omega a_\varepsilon(x)|u|^pdx\big)^{\frac{1}{p}},\quad p>1.$$
Then we have the following compact embedding result:
\bl\lab{2014-11-26-xl3}
For any $\varepsilon\in (0,s)$, the embedding $D_{0}^{1,2}(\Omega)\hookrightarrow L^{2^*(s)}(\Omega,a_\varepsilon(x)dx)$ is compact.
\el
\bp
Let $\{u_n\}\subset D_{0}^{1,2}(\Omega)$ be a bounded sequence. Up to a subsequence, we may assume that $u_n\rightharpoonup u$ in $D_{0}^{1,2}(\Omega)$ and $u_n\rightarrow u$ a.e. in $\Omega$.
Then for any $R>1$, we have
\be\lab{2014-11-26-xbue6}
\int_{\Omega\cap B_R^c}a_\varepsilon(x)|u_n|^{2^*(s)}dx\leq \frac{1}{R^\varepsilon}\int_{\Omega\cap B_R^c}a_0|u_n|^{2^*(s)}dx\rightarrow 0,
\ee
uniformly for all $n$ as $R\rightarrow +\infty$.

Noting that $2^*(s)<2^*(s-\varepsilon), 2^*(s)<2^*:=\frac{2N}{N-2}$,  by Rellich-Kondrachov compact theorem, we have
\be\lab{2014-11-26-xbue7}
\lim_{n\rightarrow \infty}\int_{\Omega\cap B_R}a_\varepsilon(x)|u_n-u|^{2^*(s)}dx=0.
\ee
By \eqref{2014-11-26-xbue6} and \eqref{2014-11-26-xbue7}, we prove this Lemma.
\ep

\bl\lab{2014-11-26-xl4}
For any $\varepsilon\in (0,s)$, $S_{\alpha,\beta,\lambda,\mu}^{\varepsilon}(\Omega)$ is attained by some extremal  $(u_\varepsilon, v_\varepsilon)$, i.e.,
$$\int_\Omega a_\varepsilon(x)\big[\lambda |u_\varepsilon|^{2^*(s)}+\mu|v_\varepsilon|^{2^*(s)}+2^*(s)\kappa |u_\varepsilon|^\alpha|v_\varepsilon|^\beta\big]dx=1$$
and
$$ \int_\Omega \big(|\nabla u_\varepsilon|^2+|\nabla v_\varepsilon|^2\big)dx=S_{\alpha,\beta,\lambda,\mu}^{\varepsilon}(\Omega).$$
Moreover, $(u_\varepsilon, v_\varepsilon)$ satisfies the following equation:
\be\lab{2014-11-26-xbue8}
\begin{cases}
-\Delta u=S_{\alpha,\beta,\lambda,\mu}^{\varepsilon}(\Omega)\Big(\lambda a_\varepsilon(x)|u|^{2^*(s)-2}u+\kappa\alpha a_\varepsilon(x)|u|^{\alpha-2}u|v|^\beta\Big)\quad &\hbox{in}\;\Omega,\\
-\Delta v=S_{\alpha,\beta,\lambda,\mu}^{\varepsilon}(\Omega)\Big(\mu a_\varepsilon(x)|v|^{2^*(s)-2}v+\kappa\beta a_\varepsilon(x)|u|^{\alpha}|v|^{\beta-2}v\Big)\quad &\hbox{in}\;\Omega,\\
u\geq 0, v\geq 0, (u,v)\in \mathscr{D}:=D_{0}^{1,2}(\Omega)\times D_{0}^{1,2}(\Omega),
\end{cases}
\ee
\el
\bp
Let $\{u_{n,\varepsilon}\}\subset \mathscr{D}$  be a minimizing sequence, i.e., $$\int_\Omega a_\varepsilon(x)\big[\lambda |u_{n,\varepsilon}|^{2^*(s)}+\mu|v_{n,\varepsilon}|^{2^*(s)}+2^*(s)\kappa |u_{n,\varepsilon}|^\alpha|v_{n,\varepsilon}|^\beta\big]dx=1$$
and
$$\int_\Omega \big(|\nabla u_{n,\varepsilon}|^2+|\nabla u_{n,\varepsilon}|^2\big)dx\rightarrow S_{\alpha,\beta,\lambda,\mu}^{\varepsilon}(\Omega)\;\hbox{as}\;n\rightarrow +\infty.$$
 Then we see that $\{u_{n,\varepsilon}\}$ and $\{v_{n,\varepsilon}\}$ are bounded in $D_{0}^{1,2}(\Omega)$. By Lemma \ref{2014-11-26-xl3}, we see that up to a subsequence $u_{n,\varepsilon}\rightarrow u_\varepsilon$ and $v_{n,\varepsilon}\rightarrow v_\varepsilon$  in $L^{2^*(s)}(\Omega,a_\varepsilon(x)dx)$. Hence
 $$\int_\Omega a_\varepsilon(x)\big[\lambda |u_\varepsilon|^{2^*(s)}+\mu|v_\varepsilon|^{2^*(s)}+2^*(s)\kappa |u_\varepsilon|^\alpha|v_\varepsilon|^\beta\big]dx=1.$$
 Then by the definition of $S_{\alpha,\beta,\lambda,\mu}^{\varepsilon}(\Omega)$, we have
$$\int_\Omega \big(|\nabla u_\varepsilon|^2+|\nabla v_\varepsilon|^2\big)dx\geq S_{\alpha,\beta,\lambda,\mu}^{\varepsilon}(\Omega).$$
On the other hand, by  the weak semi-continuous of a norm (or Fatou's Lemma), we have
$$\int_\Omega \big(|\nabla u_\varepsilon|^2+|\nabla v_\varepsilon|^2\big)dx\leq \liminf_{n\rightarrow \infty} \int_\Omega \big(|\nabla u_{n,\varepsilon}|^2+|\nabla u_{n,\varepsilon}|^2\big)dx=S_{\alpha,\beta,\lambda,\mu}^{\varepsilon}(\Omega).$$
Hence, $(u_\varepsilon, v_\varepsilon)$ is a minimizer of $S_{\alpha,\beta,\lambda,\mu}^{\varepsilon}(\Omega)$. Without loss of generality, we may assume that $u_{n,\varepsilon}\geq 0$ and $v_{n,\varepsilon}\geq 0$ for all $n$. Then since $u_{n,\varepsilon}\rightarrow u_\varepsilon, v_{n,\varepsilon}\rightarrow v_\varepsilon$ a.e. in $\Omega$, we obtain that $u_\varepsilon\geq 0, v_\varepsilon\geq 0$ a.e. in $\Omega$.
There exists some Lagrange multiplier $\Lambda\in\R$ such that
$$\begin{cases}
-\Delta u_\varepsilon=\Lambda\Big(\lambda a_\varepsilon(x)|u_\varepsilon|^{2^*(s)-2}u_\varepsilon+\kappa\alpha a_\varepsilon(x)|u_\varepsilon|^{\alpha-2}u_\varepsilon|v_\varepsilon|^\beta\Big)\quad &\hbox{in}\;\Omega,\\
-\Delta v_\varepsilon=\Lambda\Big(\mu a_\varepsilon(x)|v_\varepsilon|^{2^*(s)-2}v_\varepsilon+\kappa\beta a_\varepsilon(x)|u_\varepsilon|^{\alpha}|v_\varepsilon|^{\beta-2}v_\varepsilon\Big)\quad &\hbox{in}\;\Omega.
\end{cases}$$
Testing by $(u_\varepsilon, v_\varepsilon)$, we obtain that
$\Lambda=S_{\alpha,\beta,\lambda,\mu}^{\varepsilon}(\Omega)$.
\ep

\bl\lab{2014-11-26-xl5}
For $\varepsilon\in (0,s)$, assume $(u_\varepsilon,v_\varepsilon)$ is a solution to \eqref{2014-11-26-xbue8} given by Lemma \ref{2014-11-26-xl4}. Then:
\begin{itemize}
\item[(i)] The family $(u_\varepsilon,v_\varepsilon)$ is bounded in $\mathscr{D}$;
\item[(ii)] Up to a subsequence, we have $u_\varepsilon\rightharpoonup u, v_\varepsilon\rightharpoonup v$ in $D_{0}^{1,2}(\Omega)$ and $u_\varepsilon\rightarrow u, v_\varepsilon\rightarrow v$ a.e.  in $\Omega$ as $\varepsilon\rightarrow 0$.
\item[(iii)] $(u,v)$ given in $(ii)$ weakly solves
\be\lab{2014-11-26-xbue9}
\begin{cases}
-\Delta u=S_{\alpha,\beta,\lambda,\mu}(\Omega)\Big(\lambda a_0(x)|u|^{2^*(s)-2}u+\kappa\alpha a_0(x)|u|^{\alpha-2}u|v|^\beta\Big)\quad &\hbox{in}\;\Omega,\\
-\Delta v=S_{\alpha,\beta,\lambda,\mu}(\Omega)\Big(\mu a_0(x)|v|^{2^*(s)-2}v+\kappa\beta a_0(x)|u|^{\alpha}|v|^{\beta-2}v\Big)\quad &\hbox{in}\;\Omega,\\
u\geq 0, v\geq 0, (u,v)\in \mathscr{D},
\end{cases}
\ee
\item[(iv)] If $(u,v)\neq (0,0)$, then
 $\displaystyle\int_\Omega a_0(x)\big[\lambda |u|^{2^*(s)}+\mu|v|^{2^*(s)}+2^*(s)\kappa |u|^\alpha|v|^\beta\big]dx=1$
 and $u_\varepsilon\rightarrow u$ and  $v_\varepsilon\rightarrow v$  strongly in $D_{0}^{1,2}(\Omega)$. Moreover, $(u,v)$ is an extremal function of $S_{\alpha,\beta,\lambda,\mu}(\Omega).$
\end{itemize}
\el
\bp
$(i)$ follows by Lemma \ref{2014-11-26-xl2} and
$(ii)$  is trivial;

$(iii)$  Without loss of generality, we assume that $\varepsilon_k\downarrow 0$ as $k\rightarrow \infty$. For any $\phi\in C_c^\infty(\Omega), \psi\in C_c^\infty(\Omega)$, since $(u_{\varepsilon_k}, v_{\varepsilon_k})$ is a solution of \eqref{2014-11-26-xbue8} with $\varepsilon=\varepsilon_k$, we have
\begin{align}\lab{2014-11-26-xbue10}
&\int_\Omega \big(\nabla u_{\varepsilon_k}\cdot \nabla \phi+\nabla v_{\varepsilon_k}\cdot \nabla \psi\big)dx\\
&=S_{\alpha,\beta,\lambda,\mu}^{\varepsilon_k}(\Omega)\int_\Omega \Big(\lambda a_\varepsilon(x)|u_\varepsilon|^{2^*(s)-2}u_\varepsilon \phi+\kappa\alpha a_\varepsilon(x)|u_\varepsilon|^{\alpha-2}u_\varepsilon \phi|v_\varepsilon|^\beta\Big)dx\nonumber\\
&\quad+S_{\alpha,\beta,\lambda,\mu}^{\varepsilon_k}(\Omega)\int_\Omega \Big(\mu a_\varepsilon(x)|v_\varepsilon|^{2^*(s)-2}v_\varepsilon\psi+\kappa\beta a_\varepsilon(x)|u_\varepsilon|^{\alpha}|v_\varepsilon|^{\beta-2}v_\varepsilon\psi\Big)dx.\nonumber
\end{align}
Recalling that $a_{\varepsilon_k}(x)\rightarrow a_0(x)$ in $L^\infty(sppt(\phi)\cup sppt(\psi))$ and $S_{\alpha,\beta,\lambda,\mu}^{\varepsilon_k}(\Omega)\rightarrow S_{\alpha,\beta,\lambda,\mu}(\Omega)$ as $k\rightarrow \infty$,   we obtain that
\begin{align}\lab{2014-11-26-xbue11}
&\int_\Omega \big(\nabla u\cdot \nabla \phi+\nabla v\cdot \nabla \psi\big)dx\\
&=S_{\alpha,\beta,\lambda,\mu}(\Omega)\int_\Omega \Big(\lambda a_0(x)|u|^{2^*(s)-2}u \phi+\kappa\alpha a_0(x)|u|^{\alpha-2}u \phi|v|^\beta\Big)dx\nonumber\\
&\quad +S_{\alpha,\beta,\lambda,\mu}(\Omega)\int_\Omega \Big(\mu a_0(x)|v|^{2^*(s)-2}v\psi+\kappa\beta a_0(x)|u|^{\alpha}|v|^{\beta-2}v\psi\Big)dx.\nonumber
\end{align}
Since  $(\phi,\psi)$ is arbitrary, we  see that $(u,v)$ weakly solve \eqref{2014-11-26-xbue9}.

$(iv)$ By Fatou's lemma, we have
\begin{align}\lab{2014-11-26-xbue12}
&\int_\Omega a_0(x)\big[\lambda |u|^{2^*(s)}+\mu|v|^{2^*(s)}+2^*(s)\kappa |u|^\alpha|v|^\beta\big]dx\\
\leq&\liminf_{k\rightarrow \infty}\int_\Omega a_{\varepsilon_k}(x)\big[\lambda |u_{\varepsilon_k}|^{2^*(s)}+\mu|v_{\varepsilon_k}|^{2^*(s)}+2^*(s)\kappa |u_{\varepsilon_k}|^\alpha|v_{\varepsilon_k}|^\beta\big]dx=1.\nonumber
\end{align}
If $\int_\Omega a_0(x)\big[\lambda |u|^{2^*(s)}+\mu|v|^{2^*(s)}+2^*(s)\kappa |u|^\alpha|v|^\beta\big]dx\neq 1$, since $(u,v)\neq (0,0)$, we have
$$ 0<\int_\Omega a_0(x)\big[\lambda |u|^{2^*(s)}+\mu|v|^{2^*(s)}+2^*(s)\kappa |u|^\alpha|v|^\beta\big]dx<1.$$
Hence, by $(iii)$ and the definition of $S_{\alpha,\beta,\lambda,\mu}(\Omega)$, we have
\begin{align}\lab{2014-11-26-xbue13}
S_{\alpha,\beta,\lambda,\mu}(\Omega)=&\frac{\int_\Omega \big(|\nabla u|^2+|\nabla v|^2\big)dx}{\int_\Omega a_0(x)\big[\lambda |u|^{2^*(s)}+\mu|v|^{2^*(s)}+2^*(s)\kappa |u|^\alpha|v|^\beta\big]dx}\\
>&\frac{\int_\Omega \big(|\nabla u|^2+|\nabla v|^2\big)dx}{\Big(\int_\Omega a_0(x)\big[\lambda |u|^{2^*(s)}+\mu|v|^{2^*(s)}+2^*(s)\kappa |u|^\alpha|v|^\beta\big]dx\Big)^{\frac{2}{2^*(s)}}}\nonumber\\
\geq&S_{\alpha,\beta,\lambda,\mu}(\Omega),\nonumber
\end{align}
a contradiction. Hence, $\int_\Omega a_0(x)\big[\lambda |u|^{2^*(s)}+\mu|v|^{2^*(s)}+2^*(s)\kappa |u|^\alpha|v|^\beta\big]dx=1$. It follows that
\be\lab{2014-11-26-xbue14}
\int_\Omega \big(|\nabla u_{\varepsilon_k}|^2+|\nabla v_{\varepsilon_k}|^2\big)\rightarrow \int_\Omega \big(|\nabla u|^2+|\nabla v|^2\big)dx,
\ee
which implies that $u_{\varepsilon_k}\rightarrow u, v_{\varepsilon_k}\rightarrow v$ in $D_{0}^{1,2}(\Omega)$.
\ep

\vskip0.3in

\subsection{Pohozaev Identity and  the Proof of Theorem \ref{2014-11-26-th1}}
\bo\lab{2013-10-26-prop1}
Let $(u,v)\in D_{0}^{1,2}(\Omega)\times D_{0}^{1,2}(\Omega)$ be a solution of the system
$$\begin{cases}-\Delta u&=G_u(x, u, v),\\-\Delta v&=G_v(x,u,v), \end{cases}$$ where $\displaystyle G(x, u,v)=\int_0^uG_s(x, s,v)ds+G(x,0,v)=\int_0^v G_t(x,u,t)dt+G(x,u,0)$ is such that $ G(x,0,0)\equiv 0, G(\cdot, u(\cdot), v(\cdot))$ and that $x_iG_{x_i}(\cdot, u(\cdot), v(\cdot))$ are in $L^1(\Omega)$. Then $(u,v)$ satisfies:
\begin{align}\lab{2013-10-26-e1}
&\int_{\partial\Omega}|\nabla(u,v)|^2x\cdot \nu d\sigma\\
=&2N\int_\Omega G(x, u, v)dx+2\sum_{i=1}^{N}\int_\Omega x_iG_{x_i}(x, u, v)dx-(N-2)\int_\Omega |\nabla (u,v)|^2dx,\nonumber
\end{align}
where $\Omega$ is a regular domain in $\R^N$, $\nu$ denotes the unitary exterior normal vector to $\partial\Omega$ and $|\nabla (u,v)|^2:=|\nabla u|^2+|\nabla v|^2$. Moreover, if $\Omega=\R^N$ or a cone, then
\be\lab{2013-10-26-e2}
2N\int_{\Omega}G(x, u, v)dx+2\sum_{i=1}^{N}\int_{\Omega}x_iG_{x_i}(x, u, v)dx=(N-2)\int_{\Omega}|\nabla (u,v)|^2dx.
\ee
\eo
\bp
Since $(u,v)$ is a solution, then we have
\be\lab{2013-10-26-e3}
0=\big(-\Delta u+G_u(x, u, v)\big)x\cdot \nabla u=\big(-\Delta v+G_v(x, u, v)\big)x\cdot \nabla v.
\ee
It is clear that
$$-\Delta wx\cdot \nabla w=-div(\nabla wx\cdot\nabla w- x\frac{|\nabla w|^2}{2})-\frac{N-2}{2}|\nabla w|^2,\quad (w=u,v)$$
\begin{align*}&G_u(x, u, v)x\cdot\nabla u+G_v(x, u, v)x\cdot \nabla v\\
=&div(xG(x,u,v))-NG(x,u,v)-\sum_{i=1}^{N} x_iG_{x_i}(x,u,v).
\end{align*}
Integrating by parts, we obtain
\begin{align}\lab{2013-10-26-e4}
&\int_{\partial\Omega}\big(\sigma G(\sigma, u, v)+(\nabla u\sigma\cdot\nabla u- \sigma\frac{|\nabla u|^2}{2})+(\nabla v\sigma\cdot\nabla v-\sigma\frac{|\nabla v|^2}{2})\big)\cdot\nu d\sigma\\
=&\int_\Omega \Big( NG(x,u,v)-\frac{N-2}{2}|\nabla (u,v)|^2+\sum_{i=1}^{N}x_iG_{x_i}(x,u,v)\Big)dx.\nonumber
\end{align}
When $u=v=0$ on $\partial\Omega$, we have
\be\lab{2013-10-26-e5}
\nabla u=\nabla u\cdot \nu \nu, \nabla v=\nabla v\cdot \nu \nu.
\ee
Then by (\ref{2013-10-26-e4}), it follows that
\begin{align}\lab{2013-10-26-e6}
&\int_{\partial\Omega}\big(G(\sigma, u, v)+\frac{1}{2}|\nabla (u,v)|^2\big)\sigma\cdot\nu d\sigma\\
=&\int_\Omega \big( NG(x,u,v)-\frac{N-2}{2}|\nabla(u,v)|^2+\sum_{i=1}^{N}x_iG_{x_i}(x,u,v)\big)dx.\nonumber
\end{align}
Moreover, since $G(x, 0,0)\equiv 0$, if $u=v=0$ on $\partial \Omega$, we have $G(\sigma, u,v)\equiv0$ on $\partial\Omega$, then we obtain that
\be\lab{2013-10-26-e7}
\frac{1}{2}\int_{\partial\Omega}|\nabla (u,v)|^2\sigma\cdot\nu d\sigma
=\int_\Omega \big( NG(x,u,v)-\frac{N-2}{2}|\nabla(u,v)|^2+\sum_{i=1}^{N}x_iG_{x_i}(x,u,v)\big)dx,
\ee
which is equivalent to  (\ref{2013-10-26-e1}).  Using polar coordinate transformation, since $|\nabla (u,v)|\in L^2(\R^N)$, we have
\begin{align*}
&\int_{\R^N} |\nabla(u,v)|^2dx=\int_{\R^N} |\nabla u|^2+|\nabla v|^2dx\\
=&\int_0^\infty \int_{\partial B_r(0)}|\nabla u(r,\theta)|^2+|\nabla v(r,\theta)|^2d\theta r^{N-1}dr
=\int_0^\infty \zeta(r)dr<\infty,
\end{align*}
where $\zeta(r):=\int_{\partial B_r(0)}|\nabla u(r,\theta)|^2+|\nabla v(r,\theta)|^2d\theta r^{N-1}\geq 0$, then by the absolute continuity, there exists $r_n\rightarrow \infty$ such that $\zeta(r_n)\rightarrow 0$. Since $N\geq 3$, we have
$$ \int_{\partial B_{r_n}(0)}\big(|\nabla u(r_n,\theta)|^2+|\nabla v(r_n,\theta)|^2\big)r_nd\theta\rightarrow 0,$$
which implies that
\be\lab{2013-10-26-e8}
\int_{\partial B_{R}(0)}|\nabla(u,v)|^2\sigma\cdot \nu d\sigma\rightarrow 0\;\hbox{as}\;R\rightarrow \infty.
\ee
Since $|G(x, u.,v)|\in L^1(\R^N)$, then by the the similar arguments, we obtain that
\be\lab{2013-10-26-e9}
\int_{\partial B_{R}(0)}G(\sigma, u,v)\sigma\cdot \nu d\sigma\rightarrow 0\;\hbox{as}\;R\rightarrow \infty.
\ee
When considering  $\Omega=B_R(0)$ in (\ref{2013-10-26-e4}), it follows that $\nu=\frac{x}{|x|}, \;\sigma\cdot \nu=|x|, \;0\leq(\nabla u \sigma\cdot \nabla u)\cdot \nu\leq |\nabla u|^2 \sigma \cdot \nu$. Hence by (\ref{2013-10-26-e8}), we have
\be\lab{2014-1-13-e1}
\int_{\partial B_R(0)}(\nabla u\sigma\cdot\nabla u- \sigma\frac{|\nabla u|^2}{2})\cdot \nu d\sigma\rightarrow 0\;\hbox{as}\;R\rightarrow +\infty.
\ee
Similarly, we also have
\be\lab{2014-1-13-e2}
\int_{\partial B_R(0)}(\nabla v\sigma\cdot\nabla v- \sigma\frac{|\nabla v|^2}{2})\cdot \nu d\sigma\rightarrow 0\;\hbox{as}\;R\rightarrow +\infty.
\ee
Finally, by the Lebesgue dominated convergence theorem, when $R\rightarrow \infty$, we have that
$$\int_{B_R(0)}|\nabla (u,v)|^2 dx\rightarrow \int_{\R^N}|\nabla (u,v)|^2dx, \int_{B_R(0)}G(x, u, v)dx\rightarrow \int_{\R^N}G(x, u, v)dx,$$ and
$$\sum_{i=1}^{N}\int_{B_R(0)}x_iG_{x_i}(x, u, v)dx\rightarrow \sum_{i=1}^{N}\int_{\R^N}x_iG_{x_i}(x, u, v)dx.$$
Combining with these results and (\ref{2013-10-26-e4}), (\ref{2013-10-26-e9}), (\ref{2014-1-13-e1}), (\ref{2014-1-13-e2})
we obtain (\ref{2013-10-26-e2}). The case of that $\Omega$ is a cone, we have $x\cdot \nu\equiv 0$ for $x\in\partial\Omega$, then \eqref{2013-10-26-e2} follows by  \eqref{2013-10-26-e1} easily.
\ep

\bc\lab{2014-11-27-cro1}
Let $0<\varepsilon<s<2$ and  $a_\varepsilon(x)$ be defined by \eqref{2014-10-8-e1}. Suppose that $\alpha>1,\beta>1,\alpha+\beta=2^*(s)$ and $\Omega$ is a cone. Then  any solution $(u,v)$ of
\be\lab{2014-11-27-e1}
\begin{cases}
-\Delta u=S_{\alpha,\beta,\lambda,\mu}^{\varepsilon}(\Omega)\Big(\lambda a_\varepsilon(x)|u|^{2^*(s)-2}u+\kappa\alpha a_\varepsilon(x)|u|^{\alpha-2}u|v|^\beta\Big)\quad &\hbox{in}\;\Omega,\\
-\Delta v=S_{\alpha,\beta,\lambda,\mu}^{\varepsilon}(\Omega)\Big(\mu a_\varepsilon(x)|v|^{2^*(s)-2}v+\kappa\beta a_\varepsilon(x)|u|^{\alpha}|v|^{\beta-2}v\Big)\quad &\hbox{in}\;\Omega,\\
u\geq 0, v\geq 0, (u,v)\in \mathscr{D}:=D_{0}^{1,2}(\Omega)\times D_{0}^{1,2}(\Omega),
\end{cases}
\ee
satisfies
\begin{align}\lab{2014-11-27-e6}
&\int_{\Omega\cap B_1}
\big[\frac{\lambda}{2^*(s)}a_\varepsilon(x)|u|^{2^*(s)}
+\frac{\mu}{2^*(s)}a_\varepsilon(x)|v|^{2^*(s)}
+2^*(s)\kappa a_\varepsilon(x)|u|^\alpha|v|^\beta\big]dx\\
=&\int_{\Omega\cap B_1^c}
\big[\frac{\lambda}{2^*(s)}a_\varepsilon(x)|u|^{2^*(s)}
+\frac{\mu}{2^*(s)}a_\varepsilon(x)|v|^{2^*(s)}
+2^*(s)\kappa a_\varepsilon(x)|u|^\alpha|v|^\beta\big]dx.\nonumber
\end{align}
\ec
\bp
Let
\be\lab{2014-11-27-e2}
G(x,u,v):=S_{\alpha,\beta,\lambda,\mu}^{\varepsilon}(\Omega)
\big[\frac{\lambda}{2^*(s)}a_\varepsilon(x)|u|^{2^*(s)}
+\frac{\mu}{2^*(s)}a_\varepsilon(x)|v|^{2^*(s)}
+\kappa a_\varepsilon(x)|u|^\alpha|v|^\beta\big].
\ee
Noting that
\be\lab{2014-11-27-e2}
\frac{\partial}{\partial x_i}a_\varepsilon(x)=\begin{cases}
-(s-\varepsilon)\frac{1}{|x|^{s+2-\varepsilon}}x_i\quad &\hbox{for}\;|x|<1,\\
-(s+\varepsilon)\frac{1}{|x|^{s+2+\varepsilon}}x_i\quad&\hbox{for}\;|x|>1,
\end{cases}
\ee
we have
\begin{align}\lab{2014-11-27-e3}
&x_i\cdot G_{x_i}(x, u,v)\\
=&S_{\alpha,\beta,\lambda,\mu}^{\varepsilon}(\Omega)
\big[\frac{\lambda}{2^*(s)}|u|^{2^*(s)}
+\frac{\mu}{2^*(s)}|v|^{2^*(s)}
+\kappa |u|^\alpha|v|^\beta\big]x_i\frac{\partial}{\partial x_i}a_\varepsilon(x)\nonumber\\
=&\begin{cases}
-(s-\varepsilon)S_{\alpha,\beta,\lambda,\mu}^{\varepsilon}(\Omega)
\big[\frac{\lambda}{2^*(s)}|u|^{2^*(s)}
+\frac{\mu}{2^*(s)}|v|^{2^*(s)}
+\kappa |u|^\alpha|v|^\beta\big]\frac{1}{|x|^{s+2-\varepsilon}}x_i^2\;\\\hbox{if}\;|x|<1,\\
-(s+\varepsilon)S_{\alpha,\beta,\lambda,\mu}^{\varepsilon}(\Omega)
\big[\frac{\lambda}{2^*(s)}|u|^{2^*(s)}
+\frac{\mu}{2^*(s)}|v|^{2^*(s)}
+\kappa |u|^\alpha|v|^\beta\big]\frac{1}{|x|^{s+2+\varepsilon}}x_i^2\;\\\hbox{if}\;|x|>1.
\end{cases}\nonumber
\end{align}
Hence, by Proposition \ref{2013-10-26-prop1}, we have
\begin{align}\lab{2014-11-27-e4}
&-2(N-s)\int_\Omega S_{\alpha,\beta,\lambda,\mu}^{\varepsilon}(\Omega)
\big[\frac{\lambda}{2^*(s)}a_\varepsilon(x)|u|^{2^*(s)}
+\frac{\mu}{2^*(s)}a_\varepsilon(x)|v|^{2^*(s)}
+\kappa a_\varepsilon(x)|u|^\alpha|v|^\beta\big]dx\\
&+(N-2)\int_\Omega \big(|\nabla u|^2+|\nabla v|^2\big)dx\nonumber\\
=&2\varepsilon \int_{\Omega\cap B_1}S_{\alpha,\beta,\lambda,\mu}^{\varepsilon}(\Omega)
\big[\frac{\lambda}{2^*(s)}a_\varepsilon(x)|u|^{2^*(s)}
+\frac{\mu}{2^*(s)}a_\varepsilon(x)|v|^{2^*(s)}
+\kappa a_\varepsilon(x)|u|^\alpha|v|^\beta\big]dx\nonumber\\
&-2\varepsilon\int_{\Omega\cap B_1^c}S_{\alpha,\beta,\lambda,\mu}^{\varepsilon}(\Omega)
\big[\frac{\lambda}{2^*(s)}a_\varepsilon(x)|u|^{2^*(s)}
 +\frac{\mu}{2^*(s)}a_\varepsilon(x)|v|^{2^*(s)}
+\kappa a_\varepsilon(x)|u|^\alpha|v|^\beta\big]dx.\nonumber
\end{align}
On the other hand, since $(u,v)$ is a solution, we have
\begin{align}\lab{2014-11-27-e5}
&\int_\Omega \big(|\nabla u|^2+|\nabla v|^2\big)dx\\
&=S_{\alpha,\beta,\lambda,\mu}^{\varepsilon}(\Omega)\int_\Omega \big[\lambda a_\varepsilon(x)|u|^{2^*(s)}
+\mu a_\varepsilon(x)|v|^{2^*(s)}
+2^*(s)\kappa a_\varepsilon(x)|u|^\alpha|v|^\beta\big]dx.\nonumber
\end{align}
Recalling that $S_{\alpha,\beta,\lambda,\mu}^{\varepsilon}(\Omega)>0, \varepsilon>0$, by \eqref{2014-11-27-e4} and \eqref{2014-11-27-e5}, we obtain \eqref{2014-11-27-e6}.
\ep

\bc\lab{2014-11-27-cro2}
Let $0<\varepsilon<s<2$ and  $a_\varepsilon(x)$ be defined by \eqref{2014-10-8-e1}. Suppose that $\alpha>1,\beta>1,\alpha+\beta=2^*(s)$ and $\Omega$ is a cone.
Let $(u_\varepsilon,v_\varepsilon)$ be a solution to \eqref{2014-11-26-xbue8} given by Lemma \ref{2014-11-26-xl4}. Then up to a subsequence, there exists some $(u,v)\in \mathscr{D}$ such that $u_\varepsilon\rightarrow u, v_\varepsilon\rightarrow v$ strongly in $D_{0}^{1,2}(\Omega)$ as $\varepsilon\rightarrow 0$.
\ec
\bp
By Lemma \ref{2014-11-26-xl5}, we only need to prove that $(u,v)\neq (0,0)$. Now,we proceed by contradiction. We assume that $u=v=0$.
Let $\chi(x)\in C_c^\infty(\Omega)$ be a cut-off function such that $\chi(x)\equiv 1$ in $B_{\frac{r}{2}}\cap \Omega$, $\chi(x)\equiv 0$ in $\Omega\backslash B_{r}$, recalling the Rellich-Kondrachov compact theorem and $2<2^*(s)<2^*:=\frac{2N}{N-2}$, we have $u_\varepsilon\rightarrow 0$ in $L^t(\Omega_1)$ for all $1<t<2^*$ if $0\not\in\bar{\Omega}_1$. Hence, it is easy to see that
\begin{align}\lab{2014-11-27-e7}
&\int_\Omega \big[\lambda a_\varepsilon(x)|\chi u_\varepsilon|^{2^*(s)}
+\mu a_\varepsilon(x)|\chi v_\varepsilon|^{2^*(s)}
+2^*(s)\kappa a_\varepsilon(x)|\chi u_\varepsilon|^\alpha|\chi v_\varepsilon|^\beta\big]dx\\
&=\int_{\Omega\cap B_r} \big[\lambda a_\varepsilon(x)|u_\varepsilon|^{2^*(s)}
+\mu a_\varepsilon(x)|v_\varepsilon|^{2^*(s)}
+2^*(s)\kappa a_\varepsilon(x)|u_\varepsilon|^\alpha|v_\varepsilon|^\beta\big]dx+o(1)\nonumber\\
&=:\eta_r+o(1).\nonumber
\end{align}
On the other hand, by the triangle inequality, we have
\begin{align}\lab{2014-11-27-e8}
&\int_\Omega \big(|\nabla (\chi u_\varepsilon)|^2+|\nabla (\chi v_\varepsilon)|^2\big)dx
=\int_\Omega \big(|(\nabla \chi)u_\varepsilon+\chi \nabla u_\varepsilon |^2+|(\nabla \chi)u_\varepsilon+\chi\nabla v_\varepsilon |^2\big)dx\\
\leq&\Bigg(\big(\int_\Omega |(\nabla \chi)|^2u_\varepsilon^2dx\big)^{\frac{1}{2}}+\big(\int_\Omega |\chi|^2|\nabla u_\varepsilon|^2dx\big)^{\frac{1}{2}}
+\big(\int_\Omega |(\nabla \chi)|^2v_\varepsilon^2dx\big)^{\frac{1}{2}}+\big(\int_\Omega |\chi|^2|\nabla v_\varepsilon|^2dx\big)^{\frac{1}{2}}\Bigg)^2\nonumber\\
=&\int_{\Omega\cap B_r}\big(|\nabla u_\varepsilon|^2+|\nabla v_\varepsilon|^2\big)dx+o(1)
:=\sigma_r+o(1).\nonumber
\end{align}
By \eqref{2014-11-27-e7} and \eqref{2014-11-27-e8}, we obtain that
\be\lab{2014-11-27-e9}
S_{\alpha,\beta,\lambda,\mu}^{\varepsilon}(\Omega) \big(\eta_r+o(1)\big)^{\frac{2}{2^*(s)}}\leq \sigma_r+o(1).
\ee
Similarly, we take $\tilde{\chi}(x)\in C^\infty(\Omega)$ such that $\tilde{\chi}(x)\equiv 0$ in $B_r\cap \Omega$ and $\tilde{\chi}\equiv 1$ in $\Omega\backslash B_{2r}$. Then by repeating the above steps, we obtain that
\be\lab{2014-11-27-e10}
S_{\alpha,\beta,\lambda,\mu}^{\varepsilon}(\Omega) \big(1-\eta_r+o(1)\big)^{\frac{2}{2^*(s)}}\leq S_{\alpha,\beta,\lambda,\mu}^{\varepsilon}-\sigma_r+o(1).
\ee
By \eqref{2014-11-27-e9} and \eqref{2014-11-27-e10}, we deduce that
\be\lab{2014-11-27-e11}
\big(\eta_r+o(1)\big)^{\frac{2}{2^*(s)}}+\big(1-\eta_r+o(1)\big)^{\frac{2}{2^*(s)}}\leq 1.
\ee
Notice that $h(t):=t^{\frac{2}{2^*(s)}}+(1-t)^{\frac{2}{2^*(s)}}$ satisfying that
$\min_{t\in [0,1]}h(t)=1$ and only achieved by $t=0$ or $t=1$.
Hence, we obtain that $\eta_r\equiv0$ or $\eta_r\equiv1$ for any $r>0$.

But by Corollary \ref{2014-11-27-cro1},  for any $\varepsilon\in (0,s)$, we have
\begin{align*}
&\int_{\Omega\cap B_1}
\big[\lambda a_\varepsilon(x)|u_\varepsilon|^{2^*(s)}
+\mu a_\varepsilon(x)|v_\varepsilon|^{2^*(s)}
+2^*(s)\kappa a_\varepsilon(x)|u_\varepsilon|^\alpha|v_\varepsilon|^\beta\big]dx\nonumber\\
&=\int_{\Omega\cap B_1^c}
\big[\lambda a_\varepsilon(x)|u_\varepsilon|^{2^*(s)}
+\mu a_\varepsilon(x)|v_\varepsilon|^{2^*(s)}
+2^*(s)\kappa a_\varepsilon(x)|u_\varepsilon|^\alpha|v_\varepsilon|^\beta\big]dx.
\end{align*}
Combined with the fact of that
\begin{align*}
&\int_{\Omega\cap B_1}
\big[\lambda a_\varepsilon(x)|u_\varepsilon|^{2^*(s)}
+\mu a_\varepsilon(x)|v_\varepsilon|^{2^*(s)}
+2^*(s)\kappa a_\varepsilon(x)|u_\varepsilon|^\alpha|v_\varepsilon|^\beta\big]dx\nonumber\\
&+\int_{\Omega\cap B_1^c}
\big[\lambda a_\varepsilon(x)|u_\varepsilon|^{2^*(s)}
+\mu a_\varepsilon(x)|v_\varepsilon|^{2^*(s)}
+2^*(s)\kappa a_\varepsilon(x)|u_\varepsilon|^\alpha|v_\varepsilon|^\beta\big]dx=1,
\end{align*}
we obtain that
\begin{align}\lab{2014-11-27-e12}
&\frac{1}{2}=\int_{\Omega\cap B_1}
\big[\lambda a_\varepsilon(x)|u_\varepsilon|^{2^*(s)}
+\mu a_\varepsilon(x)|v_\varepsilon|^{2^*(s)}
+2^*(s)\kappa a_\varepsilon(x)|u_\varepsilon|^\alpha|v_\varepsilon|^\beta\big]dx\\
&= \int_{\Omega\cap B_1^c}
\big[\lambda a_\varepsilon(x)|u_\varepsilon|^{2^*(s)}
+\mu a_\varepsilon(x)|v_\varepsilon|^{2^*(s)}
+2^*(s)\kappa a_\varepsilon(x)|u_\varepsilon|^\alpha|v_\varepsilon|^\beta\big]dx.\nonumber
\end{align}
Hence, we have $\eta_r\equiv\frac{1}{2}$ for any $r>0$, a contradiction.
\ep

\noindent{\bf Proof of Theorem \ref{2014-11-26-th1}:}
Let $(u_{\varepsilon_k}, v_{\varepsilon_k})$ be a solution to \eqref{2014-11-26-xbue8} with $\varepsilon=\varepsilon_k$ given by Lemma \ref{2014-11-26-xl4} and $\varepsilon_k\rightarrow 0$ as $k\rightarrow +\infty$. Up to a subsequence, we may assume that $u_{\varepsilon_k}\rightharpoonup u, v_{\varepsilon_k}\rightharpoonup v$ in $D_{0}^{1,2}(\Omega)$ and $u_{\varepsilon}\rightarrow u, v_{\varepsilon}\rightarrow v$ a.e. in $\Omega$ (see Lemma \ref{2014-11-26-xl5}). Then if $0\not\in \bar{\Omega}$, by Rellich-Kondrachov compact theorem, it is easy to see that $u_{\varepsilon}\rightarrow u, v_{\varepsilon_k}\rightarrow v$ strongly in $D_{0}^{1,2}(\Omega)$.
When $\Omega$ is a cone, by Corollary \ref{2014-11-27-cro2}, we also obtain that $u_{\varepsilon}\rightarrow u, v_{\varepsilon_k}\rightarrow v$ strongly in $D_{0}^{1,2}(\Omega)$. By $(iv)$ of Lemma \ref{2014-11-26-xl5}, we obtain that $(u,v)$ is an  extremal of $S_{\alpha,\beta,\lambda,\mu}(\Omega)$, the proof is completed.\hfill$\Box$

\br\lab{2014-11-27-r1}
When $\Omega$ is a cone, let $(u,v)$ be the extremal obtained as limit of $(u_{\varepsilon_k}, v_{\varepsilon_k})$, the solution to \eqref{2014-11-26-xbue8} with $\varepsilon=\varepsilon_k$ given by Lemma \ref{2014-11-26-xl4}. Then by Lemma \ref{2014-11-26-xl5}, Corollary \ref{2014-11-27-cro2} and the formula \eqref{2014-11-27-e12}, we see that $(u,v)$ satisfies
\begin{align}\lab{2014-11-27-e13}
&\frac{1}{2}=\int_{\Omega\cap B_1}
\big[\lambda a_0(x)|u|^{2^*(s)}
+\mu a_0(x)|v|^{2^*(s)}
+2^*(s)\kappa a_0(x)|u|^\alpha|v|^\beta\big]dx\nonumber\\
&=\int_{\Omega\cap B_1^c}
\big[\lambda a_0(x)|u|^{2^*(s)}
+\mu a_0(x)|v|^{2^*(s)}
+2^*(s)\kappa a_0(x)|u|^\alpha|v|^\beta\big]dx.\nonumber
\end{align}
Such a property has been observed for the scalar equation. \er

\subsection{Existence of positive ground state solutions}
By \eqref{2014-11-26-xe3}, we always have $S_{\alpha,\beta,\lambda,\mu}(\Omega)\leq \big(\max\{\lambda, \mu\}\big)^{-\frac{2}{2^*(s)}}\mu_s(\Omega)$. When $\Omega$ is a cone and  $s\in (0,2), \alpha>1,\beta>1,\alpha+\beta=2^*(s)$, by Theorem \ref{2014-11-26-th1}, $S_{\alpha,\beta,\lambda,\mu}(\Omega)$ is always attained (although the extremals may be semi-trivial). Indeed, if $S_{\alpha,\beta,\lambda,\mu}(\Omega)= \big(\max\{\lambda, \mu\}\big)^{-\frac{2}{2^*(s)}}\mu_s(\Omega)$,  then $S_{\alpha,\beta,\lambda,\mu}(\Omega)$ can be achieved by semi-trivial function. To see this, we just need to plug in the pairs $(U, 0)$  or $(0, U)$, where
$U$ is an extremal function of $\mu_s(\Omega)$.  But, under some special conditions, $S_{\alpha,\beta,\lambda,\mu}(\Omega)$ can also be achieved by nontrivial function even $S_{\alpha,\beta,\lambda,\mu}(\Omega)= \big(\max\{\lambda, \mu\}\big)^{-\frac{2}{2^*(s)}}\mu_s(\Omega)$, see
Theorem \ref{2014-11-30-th2} below.  However, if
\be\lab{2014-11-27-lae1}
S_{\alpha,\beta,\lambda,\mu}(\Omega)< \big(\max\{\lambda, \mu\}\big)^{-\frac{2}{2^*(s)}}\mu_s(\Omega),
\ee
then the  extremal functions (hence the positive ground state of the system \eqref{2015-1-5-e1})  of  $S_{\alpha,\beta,\lambda,\mu}(\Omega)$ must be nontrivial.  Therefore,
next we   need to search some sufficient conditions to ensure the  above  strict inequality \eqref{2014-11-27-lae1}. We obtain the following  theorem on the existence, regularity and decay estimate.


\bt\lab{2014-11-27-wth1}
Let $\Omega$ be a cone in $\R^N$(especially, $\Omega=\R^N$ or $\Omega=\R_+^N$) or $\Omega$ be an open domain but $0\not\in \bar{\Omega}$. Assume  $s\in (0,2),\kappa>0, \alpha>1,\beta>1,\alpha+\beta=2^*(s)$. Then system  \eqref{2015-1-5-e1} possesses a positive ground state solution $(\phi,\psi)$ (i.e., $\phi>0, \psi>0$) provided  that  one of the following conditions holds:
\begin{itemize}
\item[$(a_1)$]$\lambda>\mu$ and either $ 1<\beta<2$ or $\begin{cases}\beta=2\\ \kappa>\frac{\lambda}{2^*(s)} \end{cases}$;
\item[$(a_2)$] $\lambda=\mu$ and either $ \min\{\alpha,\beta\}<2$ or $\begin{cases} \min\{\alpha,\beta\}=2,\\ \kappa>\frac{\lambda}{2^*(s)} \end{cases}$;
\item[$(a_3)$] $\lambda<\mu$ and either $ 1<\alpha<2$ or $\begin{cases}\alpha=2\\ \kappa>\frac{\mu}{2^*(s)} \end{cases}$.
\end{itemize}
Moreover, when $\Omega$ is a cone, we have the following regularity and decay properties:
\begin{itemize}
\item[$(b_1)$] if $0<s<\frac{N+2}{N}$, $\phi,\psi\in C^2(\overline{\Omega})$;
\item[$(b_2)$] if $s=\frac{N+2}{N}$, $\phi, \psi\in C^{1,\gamma}({\Omega})$ for all $0<\gamma<1$;
\item[$(b_3)$]if $s>\frac{N+2}{N}$, $\phi,\psi\in C^{1,\gamma}({\Omega})$ for all $0<\gamma<\frac{N(2-s)}{N-2}$.
\end{itemize}
When   $\Omega$ is a cone with $0\in \partial\Omega$ (e.g., $\Omega=\R_+^N$), then there exists a constant $C$ such that $$|\phi(x)|, |\psi(x)|\leq C (1+|x|^{-(N-1)}), \quad \quad |\nabla \phi(x)|, |\nabla \psi(x)|\leq C|x|^{-N}.$$
When $\Omega=\R^N$,
$$|\phi(x)|, |\psi(x)|\leq C (1+|x|^{-N}), \quad \quad |\nabla \phi(x)|, |\nabla \psi(x)|\leq C|x|^{-N-1}$$
In particular, if $\Omega=\R_+^N$, then  $(\phi(x),\psi(x))$ is axially symmetric with respect to the $x_N$-axis, i.e., $$\big(\phi(x),\psi(x)\big)=\big(\phi(x',x_N),
\psi(x',x_N)\big)=\big(\phi(|x'|,x_N),\psi(|x'|,x_N)\big).$$
\et

\br\lab{2014-11-27-wbur1}
The conditions $(a_1)-(a_3)$ imposed  in Theorem \ref{2014-11-27-wth1} are some sufficient conditions to ensure the inequality \eqref{2014-11-27-lae1} (see Lemma \ref{2014-11-27-xbul1} below).
But they are not necessary conditions. For example,  when  $\lambda>\mu$ and $1<\alpha<2$, we  can  not   exclude that
$\displaystyle S_{\alpha,\beta,\lambda,\mu}(\Omega)<\lambda ^{-\frac{2}{2^*(s)}}\mu_s(\Omega)$, i.e., \eqref{2014-11-27-lae1} may  be true.
\er

Define the functional
\begin{align}\lab{2014-11-27-we3}
\Phi(u,v)=&\frac{1}{2}\int_\Omega \big(|\nabla u|^2+|\nabla v|^2\big)dx\\
&-\frac{1}{2^*(s)}\int_\Omega \frac{1}{|x|^s}\big[\lambda |u|^{2^*(s)}+\mu|v|^{2^*(s)}+2^*(s)\kappa |u|^\alpha|v|^\beta\big]dx\nonumber
\end{align}
and the corresponding Nehari manifold
\be\lab{2014-11-27-we5}
\mathcal{N}:=\{(u,v)\in \mathscr{D}\backslash \{0,0\}:\;J(u,v)=0\}
\ee
where
\begin{align}\lab{2014-11-27-we4}
&J(u,v):=\langle \Phi'(u,v), (u,v)\rangle\\
&=\int_\Omega \big(|\nabla u|^2+|\nabla v|^2\big)dx-\int_\Omega \frac{1}{|x|^s}\big[\lambda |u|^{2^*(s)}+\mu|v|^{2^*(s)}+2^*(s)\kappa |u|^\alpha|v|^\beta\big]dx.\nonumber
\end{align}
By Lemma \ref{2014-11-23-xl1}, $\mathcal{N}$ is well defined. Define
\be\lab{2014-11-27-we6}
c_0:=\inf_{(u,v)\in \mathcal{N}}\Phi(u,v),
\ee
then basing on the results of Section 4 and Section 6, we see that
\be\lab{2014-11-27-we7}
0<c_0\leq [\frac{1}{2}-\frac{1}{2^*(s)}]\big[\mu_{s}(\Omega)\big]^{\frac{2^*(s)}{2^*(s)-2}} \big(\max\{\lambda,\mu\}\big)^{-\frac{2}{2^*(s)-2}}.
\ee
Moreover, we have the following result.
\bl\lab{2014-11-27-zl1}
Let $\Omega$ be a cone of $\R^N$ or  $\Omega$ be an open domain but $0\not\in \bar{\Omega}$. Assume that $\kappa>0, s\in (0,2), \alpha>1,\beta>1,\alpha+\beta=2^*(s)$ and let
$c_0$ be defined by \eqref{2014-11-27-we6}, then we have
\be\lab{2014-11-27-we8}
c_0< \left[\frac{1}{2}-\frac{1}{2^*(s)}\right]\big(\mu_{s}(\Omega)\big)^{\frac{2^*(s)}{2^*(s)-2}} \big(\max\{\lambda,\mu\}\big)^{-\frac{2}{2^*(s)-2}}
\ee
if and only if
\be\lab{2014-11-27-we9}
S_{\alpha,\beta,\lambda,\mu}(\Omega)< \big(\max\{\lambda, \mu\}\big)^{-\frac{2}{2^*(s)}}\mu_s(\Omega).
\ee
\el
\bp
A direct computation shows that
\be\lab{2014-11-27-we12}
c_0=\left[\frac{1}{2}-\frac{1}{2^*(s)}\right]\big(S_{\alpha,\beta,\lambda,\mu}(\Omega)\big)^{\frac{2^*(s)}{2^*(s)-2}}.
\ee
\ep

Then combining with the conclusions of Section 6, we have the following result:
\bl\lab{2014-11-27-xbul1}
Let $\Omega$ be a cone in $\R^N$(especially, $\Omega=\R^N$ and $\Omega=\R_+^N$ ) or  $\Omega$ be an open domain but $0\not\in \bar{\Omega}$. Suppose  $s\in (0,2),\kappa>0, \alpha>1,\beta>1,\alpha+\beta=2^*(s)$. Then $\displaystyle S_{\alpha,\beta,\lambda,\mu}(\Omega)< \big(\max\{\lambda, \mu\}\big)^{-\frac{2}{2^*(s)}}\mu_s(\Omega)$ if one of the following holds:
\begin{itemize}
\item[$(i)$]$\lambda>\mu$ and either $ 1<\beta<2$ or $\begin{cases}\beta=2\\ \kappa>\frac{\lambda}{2^*(s)} \end{cases}$;
\item[$(ii)$] $\lambda=\mu$ and either $ \min\{\alpha,\beta\}<2$ or $\begin{cases} \min\{\alpha,\beta\}=2,\\ \kappa>\frac{\lambda}{2^*(s)} \end{cases}$;
\item[$(iii)$] $\lambda<\mu$ and either $ 1<\alpha<2$ or $\begin{cases}\alpha=2\\ \kappa>\frac{\mu}{2^*(s)} \end{cases}$.
\end{itemize}
\el
\bp
It follows by Corollary \ref{2014-11-24-xcro1}, Corollary \ref{2014-4-12-Cro4}, Lemma \ref{2014-11-24-l3} and Lemma \ref{2014-11-27-zl1}.
\ep

\noindent{\bf Proof of Theorem \ref{2014-11-27-wth1}:} Under the assumptions of Theorem \ref{2014-11-27-wth1}, firstly by  Theorem \ref{2014-11-26-th1}, $S_{\alpha,\beta,\lambda,\mu}(\Omega)$ is attained by some   nonnegative pair $(u,v)$ such that  $(u,v)\neq (0,0)$. On the other hand, by Lemma \ref{2014-11-27-xbul1}, we have
$\displaystyle S_{\alpha,\beta,\lambda,\mu}(\Omega)< \big(\max\{\lambda, \mu\}\big)^{-\frac{2}{2^*(s)}}\mu_s(\Omega).$
Hence, we see that $u\neq 0, v\neq 0$. Hence, by Proposition \ref{zzz=111},
$$ (\phi,\psi):=\big((S_{\alpha,\beta,\lambda,\mu}(\Omega))^{\frac{1}{2^*(s)-2}}u, (S_{\alpha,\beta,\lambda,\mu}(\Omega))^{\frac{1}{2^*(s)-2}}v\big)$$ is a ground state solution of system  \eqref{2015-1-5-e1}. Then by the strong maximum principle, it is easy to see that $\phi>0,, \psi>0$ in $\Omega$. We note that the arguments in Proposition \ref{2014-11-21-prop1} and Proposition \ref{2014-11-21-Prop2} are valid for general cone. Combining with Proposition \ref{2014-11-22-Prop1}, we complete the proof. \hfill$\Box$

\subsection{Uniqueness and Nonexistence of positive ground state solutions}

In the previous subsection, in Theorem \ref{2014-11-27-wth1}, we have established  the existence
of the positive ground state solution to the  system  \eqref{2015-1-5-e1}.  Now, in the current subsection, we
obtain the uniqueness  of  the positive ground state solution  to the  system  \eqref{2015-1-5-e1}.
Define
\be\lab{2014-11-28-xbue1}
G(u,v):=\frac{\int_\Omega \big(|\nabla u|^2+|\nabla v|^2\big)dx}{\Big(\int_\Omega \big(\lambda \frac{|u|^{2^*(s)}}{|x|^s}+\mu \frac{|v|^{2^*(s)}}{|x|^s}+2^*(s)\kappa \frac{|u|^\alpha |v|^\beta}{|x|^s}\big)dx\Big)^{\frac{2}{2^*(s)}}},\;\;(u,v)\neq (0,0)
\ee
then we have
\be\lab{2014-11-28-xbue2}
S_{\alpha,\beta,\lambda,\mu}(\Omega)=\inf_{(u,v)\in \mathscr{D}\backslash\{(0,0)\}} G(u,v).
\ee
For any $u\neq 0, v\neq 0$ and $t\geq 0$,  we have
\be\lab{2014-11-28-xbue3}
G(u,tv)=\frac{\int_\Omega \big(|\nabla u|^2+|\nabla v|^2t^2\big)dx}{\Big(\int_\Omega \big(\lambda \frac{|u|^{2^*(s)}}{|x|^s}+\mu \frac{|v|^{2^*(s)}}{|x|^s}t^{2^*(s)}+2^*(s)\kappa \frac{|u|^\alpha |v|^\beta}{|x|^s}t^\beta\big) dx\Big)^{\frac{2}{2^*(s)}}}.
\ee
In particular,
\begin{align}\lab{2014-11-28-xbue7}
G(u,tu)
=&\frac{1+t^2}{\left[\lambda+\mu t^{2^*(s)}+2^*(s)\kappa t^\beta\right]^{\frac{2}{2^*(s)}}} \frac{\int_\Omega |\nabla u|^2 dx}{\left(\int_\Omega \frac{|u|^{2^*(s)}}{|x|^s}dx\right)^{\frac{2}{2^*(s)}}}\nonumber\\
:=&g(t)\frac{\int_\Omega |\nabla u|^2 dx}{\left(\int_\Omega \frac{|u|^{2^*(s)}}{|x|^s}dx\right)^{\frac{2}{2^*(s)}}}.
\end{align}
We define $\displaystyle g(+\infty)=\lim_{t\rightarrow +\infty}g(t)=\mu^{-\frac{2}{2^*(s)}}$, then we  see that
$$G(0,v)=\lim_{t\rightarrow +\infty}G(v,tv)=g(+\infty)\frac{\int_\Omega |\nabla v|^2 dx}{\left(\int_\Omega \frac{|v|^{2^*(s)}}{|x|^s}dx\right)^{\frac{2}{2^*(s)}}}.$$
Hence, we have
\begin{align}\lab{2014-11-28-xbue8}
S_{\alpha,\beta,\lambda,\mu}(\Omega)=&\inf_{(u,v)\in \mathscr{D}\backslash\{(0,0)\}}G(u,v)\\
\leq&\inf_{u\in D_{0}^{1,2}(\Omega)}\inf_{t\in [0,+\infty)}G(u,tu)\nonumber\\
=&\inf_{t\in [0,+\infty)}g(t) \inf_{u\in D_{0}^{1,2}(\Omega)\backslash \{0\}}\frac{\int_\Omega |\nabla u|^2 dx}{\left(\int_\Omega \frac{|u|^{2^*(s)}}{|x|^s}dx\right)^{\frac{2}{2^*(s)}}}\nonumber\\
=&\inf_{t\in [0,+\infty)}g(t) \mu_s(\Omega).\nonumber
\end{align}
Moreover, we can obtain the follow precise result:
\bl\lab{2014-11-28-xl1}
$\displaystyle S_{\alpha,\beta,\lambda,\mu}(\Omega)=\inf_{t\in [0,+\infty)}g(t) \mu_s(\Omega)$, where
\be\lab{2014-11-28-xbue9}
g(t):=\frac{1+t^2}{\left[\lambda+\mu t^{2^*(s)}+2^*(s)\kappa t^\beta\right]^{\frac{2}{2^*(s)}}}.
\ee
\el
\bp
By \eqref{2014-11-28-xbue8}, we only need to prove the reverse inequality. Now, let $\{(u_n,v_n)\}$ be a minimizing sequence of $S_{\alpha,\beta,\lambda,\mu}(\Omega)$. Since $G(u,v)=G(tu,tv)$ for all $t>0$, without loss of generality, we may assume that
$\displaystyle\int_\Omega \big(\frac{|u_n|^{2^*(s)}}{|x|^s}+\frac{|v_n|^{2^*(s)}}{|x|^s}\big)\equiv 1$
and $\displaystyle G(u_n,v_n)=S_{\alpha,\beta,\lambda,\mu}(\Omega)+o(1).$

\noindent{\bf Case 1:} $\displaystyle\liminf_{n\rightarrow +\infty}\int_\Omega \frac{|u_n|^{2^*(s)}}{|x|^s} dx=0$. Since $\{v_n\}$ is bounded in $L^{2^*(s)}(\Omega,\frac{dx}{|x|^s})$, by H\"older inequality, up to a subseqeunce, we see that
$$\int_\Omega \big(\lambda \frac{|u_n|^{2^*(s)}}{|x|^s}+\mu \frac{|v_n|^{2^*(s)}}{|x|^s}+2^*(s)\kappa \frac{|u_n|^\alpha |v_n|^\beta}{|x|^s}\big)dx=\mu\int_\Omega \frac{|v_n|^{2^*(s)}}{|x|^s}dx+o(1)=\mu+o(1).$$
 Hence
 $\displaystyle \lim_{n\rightarrow \infty}G(u,v_n)\geq \lim_{n\rightarrow \infty} G(0,v_n)$.
 We see that $(0,v_n)$ is also a minimizing sequence of $S_{\alpha,\beta,\lambda,\mu}(\Omega)$. Hence, it is easy to see that
\begin{align}\lab{2014-11-28-xbue10}
S_{\alpha,\beta,\lambda,\mu}(\Omega)=\mu^{-\frac{2}{2^*(s)}}\mu_s(\Omega)
=g(+\infty)\mu_s(\Omega)
\geq\inf_{t\in (0,+\infty)}g(t) \mu_s(\Omega).
\end{align}

\noindent{\bf Case 2:} $\displaystyle\liminf_{n\rightarrow +\infty}\int_\Omega \frac{|v_n|^{2^*(s)}}{|x|^s} dx=0$. Similarly to Case 1, we can obtain that
\begin{align}\lab{2014-11-28-xbue11}
S_{\alpha,\beta,\lambda,\mu}(\Omega)=\lambda^{-\frac{2}{2^*(s)}}\mu_s(\Omega)
=g(0)\mu_s(\Omega)
\geq\inf_{t\in [0,+\infty)}g(t) \mu_s(\Omega).
\end{align}

\noindent{\bf Case 3:} Up to a subseqeuce if necessary, we may assume that $\displaystyle\lim_{n\rightarrow +\infty}\int_\Omega \frac{|u_n|^{2^*(s)}}{|x|^s} dx=\delta>0$ and $\displaystyle\lim_{n\rightarrow +\infty}\int_\Omega \frac{|v_n|^{2^*(s)}}{|x|^s} dx=1-\delta>0$. Let $t_n>0$ such that $$ \int_\Omega \frac{|v_n|^{2^*(s)}}{|x|^s}dx=\int_\Omega \frac{|t_nu_n|^{2^*(s)}}{|x|^s}dx,$$ then we see that $\{t_n\}$ is bounded and away from $0$. Up to a subsequence, we may assume that $t_n\rightarrow t_0=\left(\frac{\delta}{1-\delta}\right)^{\frac{1}{2^*(s)}}$.
Now let $w_n=\frac{1}{t_n}v_n$, then we have
\be\lab{2014-11-28-xbue12}
\int_\Omega \frac{|u_n|^{2^*(s)}}{|x|^s}dx=\int_\Omega\frac{|w_n|^{2^*(s)}}{|x|^s}dx
\ee
and by Young's inequality, we have
\begin{align}\lab{2014-11-28-xbue13}
\int_\Omega \frac{|u_n|^\alpha|w_n|^\beta}{|x|^s}dx\leq& \frac{\alpha}{2^*(s)}\int_\Omega \frac{|u_n|^{2^*(s)}}{|x|^s}dx+\frac{\beta}{2^*(s)}\int_\Omega \frac{|w_n|^{2^*(s)}}{|x|^s}dx\\
=&\int_\Omega \frac{|u_n|^{2^*(s)}}{|x|^s}dx=\int_\Omega\frac{|w_n|^{2^*(s)}}{|x|^s}dx.\nonumber
\end{align}
Hence,
\begin{align}\lab{2014-11-28-xbue14}
&G(u_n,v_n)=G(u_n,t_nw_n)\\
=&\frac{\int_\Omega |\nabla u_n|^2dx}{\Big(\int_\Omega \big(\lambda \frac{|u_n|^{2^*(s)}}{|x|^s}+\mu t_{n}^{2^*(s)}\frac{|w_n|^{2^*(s)}}{|x|^s}+2^*(s)\kappa t_n^\beta\frac{|u_n|^\alpha |w_n|^\beta}{|x|^s}\big)dx\Big)^{\frac{2}{2^*(s)}}}\nonumber\\
&+\frac{\int_\Omega t_n^2|\nabla w_n|^2dx}{\Big(\int_\Omega \big(\lambda \frac{|u_n|^{2^*(s)}}{|x|^s}+\mu t_{n}^{2^*(s)}\frac{|w_n|^{2^*(s)}}{|x|^s}+2^*(s)\kappa t_n^\beta\frac{|u_n|^\alpha |w_n|^\beta}{|x|^s}\big)dx\Big)^{\frac{2}{2^*(s)}}}\nonumber\\
\geq&\frac{1}{\left[\lambda+\mu t_{n}^{2^*(s)}+2^*(s)\kappa t_n^\beta\right]^{\frac{2}{2^*(s)}}}\frac{\int_\Omega |\nabla u_n|^2 dx}{\left(\int_\Omega \frac{|u_n|^{2^*(s)}}{|x|^s}dx\right)^{\frac{2}{2^*(s)}}}\nonumber\\
&+\frac{t_n^2}{\left[\lambda+\mu t_{n}^{2^*(s)}+2^*(s)\kappa t_n^\beta\right]^{\frac{2}{2^*(s)}}}\frac{\int_\Omega |\nabla w_n|^2 dx}{\left(\int_\Omega \frac{|w_n|^{2^*(s)}}{|x|^s}dx\right)^{\frac{2}{2^*(s)}}}\nonumber\\
\geq&\frac{1}{\left[\lambda+\mu t_{n}^{2^*(s)}+2^*(s)\kappa t_n^\beta\right]^{\frac{2}{2^*(s)}}}\mu_{s}(\Omega)
 +\frac{t_n^2}{\left[\lambda+\mu t_{n}^{2^*(s)}+2^*(s)\kappa t_n^\beta\right]^{\frac{2}{2^*(s)}}}\mu_{s}(\Omega)\nonumber\\
=&g(t_n)\mu_{s}(\Omega).\nonumber
\end{align}
Let $n\rightarrow +\infty$, we obtain that
$$ S_{\alpha,\beta,\lambda,\mu}(\Omega)\geq g(t_0)\mu_s(\Omega)\geq \inf_{t\in (0,+\infty)}g(t) \mu_s(\Omega).$$
Thereby $\displaystyle S_{\alpha,\beta,\lambda,\mu}(\Omega)= \inf_{t\in (0,+\infty)}g(t) \mu_s(\Omega)$ is proved.
\ep

Basing on Lemma \ref{2014-11-28-xl1}, we can propose the ``uniqueness" type result as following:

\bt\lab{2014-11-30-th1}
Let $\Omega$ either be a cone in $\R^N$(in particular, $\Omega=\R^N$ and $\Omega=\R_+^N$) or  $\Omega$ be an open domain but $0\not\in \bar{\Omega}$. Assume  $s\in (0,2),\kappa>0, \alpha>1,\beta>1,\alpha+\beta=2^*(s)$. Let $(\phi,\psi)$ be a positive ground state solution to problem \eqref{2015-1-5-e1}, then $$ \phi=C(t_0)U,\quad  \psi=t_0C(t_0)U,$$
where $U$ is the ground state solution of
\be\lab{2014-11-30-e1}
\begin{cases}
-\Delta u=\mu_{s}(\Omega) \frac{u^{2^*(s)-1}}{|x|^{s}}\;\hbox{in}\;\Omega,\\
u=0\;\;\hbox{on}\;\partial \Omega,
\end{cases}
\ee
while $t_0>0$ satisfies that
\be\lab{2014-11-30-xe1}
g(t_0)=\inf_{t\in(0,+\infty)}g(t)
\ee
and
\be\lab{2014-11-30-xe2}
C(t_0):=\left[S_{\alpha,\beta,\lambda,\mu}(\Omega)\right]^{\frac{1}{2^*(s)-2}}
\left(\lambda+\mu t_{0}^{2^*(s)}+2^*(s)\kappa t_0^\beta\right)^{-\frac{1}{2^*(s)}},
\ee
where $g(t)$ is defined in \eqref{2014-11-28-xbue9}.
\et
\bp
By the processes of Case 3 in the proof of Lemma \ref{2014-11-28-xl1}, if $S_{\alpha,\beta,\lambda,\mu}(\Omega)$ is attained by some nontrivial function $(u,v)$, i.e., $u\neq 0,v\neq0$, then there exists some $t_0>0$ such that $v=t_0u$,  where $u$ is a minimizer of $\mu_s(\Omega)$ and $t_0$ satisfies $\displaystyle g(t_0)=\inf_{t\in (0,+\infty)}g(t)$.

Now assume that $u=CU, v=t_0CU$, then a direct computation shows that
$$ \int_\Omega \frac{1}{|x|^s}\big[\lambda |u|^{2^*(s)}+\mu|v|^{2^*(s)}+2^*(s)\kappa |u|^\alpha|v|^\beta\big]dx=1$$
if and only if
$$ C=\left(\lambda+\mu t_{0}^{2^*(s)}+2^*(s)\kappa t_0^\beta\right)^{-\frac{1}{2^*(s)}}.$$
Finally,  we see that $\phi=\left[S_{\alpha,\beta,\lambda,\mu}(\Omega)\right]^{\frac{1}{2^*(s)-2}}u, \psi=\left[S_{\alpha,\beta,\lambda,\mu}(\Omega)\right]^{\frac{1}{2^*(s)-2}} v$,  we complete the proof.
\ep

\br\lab{2014-11-30-xbur1} Under the assumption that
 $$ S_{\alpha,\beta,\lambda,\mu}(\Omega)< \big(\max\{\lambda, \mu\}\big)^{-\frac{2}{2^*(s)}}\mu_s(\Omega),$$
 we  have seen that the  problem \eqref{2015-1-5-e1} possesses a positive ground state solution.  But the converse usually is not true.
Next, we construct an example where  $S_{\alpha,\beta,\lambda,\mu}(\Omega)=\big(\max\{\lambda, \mu\}\big)^{-\frac{2}{2^*(s)}}\mu_s(\Omega)$
but  the  system  \eqref{2015-1-5-e1} still has  (multiple)  positive ground state solutions.
\er

\bt\lab{2014-11-30-th2}
Let $\Omega$ be a cone in $\R^3$ or  $\Omega$ be an open domain but $0\not\in \bar{\Omega}$.
Assume the following conditions hold:
\begin{itemize}
\item[(a)] $0<s<2$,
\item[(b)] either $\alpha>2$ or $\begin{cases} \alpha=2\\ \mu\geq 2\kappa \end{cases}$,
\item[(c)] either  $\beta>2$ or $\begin{cases} \beta=2\\ \lambda\geq 2\kappa \end{cases}$,
\item[(d)] $\alpha+\beta=2^*(s)$.
\end{itemize}
Then we have one of the following conclusion:\\
(1) If  $s=1,\alpha=\beta=2, \lambda=\mu=2\kappa>0$,  then the set of all extremal functions  of $S_{\alpha,\beta,\lambda,\mu}(\Omega)$ is  given by
\be\lab{2014-11-29-e3}
\mathcal{A}:=\left\{(t_1U,t_2U):t_1\geq 0,t_2\geq 0, (t_1,t_2)\neq (0,0)\;\hbox{and $U$ is an extremal of $\mu_s(\Omega)$}\;\right\}.
\ee
(2) Except for the item (1) above,
$S_{\alpha,\beta,\lambda,\mu}(\Omega)$ has no  nontrivial extremal  function.
\et

\br Under the hypotheses $(a)$-$(d)$, the dimension of the space $\R^N$ has to be three. Therefore, we can only
establish the above theorem in $\R^3$.\er

\bp
We proceed by contradiction. Assume that $S_{\alpha,\beta,\lambda,\mu}(\Omega)$ has a  nontrivial extremal $(u,v)$, then by Lemma \ref{2014-11-28-xl1} (see case 3 of the proof), we see that there exists some $t_0>0$ such that $v=t_0 u$ and $u$ is an extremal of $\mu_s(\Omega)$. Moreover, $t_0$ attains  the minimum of $g(t)$, where $g(t)$ is introduced in \eqref{2014-11-28-xbue9}.  By conditions $(b)$ and $(c)$, we see that $g''(0)\geq 0$ and $g'(t)<0$ for $t$ large enough. Hence, $\{t>0:\; g'(t)=0\}$ has at least $3$ solutions $\{t_1,t_2,t_3\}$ such that $0<t_1<t_2=t_0<t_3<\infty$. A direct computation shows that
\begin{align}\lab{2014-11-28-we1}
g'(t)
=\frac{-2t}{\left[\lambda+\mu t^{2^*(s)}+2^*(s)\kappa t^\beta\right]^{\frac{2}{2^*(s)}+1}} \;\big(\mu t^{2^*(s)-2}-\kappa \alpha t^{\beta}+\kappa\beta t^{\beta-2}-\lambda \big).
\end{align}
Define
\be\lab{2014-11-29-e1}
h(t):=\mu t^{2^*(s)-2}-\kappa \alpha t^{\beta}+\kappa\beta t^{\beta-2}-\lambda,
\ee
then we obtain that
$\{t>0:\; h(t)=0\}$ has at least $3$ solutions $\{t_1,t_2,t_3\}$ such that $0<t_1<t_2=t_0<t_3<\infty$.

\noindent{\bf Case 1:$\beta=2$ and $2\kappa-\lambda\leq 0$.}
For  this case, $h(t)=\mu t^{2^*(s)-2}-\kappa \alpha t^2+2\kappa -\lambda$. By Rolle's mean value theorem, $\{t>0:\;h'(t)=0\}$ has at least two solutions $\tilde{t}_1,\tilde{t}_2$ such that
$$ t_1<\tilde{t}_1<t_2=t_0<\tilde{t}_2<t_3.$$
Note that $$\{t>0:\;h'(t)=0\}=\{t>0:\; \mu[2^*(s)-2]t^{2^*(s)-4}-2\kappa \alpha=0\}.$$ In particular,  the set  $\{t>0:\; \mu[2^*(s)-2]t^{2^*(s)-4}-2\kappa \alpha=0\}$ has a unique solution if $2^*(s)\neq 4$, a contradiction. Hence, $2^*(s)=4$ and $\mu=2\kappa $. Recalling that $2^*(s)=\frac{2(N-s)}{N-2}, s\in (0,2)$, we obtain that $s=1, \alpha=2^*(s)-\beta=2$. Then
\be\lab{2014-11-29-e2}
g(t_0)=\frac{1+t_0^2}{[\lambda+2\kappa t_0^4+4\kappa t_0^2]^{\frac{1}{2}}}=\frac{1}{\sqrt{\lambda}},
\ee
which implies  that
$\displaystyle \lambda=\mu=2\kappa.$
It follows that $g(t)\equiv \frac{1}{\sqrt{\lambda}}$.
Hence, when $N=3,s=1,\alpha=\beta=2, \lambda=\mu=2\kappa$, the extremals of $S_{\alpha,\beta,\lambda,\mu}(\Omega)$ are given by \eqref{2014-11-29-e3}. In particular, $$\{(\phi,\psi)=\sqrt{\frac{\mu_s(\Omega)}{2\kappa(1+t^2)}}(U,tU):  t>0\}$$ are all the ground state solutions of
\be\lab{2014-11-29-e4}
\begin{cases}
-\Delta u-2\kappa \frac{|u|^2u}{|x|}=2\kappa \frac{uv^2}{|x|}\;\hbox{in}\;\Omega,\\
-\Delta v-2\kappa \frac{|v|^2v}{|x|}=2\kappa \frac{u^2v}{|x|}\;\hbox{in}\;\Omega,\\
\kappa>0, u,v\in D_{0}^{1,2}(\Omega),
\end{cases}
\ee
where $U$ is the ground state solution of
\eqref{2014-11-30-e1}.

\noindent{\bf Case 2:  $\beta>2$.} For this case, $h(t)=\mu t^{2^*(s)-2}-\kappa \alpha t^\beta+\kappa \beta t^{\beta-2}-\lambda$, similarly we see that the equation $h(t)=0$ ($t>0$) has at least three roots $t_1<t_2=t_0<t_3$. It follows that $\{h'(t),t>0\}$ has at least two roots $\tilde{t}_1$ and $\tilde{t}_2$, which implies that $p(t)=0$ $(t>0)$ has at least two solutions. Where $p(t)$ is defined by
\be\lab{2014-11-29-e5}
p(t):=\mu[2^*(s)-2]t^{2^*(s)-\beta}-\kappa\alpha\beta t^2+\kappa\beta(\beta-2).
\ee
A direct computation shows that
$p''(t)>0$ when $\alpha>2$.
Hence $p(t)=0$  $(t>0)$ could not have more than one solution, a contradiction.
If $\beta>2,\alpha=2$, we have $p(t)=\big[\mu[2^*(s)-2]-\kappa\alpha\beta\big]t^2+\kappa\beta(\beta-2)$,  which also has at most one positive root, a contradiction  too.

We note that for the case of $\mu>\lambda$, we will take $$ \tilde{g}(t):=g(\frac{1}{t})=\frac{1+t^2}{\left[\mu+\lambda t^{2^*(s)}+2^*(s)\kappa t^\alpha\right]^{\frac{2}{2^*(s)}}}$$
and the arguments above can repeated ($\beta$ is replaced by $\alpha$ now). We complete the proof.
\ep

\bc\lab{2014-11-30-cro1} Under the assumptions of  Theorem \ref{2014-11-30-th2}, there must hold
$$ S_{\alpha,\beta,\lambda,\mu}(\Omega)=\big(\max\{\lambda, \mu\}\big)^{-\frac{2}{2^*(s)}}\mu_s(\Omega).$$
\ec
\bp
For the special case $N=3,s=1,\alpha=\beta=2,\lambda=\mu=2\kappa$, a direct computation can deduce it.  For the other cases, if $\displaystyle S_{\alpha,\beta,\lambda,\mu}(\Omega)\neq \big(\max\{\lambda, \mu\}\big)^{-\frac{2}{2^*(s)}}\mu_s(\Omega),$
then there must hold that
$\displaystyle S_{\alpha,\beta,\lambda,\mu}(\Omega)<\big(\max\{\lambda, \mu\}\big)^{-\frac{2}{2^*(s)}}\mu_s(\Omega).$
By Theorem \ref{2014-11-26-th1} and Lemma \ref{2014-11-27-zl1}, $S_{\alpha,\beta,\lambda,\mu}(\Omega)$  can be achieved by some nontrivial extremal $(u,v)$, a contradiction to Theorem \ref{2014-11-30-th2}.
\ep

\subsection{Further results about cones}
Assume that $0<s<2,\alpha>1,\beta>1,\alpha+\beta=2^*(s)$.
Based  on the results of Theorem \ref{2014-11-26-th1}, we see that when $\Omega$ is a cone, $S_{\alpha,\beta, \lambda,\mu}(\Omega)$ is always achieved.
In this subsection, we always assume that $\Omega$ is a cone.  We shall investigate more properties about $S_{\alpha,\beta, \lambda,\mu}(\Omega)$
in terms of $\Omega$ such as monotonicity, continuity and intermediate value theorem. Let us begin with a remark.
\br\lab{2014-11-30-wr1}
Assume that $\Omega_1,\Omega_2$ are domains of $\R^N$ and $\Omega_1\subseteq \Omega_2$, then it is easy to see that $D_{0}^{1,2}(\Omega_1)\subseteq D_{0}^{1,2}(\Omega_2)$. Then by the definition of $S_{\alpha,\beta, \lambda,\mu}(\Omega)$ (see the formula \eqref{2014-11-26-e0}), we see that $S_{\alpha,\beta, \lambda,\mu}(\Omega_1)\geq S_{\alpha,\beta, \lambda,\mu}(\Omega_2)$.
\er

\bl\lab{2014-10-16-wl1}
Let $\{\Omega_n\}$ be a sequence of cones.
\begin{itemize}
\item[(i)] Assume  $\{\Omega_n\}$ is an  increasing sequence, i.e., $\Omega_n\subseteq \Omega_{n+1}$, then $$\lim_{n\rightarrow \infty}S_{\alpha,\beta, \lambda,\mu}(\Omega_n)=S_{\alpha,\beta, \lambda,\mu}(\lim_{n\rightarrow\infty}\Omega_n)=S_{\alpha,\beta, \lambda,\mu}(\Omega),$$
    where $$\Omega=\lim_{n\rightarrow\infty}\Omega_n=\cup_{n=1}^{\infty}\Omega_n.$$
\item[(ii)] Assume $\{\Omega_n\}$ is a decreasing sequence, i.e., $\Omega_n\supseteq \Omega_{n+1}$, then $$\lim_{n\rightarrow \infty}S_{\alpha,\beta, \lambda,\mu}(\Omega_n)=S_{\alpha,\beta, \lambda,\mu}(\lim_{n\rightarrow\infty}\Omega_n)=S_{\alpha,\beta, \lambda,\mu}(\Omega),$$
    where $$\Omega=\cap_{n=1}^{\infty}\Omega_n,$$
and we denote that $S_{\alpha,\beta, \lambda,\mu}(\Omega)=+\infty$ if $meas(\Omega)=0$.
\end{itemize}
\el
\bp
(i) By Remark \ref{2014-11-30-wr1}, we see that $\{S_{\alpha,\beta, \lambda,\mu}(\Omega_n)\}$ is a decreasing nonnegative sequence. Hence $\lim_{n\rightarrow\infty}S_{\alpha,\beta, \lambda,\mu}(\Omega_n)$ exists. Also by Remark \ref{2014-11-30-wr1} and $\Omega=\cup_{i=1}^{\infty}\Omega_i\supseteq \Omega_n,\forall\;n$, we have
\be\lab{2014-10-19-e1}
S_{\alpha,\beta, \lambda,\mu}(\Omega)\leq \lim_{n\rightarrow \infty}S_{\alpha,\beta, \lambda,\mu}(\Omega_n).
\ee
On the other hand, for any $\varepsilon>0$, $\exists\; (u_\varepsilon,v_\varepsilon)\in C_c^\infty(\Omega)\times C_c^\infty(\Omega)$ such that
\be\lab{2014-11-30-wze1}
\int_\Omega[|\nabla u_\varepsilon|^2 +|\nabla v_\varepsilon|^2]dx<S_{\alpha,\beta, \lambda,\mu}(\Omega)+\varepsilon
\ee
and
\be\lab{2014-11-30-wze2}
\int_\Omega \big(\lambda \frac{|u_\varepsilon|^{2^*(s)}}{|x|^s}+\mu \frac{|v_\varepsilon|^{2^*(s)}}{|x|^s}+2^*(s)\kappa \frac{|u_\varepsilon|^\alpha |v_\varepsilon|^\beta}{|x|^s}\big)dx=1.
\ee
Then there exists some $N_0$ large enough such that
\be\lab{2014-10-19-e3}
u_\varepsilon, v_\varepsilon\in C_c^\infty(\Omega_n)\;\hbox{for all}\;n\geq N_0.
\ee
Hence, by the definition of $S_{\alpha,\beta, \lambda,\mu}(\Omega_n)$ again, we have
\be\lab{2014-10-19-e4}
S_{\alpha,\beta, \lambda,\mu}(\Omega_n)<S_{\alpha,\beta, \lambda,\mu}(\Omega)+\varepsilon\;\hbox{for all}\;n\geq N_0.
\ee
Let $n$ go to infinity, we have
\be\lab{2014-10-19-e5}
\lim_{n\rightarrow \infty} S_{\alpha,\beta, \lambda,\mu}(\Omega_n)\leq S_{\alpha,\beta, \lambda,\mu}(\Omega)+\varepsilon.
\ee
By the arbitrariness  of $\varepsilon$, we have
\be\lab{2014-10-19-e6}
\lim_{n\rightarrow \infty}S_{\alpha,\beta, \lambda,\mu}(\Omega_n)\leq S_{\alpha,\beta, \lambda,\mu}(\Omega).
\ee
Now, \eqref{2014-10-19-e1} and \eqref{2014-10-19-e6} say that
\be\lab{2014-10-19-e7}
\lim_{n\rightarrow \infty}S_{\alpha,\beta, \lambda,\mu}(\Omega_n)=S_{\alpha,\beta, \lambda,\mu}(\lim_{n\rightarrow\infty}\Omega_n)=S_{\alpha,\beta, \lambda,\mu}(\Omega).
\ee
(ii) By Remark \ref{2014-11-30-wr1}, we see that $\{S_{\alpha,\beta, \lambda,\mu}(\Omega_n)\}$ is an increasing sequence. Let us denote $\displaystyle\bar{S}:=\lim_{n\rightarrow \infty}S_{\alpha,\beta, \lambda,\mu}(\Omega_n).$  For any $n$, let $(u_n, v_n)$ be the extremal function to $S_{\alpha,\beta, \lambda,\mu}(\Omega)$ by Theorem  \ref{2014-11-26-th1}. We can extend $u_n$ and $v_n$ by $0$ out side $\Omega_n$. By Remark \ref{2014-11-27-r1}, we have $\int_{\Omega_{1}} \big(\lambda \frac{|u_n|^{2^*(s)}}{|x|^s}+\mu \frac{|v_n|^{2^*(s)}}{|x|^s}+2^*(s)\kappa \frac{|u_n|^\alpha |v_n|^\beta}{|x|^s}\big)dx\equiv 1$ and
\begin{align}\lab{2014-10-19-e8}
&\frac{1}{2}=\int_{\Omega_{1}\cap B_1} \big(\lambda \frac{|u_n|^{2^*(s)}}{|x|^s}+\mu \frac{|v_n|^{2^*(s)}}{|x|^s}+2^*(s)\kappa \frac{|u_n|^\alpha |v_n|^\beta}{|x|^s}\big)dx\\
=&\int_{\Omega_{1}\backslash B_1} \big(\lambda \frac{|u_n|^{2^*(s)}}{|x|^s}+\mu \frac{|v_n|^{2^*(s)}}{|x|^s}+2^*(s)\kappa \frac{|u_n|^\alpha |v_n|^\beta}{|x|^s}\big)dx\nonumber
\end{align}

\noindent
{\bf Case 1--$meas(\Omega)=0$:} In this case, we shall prove that $\bar{S}=\infty$.
Now,we proceed by contradiction. If $\bar{S}<\infty$,
$\{u_n\}, \{v_n\}$ are  bounded sequences in $D_{0}^{1,2}(\Omega_{1})$.  Then up to a subsequence, we may assume that $u_n\rightharpoonup u, v_n\rightharpoonup v$ in $D_{0}^{1,2}(\Omega_{1})$ and $u_n\rightarrow u, v_n\rightarrow v$ a.e. in $\Omega_{1}$. Since $meas(\cap_{n=1}^{\infty}\Omega_n)=0$, we get $u=0, v=0$. On the other hand, by applying  the same argument as  Corollary \ref{2014-11-27-cro2}, we can obtain that $u_n\rightarrow u, v_n\rightarrow v$ in $L^{2^*(s)}(\Omega_1,\frac{dx}{|x|^s})$. Then we have
\be\lab{2014-10-19-e9}
\int_{\R^N} \big(\lambda \frac{|u|^{2^*(s)}}{|x|^s}+\mu \frac{|v|^{2^*(s)}}{|x|^s}+2^*(s)\kappa \frac{|u|^\alpha |v|^\beta}{|x|^s}\big)dx= 1,
\ee
a contradiction. Hence $\bar{S}=\infty$.\\

\noindent {\bf Case 2--$\Omega$ is a cone:} In this case, by Theorem \ref{2014-11-26-th1}, $S_{\alpha,\beta,\lambda,\mu}(\Omega)$ is well defined and can be achieved. Notice that for any $n$, we have $\Omega\subseteq \Omega_n$, by Remark \ref{2014-11-30-wr1} again, we have $S_{\alpha,\beta,\lambda,\mu}(\Omega_n)\leq S_{\alpha,\beta,\lambda,\mu}(\Omega)$. Hence
\be\lab{2014-10-19-bue3}
\bar{S}\leq S_{\alpha,\beta,\lambda,\mu}(\Omega).
\ee
 Thus, $\{u_n\},\{v_n\}$ are  bounded in $D_{0}^{1,2}(\Omega_1)$ for this case. Arguing as before, it is easy to see the weak limit $u\neq 0, v\neq 0$ and
\be\lab{2014-10-19-e10}
0<\int_{\Omega} \big(\lambda \frac{|u|^{2^*(s)}}{|x|^s}+\mu \frac{|v|^{2^*(s)}}{|x|^s}+2^*(s)\kappa \frac{|u|^\alpha |v|^\beta}{|x|^s}\big)dx\leq 1.
\ee
We claim that $(u, v)$ weakly  solves
\be\lab{2014-10-19-e11}
\begin{cases}
-\Delta u=\bar{S}\Big(\lambda \frac{1}{|x|^s}|u|^{2^*(s)-2}u+\kappa\alpha \frac{1}{|x|^s}|u|^{\alpha-2}u|v|^\beta\Big)\quad &\hbox{in}\;\Omega,\\
-\Delta v=\bar{S}\Big(\mu \frac{1}{|x|^s}|v|^{2^*(s)-2}v+\kappa\beta \frac{1}{|x|^s}|u|^{\alpha}|v|^{\beta-2}v\Big)\quad &\hbox{in}\;\Omega,\\
u\geq 0, v\geq 0, (u,v)\in \mathscr{D}:=D_{0}^{1,2}(\Omega)\times D_{0}^{1,2}(\Omega),
\end{cases}
\ee
Since $C_c^\infty(\Omega)$ is dense in $D_{0}^{1,2}(\Omega)$, we only need to prove that for all $(\phi,\psi)\in C_c^\infty(\Omega)\times C_c^\infty(\Omega)$, there holds
\begin{align}\lab{2014-10-19-e12}
&\int_\Omega \nabla u\cdot \nabla \phi+\nabla v\cdot \nabla\psi\\
=&\bar{S}\int_\Omega \big(\lambda \frac{|u|^{2^*(s)-1}\phi}{|x|^s}+\mu \frac{|v|^{2^*(s)-1}\psi}{|x|^s}+\kappa \alpha\frac{|u|^{\alpha-2}u\phi |v|^\beta}{|x|^s}+\kappa \beta\frac{|u|^{\alpha} |v|^{\beta-2}v\psi}{|x|^s}\big)dx.\nonumber
\end{align}
Now, let $(\phi,\psi)\in C_c^\infty(\Omega)\times C_c^\infty(\Omega)$ be fixed. Notice that $\Omega\subseteq\Omega_n$, we have $\phi,\psi\in D_{0}^{1,2}(\Omega_n), \forall\;n$. Then
\begin{align}\lab{2014-10-19-e13}
&\int_{\Omega_n} \nabla u_n\cdot \nabla \phi+\nabla v_n\cdot \nabla\psi\\
=&S_{\alpha,\beta,\lambda,\mu}(\Omega_n)
\int_{\Omega_n} \Big(\lambda \frac{|u_n|^{2^*(s)-1}\phi}{|x|^s}+\mu \frac{|v_n|^{2^*(s)-1}\psi}{|x|^s}\nonumber\\
&+\kappa \alpha\frac{|u_n|^{\alpha-2}u_n\phi |v_n|^\beta}{|x|^s}+\kappa \beta\frac{|u_n|^{\alpha} |v_n|^{\beta-2}v_n\psi}{|x|^s}\Big)dx.\nonumber
\end{align}
Since $supp(\phi),\;supp(\psi)\subseteq\Omega$, we have
\begin{align}\lab{2014-10-19-e14}
&\int_{\Omega_n} \nabla u_n\cdot \nabla \phi+\nabla v_n\cdot \nabla\psi\\
=&S_{\alpha,\beta,\lambda,\mu}(\Omega_n)
\int_{\Omega_n} \Big(\lambda \frac{|u_n|^{2^*(s)-1}\phi}{|x|^s}+\mu \frac{|v_n|^{2^*(s)-1}\psi}{|x|^s}\nonumber\\
&+\kappa \alpha\frac{|u_n|^{\alpha-2}u_n\phi |v_n|^\beta}{|x|^s}+\kappa \beta\frac{|u_n|^{\alpha} |v_n|^{\beta-2}v_n\psi}{|x|^s}\Big)dx\;\hbox{for all}\;n.\nonumber
\end{align}
Then apply the similar arguments as Corollary \ref{2014-11-27-cro2}, we have
\begin{align}\lab{2014-10-19-e15}
&\int_\Omega \nabla u\cdot \nabla \phi+\nabla v\cdot \nabla\psi\\
=&\bar{S}\int_\Omega \big(\lambda \frac{|u|^{2^*(s)-1}\phi}{|x|^s}+\mu \frac{|v|^{2^*(s)-1}\psi}{|x|^s}+\kappa \alpha\frac{|u|^{\alpha-2}u\phi |v|^\beta}{|x|^s}+\kappa \beta\frac{|u|^{\alpha} |v|^{\beta-2}v\psi}{|x|^s}\big)dx.\nonumber
\end{align}
Thereby the claim is proved.
By \eqref{2014-10-19-e10} and $2<2^*(s)$, we have
\begin{align}\lab{2014-10-19-e16}
&S_{\alpha,\beta,\lambda,\mu}(\Omega)
\leq \frac{\int_\Omega [|\nabla u|^2+|\nabla v|^2]dx}{\big(\int_{\Omega} \big(\lambda \frac{|u|^{2^*(s)}}{|x|^s}+\mu \frac{|v|^{2^*(s)}}{|x|^s}+2^*(s)\kappa \frac{|u|^\alpha |v|^\beta}{|x|^s}\big)dx\big)^{\frac{2}{2^*(s)}}}\\
\leq&  \frac{\int_\Omega [|\nabla u|^2+|\nabla v|^2]dx}{\int_{\Omega} \big(\lambda \frac{|u|^{2^*(s)}}{|x|^s}+\mu \frac{|v|^{2^*(s)}}{|x|^s}+2^*(s)\kappa \frac{|u|^\alpha |v|^\beta}{|x|^s}\big)dx}=\bar{S}.\nonumber
\end{align}
It follows from \eqref{2014-10-19-bue3} and \eqref{2014-10-19-e16}  that
$$ S_{\alpha,\beta,\lambda,\mu}(\Omega)=\bar{S}=\lim_{n\rightarrow\infty}S_{\alpha,\beta,\lambda,\mu}(\Omega_n)$$
and
$$ \int_{\Omega} \big(\lambda \frac{|u|^{2^*(s)}}{|x|^s}+\mu \frac{|v|^{2^*(s)}}{|x|^s}+2^*(s)\kappa \frac{|u|^\alpha |v|^\beta}{|x|^s}\big)dx=1.$$
Hence, $(u,v)$ is an extremal function of $S_{\alpha,\beta,\lambda,\mu}(\Omega)$. The proof is completed.
\ep

Define
\be\lab{2014-10-16-we1}
\underline{S}_{\alpha,\beta,\lambda,\mu}:=\inf\big\{S_{\alpha,\beta,\lambda,\mu}(\Omega):\Omega\;\hbox{is a cone properly contained in $\R^N\backslash\{0\}$}\big\}.
\ee
For any given unit versor $\nu$ in $\R^N$, let
\be\lab{2014-10-19-e17}
\Omega_\theta:=\{x\in \R^N:x\cdot \nu>|x|\cos\theta\},\;\;\;\; \theta\in (0,\pi].
\ee

\bd Assume $1<p<N, -\infty<t<N-p$,  we denote by $D_{t}^{1,p}(\Omega)$ the completion of $C_0^\infty(\Omega)$ with respect to the norm
\be\lab{zz=z}\|u\|:= \big(\int_\Omega\frac{|\nabla u|^p}{|x|^t}dx\big)^{\frac{1}{p}}\ee
\ed
Then we have the following result:
\bo\lab{2014-11-28-l1}
If $F$ is a closed subset of a $k-$dimensional subspace of $\R^N$ with $k<N-t-p$, then $D_{t}^{1,p}(\Omega)=D_{t}^{1,p}(\Omega\backslash F)$.
In particular,   $D_{t}^{1,p}(\R^N)=D_{t}^{1,p}(\R^N\backslash \{0\})$ provided $N-t-p>0$.
\eo
\bp
Without loss of generality, we assume that $\Omega=\R^N$. Notice that $C_c^\infty(\R^N\backslash F)\subseteq C_c^\infty(\R^N)$, by the definition, it is easy to see that $D_{t}^{1,p}(\R^N\backslash F)\subseteq D_{t}^{1,p}(\R^N)$.

On the other hand, for any $u\in D_{t}^{1,p}(\Omega)$, there exists a sequence $\{\varphi_n\}\subset C_c^\infty(\R^N)$ such that
\be\lab{2014-10-21-e1}
\|\varphi_n-u\|^p=\int_{\R^N}\frac{|\nabla (\varphi_n-u)|^p}{|x|^t}dx\rightarrow 0\;\hbox{as}\;n\rightarrow \infty.
\ee
If there exists a subsequence of $\{\varphi_{n_k}\}$ such that $supp(\varphi_{n_k})\cap F=\emptyset$, then $\{\varphi_{n_k}\}\subseteq C_c^\infty(\R^N\backslash F)$, and it follows that $u\in D_{t}^{1,p}(\R^N\backslash F)$ and the proof is completed.
Hence, we may assume that $supp(\varphi_{n})\cap F\neq\emptyset$ for any $n$ without loss of generality.
Now, for any fixed $n$, we may choose a suitable cutoff function $\chi_\delta$ such that $\chi_\delta=0$ in $F_\delta, \chi_\delta=0$ in $\R^N\backslash F_{2\delta}, \chi_\delta\in (0,1)$ in $(\R^N\backslash F_\delta)\cap F_{2\delta}, |\nabla \chi_\delta|\leq \frac{2}{\delta}$, where $F\subset \R^N$ and
$\displaystyle F_\delta:=\{x\in \R^N:dist(x, F)<\delta\}.$
We note that $\chi_\delta\varphi_n\in C_c^\infty(\R^N\backslash F)$ for all $\delta>0$.

Now, we estimate $\|\chi_\delta \varphi_n-\varphi_n\|^p$.
\begin{align*}
&\int_{\R^N}\frac{|\nabla (\chi_\delta \varphi_n-\varphi_n)|^p}{|x|^t}dx\\
=&\int_{\R^N}\frac{\big[|\nabla(\chi_\delta-1)^2\varphi_n^2|+2\nabla (\chi_\delta-1)\cdot \nabla \varphi_n (\chi_\delta-1)\varphi_n+(\chi_\delta-1)^2|\nabla \varphi_n|^2\big]^{\frac{p}{2}}}{|x|^t}dx\\
\leq&\int_{\R^N}\frac{\Big[2\big(|\nabla(\chi_\delta-1)^2\varphi_n^2|+(\chi_\delta-1)^2|\nabla \varphi_n|^2\big)\Big]^{\frac{p}{2}}}{|x|^t}dx\\
\leq&2^p\int_{supp(\varphi_n)\cap F_{2\delta}}\frac{|\nabla \chi_\delta|^p}{|x|^t}|\varphi_n^p|dx+2^p\int_{supp(\varphi_n)}|\chi_\delta-1|^p\frac{|\nabla \varphi_n|^p}{|x|^t}dx:=I+II.
\end{align*}
By the Lebesgue's dominated convergence theorem, it is easy to see that
\be\lab{2014-10-21-e2}
II=2^p\int_{supp(\varphi_n)}|\chi_\delta-1|^p\frac{|\nabla \varphi_n|^p}{|x|^t}dx\rightarrow 0\;\hbox{as}\;\delta\rightarrow 0.
\ee
Recalling that $|\nabla \chi_\delta|\leq \frac{2}{\delta}$, there exists some $c_n>0$ independent of $\delta$ such that
\be\lab{2014-10-21-e3}
I=2^p\int_{supp(\varphi_n)\cap F_{2\delta}}\frac{|\nabla \chi_\delta|^p}{|x|^\mu}|\varphi_n^p|dx\leq c_n (\frac{2}{\delta})^p \delta^{N-k-\mu}.
\ee
Hence, when $k<N-t-p$, we also have
\be\lab{2014-10-21-e4}
I\rightarrow 0\;\hbox{as}\;\delta\rightarrow 0.
\ee
Hence, we can take some $\delta_n$ small enough such that
\be\lab{2014-10-21-e5}
\|\chi_{\delta_n} \varphi_n-\varphi_n\|^p\leq \frac{1}{n}.
\ee
Now, we let $u_n:=\chi_{\delta_n} \varphi_n\in C_c^\infty(\R^N\backslash F)$, we see that $\|u_n-u\|^p\rightarrow 0$ as $n\rightarrow \infty$. Hence, $u\in D_{t}^{1,p}(\R^N\backslash F)$. Thus, $D_{t}^{1,p}(\R^N)\subseteq D_{t}^{1,p}(\R^N\backslash F)$.
Especially, when $N-t-p>0$, take $k=0$, we see that $D_{t}^{1,p}(\R^N)=D_{t}^{1,p}(\R^N\backslash \{0\})$ and the proof is completed.
\ep

\bl\lab{2014-10-21-l2}
$\displaystyle S_{\alpha,\beta,\lambda,\mu}(\Omega_\pi)=S_{\alpha,\beta,\lambda,\mu}(\R^N)\;\hbox{for}\;N\geq 4.$
\el
\bp
Take $F=\R^N\backslash \Omega_\pi$, $F_n:=F\cap \overline{B_n(0)}$. We note that $F_n$ is a closed subset of a $1-$dimensional subspace of $\R^N$ and $\displaystyle\lim_{n\rightarrow \infty}F_n=F$. Then by Proposition \ref{2014-11-28-l1}, $D_{0}^{1,2}(\R^N\backslash F_n)=D_{0}^{1,2}(\R^N)$ for any $n$. Then it follows that $D_{0}^{1,2}(\R^N\backslash F)=D_{0}^{1,2}(\R^N)$. That is,  $D_{0}^{1,2}(\Omega_\pi)=D_{0}^{1,2}(\R^N)$. Hence, $S_{p,a,b}(\Omega_\pi)=S_{p,a,b}(\R^N)$.
\ep

\bt\lab{2014-11-30-wth1}
For every $\tau\geq \underline{S}_{\alpha,\beta,\lambda,\mu}$, there exists a cone $\Omega$ in $\R^N$ such that $S_{\alpha,\beta,\lambda,\mu}(\Omega)=\tau$. Moreover, when $N\geq 4$, we have $\underline{S}_{\alpha,\beta,\lambda,\mu}=S_{\alpha,\beta,\lambda,\mu}(\R^N)$.
\et
\bp
Define a mapping $\tau:(0,\pi]\mapsto \R_+\cup\{0\}$ with $\tau(\theta)= S_{\alpha,\beta,\lambda,\mu}(\Omega_\theta)$. Then by Remark \ref{2014-11-30-wr1}, we see that the mapping $\tau$ is decreasing with related to $\theta$. Evidently,  $\tau(\theta)$ is continuous for a.e.  $\theta\in (0,\pi]$. Furthermore, we can strengthen the conclusion. Indeed, let $\theta\in (0,\pi)$ be fixed. For any $\theta_n\uparrow \theta$, by (i) of Lemma \ref{2014-10-16-wl1}, we have
\be\lab{2014-10-19-xe1}
\lim_{n\rightarrow \infty}\tau(\theta_n)=\tau(\theta).
\ee
On the other hand, for any $\theta_n\downarrow \theta$, by (ii) of Lemma \ref{2014-10-16-wl1}, we also obtain \eqref{2014-10-19-xe1}.
Hence, $\tau$ is continuous in $(0,\pi)$. In addition, $\tau(\theta_n)\downarrow S_{\alpha,\beta,\lambda,\mu}(\Omega_\pi)=\underline{S}_{\alpha,\beta,\lambda,\mu}$ as $\theta_n\uparrow \pi$ and $\tau(\theta_n)\uparrow +\infty$ as $\theta_n\downarrow 0$.
Especially, when $N\geq 4$, by Lemma \ref{2014-10-21-l2}, we have that
$$ \underline{S}_{\alpha,\beta,\lambda,\mu}=S_{\alpha,\beta,\lambda,\mu}(\Omega_\pi)=S_{\alpha,\beta,\lambda,\mu}(\R^N).$$
\ep

\subsection{Existsence of infinitely many sign-changing solutions}
In this subsection, we will study the existence of infinitely many  sign-changing  solutions as an application of Theorem \ref{2014-11-27-wth1}.
\bt\lab{2014-10-21-xth1}
Assume $s\in (0,2),\kappa>0, \alpha>1,\beta>1,\alpha+\beta=2^*(s)$. Let $\Omega_\theta$ be  defined by \eqref{2014-10-19-e17} for some fixed $\theta\in (0,\pi]$.   Suppose that   one of the following conditions holds:
\begin{itemize}
\item[$(a_1)$]$\lambda>\mu $ and  either $ 1<\beta<2$ or $\begin{cases}\beta=2\\ \kappa>\frac{\lambda}{2^*(s)} \end{cases}$;
\item[$(a_2)$] $\lambda=\mu$ and either $\min\{\alpha,\beta\}<2$ or $\begin{cases} \min\{\alpha,\beta\}=2,\\ \kappa>\frac{\lambda}{2^*(s)} \end{cases}$;
\item[$(a_3)$] $\lambda<\mu$ and either $ 1<\alpha<2$ or $\begin{cases}\alpha=2\\ \kappa>\frac{\mu}{2^*(s)} \end{cases}$.
\end{itemize}
 Then the problem
\be\lab{2014-10-21-xe1}
\begin{cases}
-\Delta u-\lambda \frac{1}{|x|^s}|u|^{2^*(s)-2}u=\kappa\alpha \frac{1}{|x|^s}|u|^{\alpha-2}u|v|^\beta\quad &\hbox{in}\;\Omega_\theta,\\
-\Delta v-\mu \frac{1}{|x|^s}|v|^{2^*(s)-2}v=\kappa\beta \frac{1}{|x|^s}|u|^{\alpha}|v|^{\beta-2}v\quad &\hbox{in}\;\Omega_\theta,\\
(u,v)\in \mathscr{D}:=D_{0}^{1,2}(\Omega_\theta)\times D_{0}^{1,2}(\Omega_\theta),
\end{cases}
\ee
possesses a sequence of sign changing solutions $\{(u_k, v_k)\}$ which are  distinct under the  modulo rotations around $\nu$.
Moreover,  their energies  $c_k$  satisfies
$\displaystyle \frac{c_k}{2^{k(N-1)}}\rightarrow +\infty$ as $k\rightarrow \infty$, where
\begin{align}\lab{yingqian-e1}
c_k:=&\frac{1}{2}\int_{\Omega_\theta}[|\nabla u_k|^2+|\nabla v_k|^2]dx\\
&-\frac{1}{2^*(s)}\int_{\Omega_\theta}\Big(\lambda \frac{|u_k|^{2^*(s)}}{|x|^s}+\mu \frac{|v_k|^{2^*(s)}}{|x|^s}+2^*(s)\kappa \frac{|u_k|^\alpha |v_k|^\beta}{|x|^s}\Big)dx.\nonumber
\end{align}
\et
\bp
The idea is inspired by \cite{CaldiroliMusina.1999}. We will construct a solution on $\Omega_\theta$ by gluing together suitable signed solutions corresponding to each sub-cone. Using the spherical coordinates, we write  $S^{n-1}=\{\theta_1,\cdots, \theta_{N-1}:\theta_i\in S^1, i=1,\cdots,N-1\}$. For any fixed $k\in \NN$, we set
$$ \Sigma_{j}^{(k)}=(\frac{j}{2^{k-1}}\theta-\theta, \frac{j+1}{2^{k-1}}\theta-\theta)\quad j=0,1,\cdots, 2^k-1$$
and for every choice of $(j_1,j_2,\cdots,j_{N-1})\in \{0,1,2,\cdots,2^k-1\}^{N-1}$,
$$ \Omega_{j_1,\cdots,j_{N-1}}^{(k)}:=\{x\in \Omega_\theta:\frac{x}{|x|}\in \Sigma_{j_1}^{(k)}\times\cdots\times \Sigma_{j_1,\cdots,j_{N-1}}^{(k)}\}.$$
Due to Theorem \ref{2014-11-27-wth1}, we can take $\big(u_{j_1,\cdots,j_{N-1}}^{(k)},v_{j_1,\cdots,j_{N-1}}^{(k)}\big)\in D_{0}^{1,2}(\Omega_{j_1,\cdots,j_{N-1}}^{(k)})\times D_{0}^{1,2}(\Omega_{j_1,\cdots,j_{N-1}}^{(k)})$ as the positive ground state solution to
\begin{align*}
&\begin{cases}
-\Delta u=S_{\alpha,\beta,\lambda,\beta}(\Omega_{j_1,\cdots,j_{N-1}}^{(k)})\Big(\lambda \frac{1}{|x|^s}|u|^{2^*(s)-2}u+\kappa\alpha \frac{1}{|x|^s}|u|^{\alpha-2}u|v|^\beta\Big)\\
-\Delta v=S_{\alpha,\beta,\lambda,\beta}(\Omega_{j_1,\cdots,j_{N-1}}^{(k)})\Big(\mu \frac{1}{|x|^s}|v|^{2^*(s)-2}v+\kappa\beta \frac{1}{|x|^s}|u|^{\alpha}|v|^{\beta-2}v\Big)
\end{cases}\;\hbox{in}\;\Omega_{j_1,\cdots,j_{N-1}}^{(k)}
\end{align*}
where $(u,v)\in \mathscr{D}:=D_{0}^{1,2}(\Omega_{j_1,\cdots,j_{N-1}}^{(k)})\times D_{0}^{1,2}(\Omega_{j_1,\cdots,j_{N-1}}^{(k)})$.
We can extend every $u_{j_1,\cdots,j_{N-1}}^{(k)}$ and $v_{j_1,\cdots,j_{N-1}}^{(k)}$ outside $\Omega_{j_1,\cdots,j_{N-1}}^{(k)}$ by $0$ and now we set
$$u^{(k)}:=\sum_{j_1=0}^{2^k-1}\cdots \sum_{j_{N-1}=0}^{2^k-1} (-1)^{j_1+\cdots+j_{N-1}} u_{j_1,\cdots,j_{N-1}}^{(k)}\in D_{0}^{1,2}(\Omega_\theta)$$
and
$$v^{(k)}:=\sum_{j_1=0}^{2^k-1}\cdots \sum_{j_{N-1}=0}^{2^k-1} (-1)^{j_1+\cdots+j_{N-1}} v_{j_1,\cdots,j_{N-1}}^{(k)}\in D_{0}^{1,2}(\Omega_\theta).$$
Notice that for any two different choices $(j_1,j_2,\cdots,j_{N-1})\neq (\tilde{j}_1,\tilde{j}_2,\cdots,\tilde{j}_{N-1})$, there exists some rotation $R\in O(\R^N)$, the orthogonal transformation, such that
$$ \Omega_{\tilde{j}_1,\cdots,\tilde{j}_{N-1}}^{(k)}=R\big(\Omega_{j_1,\cdots,j_{N-1}}^{(k)}\big).$$
Hence, we have
$$ S_{\alpha,\beta,\lambda,\beta}(\Omega_{j_1,\cdots,j_{N-1}}^{(k)})=S_{\alpha,\beta,\lambda,\beta}(\Omega_{\tilde{j}_1,\cdots,\tilde{j}_{N-1}}^{(k)}).$$
Then it follows that $(u^{(k)}, v^{(k)})$ weakly  solves
$$
\begin{cases}
-\Delta u=S_{\alpha,\beta,\lambda,\beta}(\Omega_{j_1,\cdots,j_{N-1}}^{(k)})\Big(\lambda \frac{1}{|x|^s}|u|^{2^*(s)-2}u+\kappa\alpha \frac{1}{|x|^s}|u|^{\alpha-2}u|v|^\beta\Big)\quad &\hbox{in}\;\Omega_{\theta},\\
-\Delta v=S_{\alpha,\beta,\lambda,\beta}(\Omega_{j_1,\cdots,j_{N-1}}^{(k)})\Big(\mu \frac{1}{|x|^s}|v|^{2^*(s)-2}v+\kappa\beta \frac{1}{|x|^s}|u|^{\alpha}|v|^{\beta-2}v\Big)\quad &\hbox{in}\;\Omega_{\theta},\\
  (u,v)\in \mathscr{D}:=D_{0}^{1,2}(\Omega_{\theta})\times D_{0}^{1,2}(\Omega_{\theta}).
\end{cases}
$$
 Noting that $2^*(s)>2$,   after a scaling, let
$$u_k:=\Big(S_{\alpha,\beta,\lambda,\beta}(\Omega_{j_1,\cdots,j_{N-1}}^{(k)})\Big)^{\frac{1}{2^*(s)-2}}u^{(k)},\;v_k:=\Big(S_{\alpha,\beta,\lambda,\beta}(\Omega_{j_1,\cdots,j_{N-1}}^{(k)})\Big)^{\frac{1}{2^*(s)-2}}v^{(k)}$$
then $(u_k,v_k)$ weakly solves \eqref{2014-10-21-xe1}.
By the construction of $u_k$ and $v_k$, it is easy to see that $(u_k,v_k)$ is  a sign changing solution  and  $\{(u_k,v_k)\}$  are distinct under  modulo rotations around $\nu$.
Moveover, we have
{\allowdisplaybreaks  \begin{align*}
c_k:=&\frac{1}{2}\int_{\Omega_\theta}[|\nabla u_k|^2+|\nabla v_k|^2]dx\nonumber\\
&-\frac{1}{2^*(s)}\int_{\Omega_\theta}\Big(\lambda \frac{|u_k|^{2^*(s)}}{|x|^s}+\mu \frac{|v_k|^{2^*(s)}}{|x|^s}+2^*(s)\kappa \frac{|u_k|^\alpha |v_k|^\beta}{|x|^s}\Big)dx\\
=&\big(\frac{1}{2}-\frac{1}{2^*(s)}\big)\int_{\Omega_\theta}[|\nabla u_k|^2+|\nabla v_k|^2]dx\\
=&\big(\frac{1}{2}-\frac{1}{2^*(s)}\big)\Big(S_{\alpha,\beta,\lambda,\mu}(\Omega_{j_1,\cdots,j_{N-1}}^{(k)})\Big)^{\frac{2}{2^*(s)-2}}\int_{\Omega_\theta}[|\nabla u^{(k)}|^2+|\nabla v^{(k)}|^2]dx\\
=&\big(\frac{1}{2}-\frac{1}{2^*(s)}\big)\Big(S_{\alpha,\beta,\lambda,\mu}(\Omega_{j_1,\cdots,j_{N-1}}^{(k)})\Big)^{\frac{2}{2^*(s)-2}}\sum_{j_1=0}^{2^k-1}\cdots \sum_{j_{N-1}=0}^{2^k-1}\int_{\Omega_{j_1,\cdots,j_{N-1}}^{(k)}}\nonumber\\
&[|\nabla u_{j_1,\cdots,j_{N-1}}^{(k)}|^2+|\nabla v_{j_1,\cdots,j_{N-1}}^{(k)}|^2]dx\\
=&\big(\frac{1}{2}-\frac{1}{2^*(s)}\big)\Big(S_{\alpha,\beta,\lambda,\mu}(\Omega_{j_1,\cdots,j_{N-1}}^{(k)})\Big)^{\frac{2}{2^*(s)-2}}
2^{k(N-1)}\\
&\quad \quad \cdot \int_{\Omega_{0,\cdots,0}} [|\nabla u_{0,\cdots,0}^{(k)}|^2+|\nabla v_{0,\cdots,0}^{(k)}|^2]dx\\
=&\big(\frac{1}{2}-\frac{1}{2^*(s)}\big)\Big(S_{\alpha,\beta,\lambda,\mu}(\Omega_{j_1,\cdots,j_{N-1}}^{(k)})\Big)^{\frac{2}{2^*(s)-2}}
2^{k(N-1)}.
\end{align*}}
By Lemma \ref{2014-10-16-wl1}, $\displaystyle S_{\alpha,\beta,\lambda,\mu}(\Omega_{j_1,\cdots,j_{N-1}}^{(k)})\rightarrow +\infty$ as $k\rightarrow \infty$. Recalling that $2^*(s)>2$ again, we obtain that
$\displaystyle\frac{c_k}{2^{k(N-1)}}\rightarrow +\infty.$
\ep

Apply the same argument as in the proof of Theorem \ref{2014-10-21-xth1}, we can obtain  the following result for the system
defined on $\R^N$:
\bt\lab{2014-10-21-xth2}
Assume $s\in (0,2), N\geq 4, \kappa>0, \alpha>1,\beta>1,\alpha+\beta=2^*(s)$. Suppose that one of the following conditions  holds:
\begin{itemize}
\item[$(a_1)$]$\lambda>\mu$ and either $ 1<\beta<2$ or $\begin{cases}\beta=2\\ \kappa>\frac{\lambda}{2^*(s)} \end{cases}$;
\item[$(a_2)$] $\lambda=\mu $ and either $  \min\{\alpha,\beta\}<2$ or $\begin{cases} \min\{\alpha,\beta\}=2,\\ \kappa>\frac{\lambda}{2^*(s)} \end{cases}$;
\item[$(a_3)$] $\lambda<\mu $ and either $ 1<\alpha<2$ or $\begin{cases}\alpha=2\\ \kappa>\frac{\mu}{2^*(s)} \end{cases}$.
\end{itemize}
Then the problem
\be\lab{2014-10-21-xe2}
\begin{cases}
-\Delta u-\lambda \frac{1}{|x|^s}|u|^{2^*(s)-2}u=\kappa\alpha \frac{1}{|x|^s}|u|^{\alpha-2}u|v|^\beta\quad &\hbox{in}\;\R^N,\\
-\Delta v-\mu \frac{1}{|x|^s}|v|^{2^*(s)-2}v=\kappa\beta \frac{1}{|x|^s}|u|^{\alpha}|v|^{\beta-2}v\quad &\hbox{in}\;\R^N,\\ (u,v)\in \mathscr{D}:=D_{0}^{1,2}(\R^N)\times D_{0}^{1,2}(\R^N),
\end{cases}
\ee
possesses a sequence of sign changing solutions $\{(u_k,v_k)\}$ whose energies  $\displaystyle\frac{c_k}{2^{k(N-1)}}\rightarrow +\infty$ as $k\rightarrow \infty$,  where
\begin{align*}
c_k:=&\frac{1}{2}\int_{\R^N}[|\nabla u_k|^2+|\nabla v_k|^2]dx\nonumber\\
&-\frac{1}{2^*(s)}\int_{\R^N}\Big(\lambda \frac{|u_k|^{2^*(s)}}{|x|^s}+\mu \frac{|v_k|^{2^*(s)}}{|x|^s}+2^*(s)\kappa \frac{|u_k|^\alpha |v_k|^\beta}{|x|^s}\Big)dx.
\end{align*}
\et
 \bp  It is a straightforward consequence of  Theorems \ref{2014-10-21-xth1} and \ref{2014-11-30-wth1}. We just keep in mind  that  when $N\geq 4$,   we have that
$S_{\alpha,\beta,\lambda,\mu}(\Omega_\pi)=S_{\alpha,\beta,\lambda,\mu}(\R^N).$
\ep
\br\lab{2014-10-21-xr1}
It is clear that this kind arguments  used in the proof of Theorem \ref{2014-10-21-xth1} can be adapted to  other cones  with suitable symmetry.
\er

\subsection{Further results on more  general domain $\Omega $ and on  the sharp constant  $S_{\alpha,\beta,\lambda,\mu}(\Omega)$}
\br\lab{2014-10-23-xr1}
Given a general open domain $\Omega$ (not necessarily  a cone), we let $S_{\alpha,\beta,\lambda,\mu}(\Omega)$ be defined by \eqref{2014-11-26-e0} if $\Omega\neq \emptyset$ and $S_{\alpha,\beta,\lambda,\mu}(\emptyset)=+\infty$.
In this subsection, we are concerned with    whether $S_{\alpha,\beta,\lambda,\mu}(\Omega)$ can be achieved  or not and we give some operational way to compute the value of $S_{\alpha,\beta,\lambda,\mu}(\Omega)$.
\er

 We note that $\Omega$ can be written as a union of a sequence of domains, $\Omega=\bigcup_{n=1}^{\infty}\Omega_n$. \bl\lab{2014-10-21-wl1}
Assume $\Omega_i\cap \Omega_j=\emptyset\;\forall\;i\neq j$, then we have
$\displaystyle S_{\alpha,\beta,\lambda,\mu}(\Omega)=\inf_{n\geq 1} S_{\alpha,\beta,\lambda,\mu}(\Omega_n).$
\el
\bp
For any $n$, since $\Omega_n\subseteq \Omega$, we have
$\displaystyle S_{\alpha,\beta,\lambda,\mu}(\Omega_n)\geq S_{\alpha,\beta,\lambda,\mu}(\Omega)\;\hbox{for all}\;n.$
Hence,
\be\lab{2014-11-30-wze3}
\inf_{n\geq 1}S_{\alpha,\beta,\lambda,\mu}(\Omega_n)\geq S_{\alpha,\beta,\lambda,\mu}(\Omega).
\ee
On the other hand, for any $\varepsilon>0$, there exists  a pair $(u, v)$ such that
\be\lab{2014-11-30-wze2}
\int_\Omega \big(\lambda \frac{|u|^{2^*(s)}}{|x|^s}+\mu \frac{|v|^{2^*(s)}}{|x|^s}+2^*(s)\kappa \frac{|u|^\alpha |v|^\beta}{|x|^s}\big)dx=1
\ee
 and
\be\lab{2014-10-21-we1}
\int_\Omega [|\nabla u|^2+|\nabla v|^2]dx<S_{\alpha,\beta,\lambda,\mu}(\Omega)+\varepsilon.
\ee
Set $u_n=u\big|_{\Omega_n}, v_n=v\big|_{\Omega_n}$, since $\Omega_i\cap \Omega_j=\emptyset$ for all $i\neq j$, we have
$(u_n,v_n)\in D_{0}^{1,2}(\Omega_n)\times D_{0}^{1,2}(\Omega_n)$ and $u=\sum_{n=1}^{\infty} u_n, v=\sum_{n=1}^{\infty} v_n$.
Then
{\allowdisplaybreaks  \begin{align*}
&\int_{\Omega_n}[|\nabla u_n|^2+|\nabla v_n|^2]dx\\
\geq&S_{\alpha,\beta,\lambda,\mu}(\Omega_n)\Big(\int_{\Omega_n} \big(\lambda \frac{|u_n|^{2^*(s)}}{|x|^s}+\mu \frac{|v_n|^{2^*(s)}}{|x|^s}+2^*(s)\kappa \frac{|u_n|^\alpha |v_n|^\beta}{|x|^s}\big)dx\Big)^{\frac{2}{2^*(s)}}\\
\geq&\big(\inf_{n\geq 1}S_{\alpha,\beta,\lambda,\mu}(\Omega_n)\big)\Big(\int_{\Omega_n} \big(\lambda \frac{|u_n|^{2^*(s)}}{|x|^s}+\mu \frac{|v_n|^{2^*(s)}}{|x|^s}+2^*(s)\kappa \frac{|u_n|^\alpha |v_n|^\beta}{|x|^s}\big)dx\Big)^{\frac{2}{2^*(s)}}\\
\geq&\big(\inf_{n\geq 1}S_{\alpha,\beta,\lambda,\mu}(\Omega_n)\big)\int_{\Omega_n} \big(\lambda \frac{|u_n|^{2^*(s)}}{|x|^s}+\mu \frac{|v_n|^{2^*(s)}}{|x|^s}+2^*(s)\kappa \frac{|u_n|^\alpha |v_n|^\beta}{|x|^s}\big)dx,
\end{align*}}
here we use  $\frac{2}{2^*(s)}<1$ and $\displaystyle \int_{\Omega_n} \big(\lambda \frac{|u_n|^{2^*(s)}}{|x|^s}+\mu \frac{|v_n|^{2^*(s)}}{|x|^s}+2^*(s)\kappa \frac{|u_n|^\alpha |v_n|^\beta}{|x|^s}\big)dx\leq 1$.
It follows that
\begin{align*}
&\int_{\Omega}[|\nabla u|^2+|\nabla v|^2]dx
=\sum_{n=1}^{\infty}\int_{\Omega_n}[|\nabla u_n|^2+|\nabla v_n|^2]dx\\
\geq&\big(\inf_{n\geq 1}S_{\alpha,\beta,\lambda,\mu}(\Omega_n)\big)\sum_{n=1}^{\infty}\int_{\Omega_n} \big(\lambda \frac{|u_n|^{2^*(s)}}{|x|^s}+\mu \frac{|v_n|^{2^*(s)}}{|x|^s}+2^*(s)\kappa \frac{|u_n|^\alpha |v_n|^\beta}{|x|^s}\big)dx\\
=&\big(\inf_{n\geq 1}S_{\alpha,\beta,\lambda,\mu}(\Omega_n)\big)\int_{\Omega} \big(\lambda \frac{|u|^{2^*(s)}}{|x|^s}+\mu \frac{|v|^{2^*(s)}}{|x|^s}+2^*(s)\kappa \frac{|u|^\alpha |v|^\beta}{|x|^s}\big)dx\\
=&\inf_{n\geq 1}S_{\alpha,\beta,\lambda,\mu}(\Omega_n).
\end{align*}
Hence,
$ \inf_{n\geq 1} S_{\alpha,\beta,\lambda,\mu}(\Omega_n)\leq S_{\alpha,\beta,\lambda,\mu}(\Omega)+\varepsilon.$
Therefore,
\be\lab{2014-11-30-wze4}
\inf_{n\geq 1} S_{\alpha,\beta,\lambda,\mu}(\Omega_n)\leq S_{\alpha,\beta,\lambda,\mu}(\Omega).
\ee
By \eqref{2014-11-30-wze3} and \eqref{2014-11-30-wze4}, we complete the proof.
\ep

Next, for $r>0$, we set
\be\lab{2014-10-23-e1}
\Omega_r:=\Omega\cap B_r,\quad
\Omega^r:=\Omega\backslash B_r.
\ee
By Remark \ref{2014-11-30-wr1}, we see that
the mapping $r\mapsto S_{\alpha,\beta,\lambda,\mu}(\Omega_r)$ is non increasing and the mapping $r\mapsto S_{\alpha,\beta,\lambda,\mu}(\Omega^r)$ is non decreasing. Hence, we can define
$\displaystyle S_{\alpha,\beta,\lambda,\mu}^{0}(\Omega):=\lim_{r\rightarrow 0}S_{\alpha,\beta,\lambda,\mu}(\Omega_r)$
and
$\displaystyle S_{\alpha,\beta,\lambda,\mu}^{\infty}(\Omega):=\lim_{r\rightarrow \infty}S_{\alpha,\beta,\lambda,\mu}(\Omega^r).$
\br\lab{2014-10-21-wbur1}
It is easy to see that
$S_{\alpha,\beta,\lambda,\mu}^{0}(\Omega)$ and $S_{\alpha,\beta,\lambda,\mu}^{\infty}(\Omega)$ still have the monotonicity property. Precisely, if $\Omega_1\subseteq \Omega_2$, then we have  $$ S_{\alpha,\beta,\lambda,\mu}^{0}(\Omega_1)\geq S_{\alpha,\beta,\lambda,\mu}^{0}(\Omega_2),\quad S_{\alpha,\beta,\lambda,\mu}^{\infty}(\Omega_1)\geq S_{\alpha,\beta,\lambda,\mu}^{\infty}(\Omega_2).$$
\er

\bt\lab{2014-10-21-wth1}
Assume that $s\in (0,2),\kappa>0, \alpha>1,\beta>1,\alpha+\beta=2^*(s)$ and $\Omega$ is an open  domain of $\R^N$. Let $\{(u_n,v_n)\}$ be a minimizing sequence, i.e.,
$$\int_\Omega \big(\lambda \frac{|u_n|^{2^*(s)}}{|x|^s}+\mu \frac{|v_n|^{2^*(s)}}{|x|^s}+2^*(s)\kappa \frac{|u_n|^\alpha |v_n|^\beta}{|x|^s}\big)dx\equiv 1$$
and $$ \int_\Omega [|\nabla u_n|^2+|\nabla v_n|^2]dx\rightarrow S_{\alpha,\beta,\lambda,\mu}(\Omega)\;\hbox{as}\;n\rightarrow \infty.$$ Then one of the following  cases happens:
\begin{itemize}
\item[(a)] There exists some $(u,v)\in D_{0}^{1,2}(\Omega)\times D_{0}^{1,2}(\Omega)$ such that $u_n\rightarrow u, v_n\rightarrow v$ strongly in $D_{0}^{1,2}(\Omega)$ and $(u,v)\neq (0,0)$ is an extremal function of $S_{\alpha,\beta,\lambda,\mu}(\Omega)$;
\item[(b)]Going to a subsequence if necessary, we set
$$\eta:=\lim_{r\rightarrow 0}\lim_{n\rightarrow \infty}\int_{\Omega_r}\big(\lambda \frac{|u_n|^{2^*(s)}}{|x|^s}+\mu \frac{|v_n|^{2^*(s)}}{|x|^s}+2^*(s)\kappa \frac{|u_n|^\alpha |v_n|^\beta}{|x|^s}\big)dx.$$
Then
\be\lab{2014-10-21-we2}
S_{\alpha,\beta,\lambda,\mu}^{0}(\Omega)\eta^{\frac{2}{2^*(s)}}+S_{\alpha,\beta,\lambda,\mu}^{\infty}(\Omega)(1-\eta)^{\frac{2}{2^*(s)}}\leq S_{\alpha,\beta,\lambda,\mu}(\Omega).
\ee
\end{itemize}
\et
\bp
It is easy to see that $(u,v)$ is an  extremal function of $S_{\alpha,\beta,\lambda,\mu}(\Omega)$ if and only if $(u,v)$ is a ground state solution of
\be\lab{2014-10-22-e1}
\begin{cases}
-\Delta u=S_{\alpha,\beta,\lambda,\beta}(\Omega)\Big(\lambda \frac{1}{|x|^s}|u|^{2^*(s)-2}u+\kappa\alpha \frac{1}{|x|^s}|u|^{\alpha-2}u|v|^\beta\Big)\quad &\hbox{in}\;\Omega,\\
-\Delta v=S_{\alpha,\beta,\lambda,\beta}(\Omega)\Big(\mu \frac{1}{|x|^s}|v|^{2^*(s)-2}v+\kappa\beta \frac{1}{|x|^s}|u|^{\alpha}|v|^{\beta-2}v\Big)\quad &\hbox{in}\;\Omega,\\
u\geq 0, v\geq 0, (u,v)\in \mathscr{D}:=D_{0}^{1,2}(\Omega)\times D_{0}^{1,2}(\Omega),
\end{cases}
\ee
 Since $\{(u_n,v_n)\}$ is a minimizing sequence, we have that $\{(u_n,v_n)\}$ is a bounded $(PS)_{d}$ sequence with $d=(\frac{1}{2}-\frac{1}{2^*(s)})S_{\alpha,\beta,\lambda,\beta}(\Omega)$. Without loss of generality, we assume that $u_n\rightharpoonup u, v_n\rightharpoonup v$ in $D_{0}^{1,2}(\Omega)$ and $u_n\rightarrow u, v_n\rightarrow v$ a.e. in $\Omega$. Then it is easy to see that  $(u,v)$ is a weak solution to  \eqref{2014-10-22-e1} and
 $$ 0\leq \int_\Omega \big(\lambda \frac{|u|^{2^*(s)}}{|x|^s}+\mu \frac{|v|^{2^*(s)}}{|x|^s}+2^*(s)\kappa \frac{|u|^\alpha |v|^\beta}{|x|^s}\big)dx\leq 1,$$
 $$
\int_\Omega [|\nabla u|^2+|\nabla v|^2]dx\leq S_{\alpha,\beta,\lambda,\mu}(\Omega).$$

\noindent
{\bf Case 1:} If $(u,v)\neq (0,0)$, we shall prove that $(a)$ happens.
In this case, $(u,v)$ is a nontrivial solution or semi-trivial solution of \eqref{2014-10-22-e1}.
We claim $$\int_\Omega \big(\lambda \frac{|u|^{2^*(s)}}{|x|^s}+\mu \frac{|v|^{2^*(s)}}{|x|^s}+2^*(s)\kappa \frac{|u|^\alpha |v|^\beta}{|x|^s}\big)dx=1.$$  If not,
$0<\int_\Omega \big(\lambda \frac{|u|^{2^*(s)}}{|x|^s}+\mu \frac{|v|^{2^*(s)}}{|x|^s}+2^*(s)\kappa \frac{|u|^\alpha |v|^\beta}{|x|^s}\big)dx<1$. Then
\begin{align*}
&S_{\alpha,\beta,\lambda,\mu}(\Omega)
=\frac{\int_\Omega [|\nabla u|^2+|\nabla v|^2]dx}{\int_\Omega \big(\lambda \frac{|u|^{2^*(s)}}{|x|^s}+\mu \frac{|v|^{2^*(s)}}{|x|^s}+2^*(s)\kappa \frac{|u|^\alpha |v|^\beta}{|x|^s}\big)dx}\\
>&\frac{\int_\Omega [|\nabla u|^2+|\nabla v|^2]dx}{\Big(\int_\Omega \big(\lambda \frac{|u|^{2^*(s)}}{|x|^s}+\mu \frac{|v|^{2^*(s)}}{|x|^s}+2^*(s)\kappa \frac{|u|^\alpha |v|^\beta}{|x|^s}\big)dx\Big)^{\frac{2}{2^*(s)}}}
\geq S_{\alpha,\beta,\lambda,\mu}(\Omega),
\end{align*}
a contradiction. Hence, $$ \int_\Omega \big(\lambda \frac{|u|^{2^*(s)}}{|x|^s}+\mu \frac{|v|^{2^*(s)}}{|x|^s}+2^*(s)\kappa \frac{|u|^\alpha |v|^\beta}{|x|^s}\big)dx=1,$$
$$\int_\Omega [|\nabla u|^2+|\nabla v|^2]dx=S_{\alpha,\beta,\lambda,\mu}(\Omega).$$ That is, $(u,v)$ is an extremal function of $S_{\alpha,\beta,\lambda,\mu}(\Omega)$. Note that  $\|u_n-u\|=\|u_n\|-\|u\|+o(1), \|v_n-v\|=\|v_n\|-\|v\|+o(1)$, we see that $u_n\rightarrow u, v_n\rightarrow v$ in $D_{0}^{1,2}(\Omega)$.

\noindent
{\bf Case 2: } If $(u,v)=(0,0)$, we shall prove that (b) happens. The idea is similar to the proof of Lemma \ref{2014-11-26-xl5} and Corollary \ref{2014-11-27-cro2}.
Going to a subsequence if necessary, we set
$$ \Lambda_0:=\lim_{r\rightarrow 0}\lim_{n\rightarrow \infty}\int_{\Omega_r}[|\nabla u_n|^2+|\nabla v_n|^2]dx$$
and
$$ \Lambda^\infty:=\lim_{r\rightarrow \infty}\lim_{n\rightarrow \infty}\int_{\Omega^r}[|\nabla u_n|^2+|\nabla v_n|^2]dx.$$
Recalling the Rellich-Kondrachov compact theorem and $2<2^*(s)<2^*:=\frac{2N}{N-2}$, we have $(u_n,v_n)\rightarrow (0,0)$ in $L_{loc}^t(\Omega)\times L_{loc}^t(\Omega)$ for all $1<t<2^*$. Hence,
\be\lab{2014-10-22-e2}
\int_{\tilde{\Omega}} \frac{|u_n|^{2^*(s)}}{|x|^{s}}dx=o(1), \int_{\tilde{\Omega}} \frac{|v_n|^{2^*(s)}}{|x|^{s}}dx=o(1)\;
\ee
for any bounded domain $\tilde{\Omega}\subset \Omega$ such that $0\not\in \overline{\tilde{\Omega}}$.
Hence, we obtain that
\be\lab{2014-10-22-e3}
\lim_{r\rightarrow \infty}\lim_{n\rightarrow \infty}\int_{\Omega^r}\big(\lambda \frac{|u_n|^{2^*(s)}}{|x|^s}+\mu \frac{|v_n|^{2^*(s)}}{|x|^s}+2^*(s)\kappa \frac{|u_n|^\alpha |v_n|^\beta}{|x|^s}\big)dx=1-\eta.
\ee
Similar to the formula  \eqref{2014-11-27-e9}, we have
\be\lab{2014-10-22-e4}
S_{\alpha,\beta,\lambda,\mu}^{0}(\Omega)\eta^{\frac{2}{2^*(s)}}\leq \Lambda_0.
\ee
and similar to the formula  \eqref{2014-11-27-e10}, we have
\be\lab{2014-10-22-e5}
S_{\alpha,\beta,\lambda,\mu}^{\infty}(\Omega)\big(1-\eta\big)^{\frac{2}{2^*(s)}}\leq \Lambda^\infty.
\ee
Then by \eqref{2014-10-22-e4} and \eqref {2014-10-22-e5}, we have
\be\lab{2014-10-22-e6}
S_{\alpha,\beta,\lambda,\mu}^{0}(\Omega)\eta^{\frac{2}{2^*(s)}}+S_{\alpha,\beta,\lambda,\mu}^{\infty}(\Omega)\big(1-\eta\big)^{\frac{2}{2^*(s)}}\leq \Lambda_0+\Lambda^\infty\leq S_{\alpha,\beta,\lambda,\mu}(\Omega).
\ee
\ep

Theorem  \ref{2014-10-21-wth1} is a kind of concentration compactness principle, original spirit we refer to \cite{Lions.1985a}. Basing on Theorem \ref{2014-10-21-wth1}, we have the following using result:
\bc\lab{2014-10-22-cor1}
Assume that $s\in (0,2),\kappa>0, \alpha>1,\beta>1,\alpha+\beta=2^*(s)$  and $\Omega$ is an open  domain of $\R^N$. Then we always have
\be\lab{2014-10-22-e7}
S_{\alpha,\beta,\lambda,\mu}(\Omega)\leq min\big\{S_{\alpha,\beta,\lambda,\mu}^{0}(\Omega),\;\; S_{\alpha,\beta,\lambda,\mu}^{\infty}(\Omega)\big\}.
\ee
Moreover, if $S_{\alpha,\beta,\lambda,\mu}(\Omega)< min\{S_{\alpha,\beta,\lambda,\mu}^{0}(\Omega),S_{\alpha,\beta,\lambda,\mu}^{\infty}(\Omega)\}$, then $S_{\alpha,\beta,\lambda,\mu}(\Omega)$ can be achieved.
\ec
\bp
We note that $\Omega_r\subseteq \Omega, \Omega^r\subseteq \Omega$, by the monotonicity property, for any $r>0$, we have
$$ S_{\alpha,\beta,\lambda,\mu}(\Omega)\leq S_{\alpha,\beta,\lambda,\mu}(\Omega_r), \;\; \;S_{\alpha,\beta,\lambda,\mu}(\Omega)\leq S_{\alpha,\beta,\lambda,\mu}(\Omega^r)$$
which deduce  \eqref{2014-10-22-e7}.
Moreover, if $S_{\alpha,\beta,\lambda,\mu}(\Omega)< min\{S_{\alpha,\beta,\lambda,\mu}^{0}(\Omega),S_{\alpha,\beta,\lambda,\mu}^{\infty}(\Omega)\}$, then
\eqref{2014-10-21-we2} will never meet. Hence, only case (a) of Theorem\ref{2014-10-21-wth1} happens. Thus, $S_{\alpha,\beta,\lambda,\mu}(\Omega)$ is achieved.
\ep

Furthermore, we have the following result:
\bc\lab{2014-10-22-cor2}
Assume that $s\in (0,2),\kappa>0, \alpha>1,\beta>1,\alpha+\beta=2^*(s)$ and $\Omega$ is an open  domain of $\R^N$. If $S_{\alpha,\beta,\lambda,\mu}(\Omega)$ can not be achieved, then at least one of the following holds:
\begin{itemize}
\item[(i)]$S_{\alpha,\beta,\lambda,\mu}(\Omega)=S_{\alpha,\beta,\lambda,\mu}(\Omega_r)$ for $\forall\; r>0$.
\item[(ii)]$S_{\alpha,\beta,\lambda,\mu}(\Omega)=S_{\alpha,\beta,\lambda,\mu}(\Omega^r)$ for $\forall\; r>0$.
\end{itemize}
\ec
\bp
When $S_{\alpha,\beta,\lambda,\mu}(\Omega)$ is not attained,
by Corollary  \ref{2014-10-22-cor1}, we have
$$ S_{\alpha,\beta,\lambda,\mu}(\Omega)= min\{S_{\alpha,\beta,\lambda,\mu}^{0}(\Omega),S_{\alpha,\beta,\lambda,\mu}^{\infty}(\Omega)\}.$$
Without loss of generality, we assume that
$S_{\alpha,\beta,\lambda,\mu}^{0}(\Omega)\leq S_{\alpha,\beta,\lambda,\mu}^{\infty}(\Omega)$, then we have
$$ S_{\alpha,\beta,\lambda,\mu}(\Omega)=S_{\alpha,\beta,\lambda,\mu}^{0}(\Omega).$$
Next, we shall prove that case (i) holds. If not, assume that there exists some $r_0>0$ such that $S_{\alpha,\beta,\lambda,\mu}(\Omega)\neq S_{\alpha,\beta,\lambda,\mu}(\Omega_{r_0})$, then by the monotonicity property, we have
$S_{\alpha,\beta,\lambda,\mu}(\Omega)<S_{\alpha,\beta,\lambda,\mu}(\Omega_{r_0})\leq S_{\alpha,\beta,\lambda,\mu}^0(\Omega_{r_0})=S_{\alpha,\beta,\lambda,\mu}^{0}(\Omega)$, a contradiction.
\ep
\br\lab{2014-10-22-xr1}
We note that the  inverse  statement  of Corollary \ref{2014-10-22-cor2} is not true. For example, by Theorem \ref{2014-11-26-th1}, when $\Omega$ is a cone, $S_{\alpha,\beta,\lambda,\mu}(\Omega)$ is attained provided $1<\alpha,1<\beta,\alpha+\beta=2^*(s), s\in (0,2),\kappa>0$. However, we still have the following result.
\er
\bl\lab{2014-10-22-xl1}
Assume that $s\in (0,2),\kappa>0, \alpha>1,\beta>1,\alpha+\beta=2^*(s)$, then
$$S_{\alpha,\beta,\lambda,\mu}(\Omega)=S_{\alpha,\beta,\lambda,\mu}(\Omega_r)=S_{\alpha,\beta,\lambda,\mu}(\Omega^r)\;\hbox{for any}\;r>0\;\hbox{if $\Omega$ is a cone of $\R^N$}.$$
In particular,
$$ S_{\alpha,\beta,\lambda,\mu}(\R^N)=S_{\alpha,\beta,\lambda,\mu}(B_r)=S_{\alpha,\beta,\lambda,\mu}(\R^N\backslash B_r)\;\hbox{for any}\;r>0.$$
Further,
$S_{\alpha,\beta,\lambda,\mu}(\Omega)=S_{\alpha,\beta,\lambda,\mu}(\R^N)$  provided that either $\Omega$ is a general open  domain   with $0\in \Omega$ or $\Omega$ is an exterior domain.
Furthermore, let $A$ be a cone of $\R^N$,
and $\Omega=A\backslash F$, where $F$ is a closed subset of $A$ such that $0\not\in F$ or $F$ is bounded. Then $$ S_{\alpha,\beta,\lambda,\mu}(\Omega)=S_{\alpha,\beta,\lambda,\mu}(A).$$
\el
\bp
We only prove $S_{\alpha,\beta,\lambda,\mu}(\Omega)=S_{\alpha,\beta,\lambda,\mu}(\Omega_r)$ and the others  are  similar. By the monotonicity property, we see that $S_{\alpha,\beta,\lambda,\mu}(\Omega)\leq S_{\alpha,\beta,\lambda,\mu}(\Omega_r)$. Next, we shall prove the opposite inequality.
By Theorem \ref{2014-11-26-th1}, $S_{\alpha,\beta,\lambda,\mu}(\Omega)$ is attained. Now, let $(u,v)\in D_{0}^{1,2}(\Omega)\times D_{0}^{1,2}(\Omega)$ be an  extremal function of $S_{\alpha,\beta,\lambda,\mu}(\Omega)$. We also let $\chi_\rho(x)\in C_c^\infty(\R^N)$ be a cut-off function such that $\chi_\rho(x)\equiv 1$ in $B_{\frac{\rho}{2}}$, $\chi_\rho(x)\equiv 0$ in $\R^N\backslash B_{\rho}, |\nabla \chi_\rho(x)|\leq \frac{4}{\rho}$ and define $\phi_\rho:=\chi_\rho(x)u(x), \psi_\rho:=\chi_\rho(x)v(x)\in C_c^\infty(\Omega_{\rho})$. It is easy to see that
$$\int_{\Omega}[|\nabla \phi_\rho|^2+|\nabla \psi_\rho|^2]dx\rightarrow \int_\Omega [|\nabla u|^2+|\nabla v|^2]dx=S_{\alpha,\beta,\lambda,\mu}(\Omega)$$
and
\begin{align*}
&\int_\Omega \big(\lambda \frac{|\phi_\rho|^{2^*(s)}}{|x|^s}+\mu \frac{|\psi_\rho|^{2^*(s)}}{|x|^s}+2^*(s)\kappa \frac{|\phi_\rho|^\alpha |\psi_\rho|^\beta}{|x|^s}\big)dx\\
\xlongrightarrow{\rho\rightarrow +\infty}&\int_\Omega \big(\lambda \frac{|u|^{2^*(s)}}{|x|^s}+\mu \frac{|v|^{2^*(s)}}{|x|^s}+2^*(s)\kappa \frac{|u|^\alpha |v|^\beta}{|x|^s}\big)dx
=1.
\end{align*}
Then $\forall\;\varepsilon>0$, there exists some $\rho_0>0$ such that
$$\frac{\int_{\Omega}[|\nabla \phi_{\rho_0}|^2+|\nabla \psi_{\rho_0}|^2]dx}{\Big(\int_\Omega \big(\lambda \frac{|\phi_{\rho_0}|^{2^*(s)}}{|x|^s}+\mu \frac{|\psi_{\rho_0}|^{2^*(s)}}{|x|^s}+2^*(s)\kappa \frac{|\phi_{\rho_0}|^\alpha |\psi_{\rho_0}|^\beta}{|x|^s}\big)dx\Big)^{\frac{2}{2^*(s)}}}\leq S_{\alpha,\beta,\lambda,\mu}(\Omega)+\varepsilon.$$
Now, consider $\tilde{u}_r(x):=\phi_{\rho_0}(\frac{\rho_0}{r}x), \tilde{v}_r(x):=\psi_{\rho_0}(\frac{\rho_0}{r}x)\in C_c^\infty(\Omega_r)$ and
\begin{align*}
&S_{p,a,b}(\Omega_r)\leq \frac{\int_{\Omega_r}[|\nabla \tilde{u}_r(x)|^2+|\nabla \tilde{v}_r(x)|^2]dx}{\Big(\int_{\Omega_r} \big(\lambda \frac{|\tilde{u}_r(x)|^{2^*(s)}}{|x|^s}+\mu \frac{|\tilde{v}_r(x)|^{2^*(s)}}{|x|^s}+2^*(s)\kappa \frac{|\tilde{u}_r(x)|^\alpha |\tilde{v}_r(x)|^\beta}{|x|^s}\big)dx\Big)^{\frac{2}{2^*(s)}}}\\
=&\frac{\int_{\Omega}[|\nabla \phi_{\rho_0}|^2+|\nabla \psi_{\rho_0}|^2]dx}{\Big(\int_\Omega \big(\lambda \frac{|\phi_{\rho_0}|^{2^*(s)}}{|x|^s}+\mu \frac{|\psi_{\rho_0}|^{2^*(s)}}{|x|^s}+2^*(s)\kappa \frac{|\phi_{\rho_0}|^\alpha |\psi_{\rho_0}|^\beta}{|x|^s}\big)dx\Big)^{\frac{2}{2^*(s)}}}
\leq S_{\alpha,\beta,\lambda,\mu}(\Omega)+\varepsilon.
\end{align*}
Hence,
$$ S_{\alpha,\beta,\lambda,\mu}(\Omega_r)\leq S_{\alpha,\beta,\lambda,\mu}(\Omega).$$
Especially, take $\Omega=\R^N$, we see that
$$ S_{\alpha,\beta,\lambda,\mu}(\R^N)=S_{\alpha,\beta,\lambda,\mu}(B_r)=S_{\alpha,\beta,\lambda,\mu}(\R^N\backslash B_r)\;\hbox{for any}\;r>0.$$
Hence, when $0\in \Omega$, then there exists some $r>0$ such that $B_r\subset \Omega$, then
$$ S_{\alpha,\beta,\lambda,\mu}(B_r)\geq S_{\alpha,\beta,\lambda,\mu}(\Omega)\geq S_{\alpha,\beta,\lambda,\mu}(\R^N)=S_{\alpha,\beta,\lambda,\mu}(B_r).$$
If $\Omega$ is an exterior domain, there exists some $r>0$ such that $(\R^N\backslash B_r)\subset \Omega$, by the monotonicity property again, we have
$$ S_{\alpha,\beta,\lambda,\mu}(\R^N\backslash B_r)\geq S_{\alpha,\beta,\lambda,\mu}(\Omega)\geq S_{\alpha,\beta,\lambda,\mu}(\R^N)=S_{\alpha,\beta,\lambda,\mu}(\R^N\backslash B_r).$$
Furthermore, $A$ is a cone and $\Omega=A\backslash F\subset A$, then we have
$S_{\alpha,\beta,\lambda,\mu}(\Omega)\geq S_{\alpha,\beta,\lambda,\mu}(A).$
If $0\not\in F$, then there exists some $r>0$ such that $\Omega_r=A_r$, then
$$ S_{\alpha,\beta,\lambda,\mu}(\Omega)\leq S_{\alpha,\beta,\lambda,\mu}(\Omega_r)=S_{\alpha,\beta,\lambda,\mu}(A_r)=S_{\alpha,\beta,\lambda,\mu}(A).$$
Hence, we have
$S_{\alpha,\beta,\lambda,\mu}(\Omega)=S_{\alpha,\beta,\lambda,\mu}(A).$
If $F$ is bounded, then there exists  some $r>0$ such that $\Omega^r=A^r$, then
$$ S_{\alpha,\beta,\lambda,\mu}(\Omega)\leq S_{\alpha,\beta,\lambda,\mu}(\Omega^r)=S_{\alpha,\beta,\lambda,\mu}(A^r)=S_{\alpha,\beta,\lambda,\mu}(A),$$
it follows that
$S_{\alpha,\beta,\lambda,\mu}(\Omega)=S_{\alpha,\beta,\lambda,\mu}(A).$
\ep

To search the results on   general domains,    we introduce the following marks:
$$\mathcal{A}^0(\Omega):=\{A: \hbox{$A$ is a cone and there exists some $r>0$ such that $\Omega_r\subseteq A$}\}$$
and
$$\mathcal{A}^\infty(\Omega):=\{A: \hbox{$A$ is a cone and there exists some $r>0$ such that $\Omega^r\subseteq A$}\}.$$
Notice that $\R^N\in\mathcal{A}^0(\Omega)\cap \mathcal{A}^\infty(\Omega) $, $\mathcal{A}^0(\Omega)\neq \emptyset, \mathcal{A}^\infty(\Omega)\neq \emptyset$. Then we can define
$$\tilde{S}_{\alpha,\beta,\lambda,\mu}^{0}(\Omega):=\sup \{S_{\alpha,\beta,\lambda,\mu}(A):A\in \mathcal{A}^0(\Omega)\}$$ and
$$ \tilde{S}_{\alpha,\beta,\lambda,\mu}^{\infty}(\Omega):=\sup \{S_{\alpha,\beta,\lambda,\mu}(A):A\in \mathcal{A}^\infty(\Omega)\}.$$
\bl\lab{2014-10-22-wl1}
$\displaystyle\tilde{S}_{\alpha,\beta,\lambda,\mu}^{0}(\Omega)\leq S_{\alpha,\beta,\lambda,\mu}^{0}(\Omega), \quad \tilde{S}_{\alpha,\beta,\lambda,\mu}^{\infty}(\Omega)\leq S_{\alpha,\beta,\lambda,\mu}^{\infty}(\Omega).$
\el
\bp
We only prove $\tilde{S}_{\alpha,\beta,\lambda,\mu}^{0}(\Omega)\leq S_{\alpha,\beta,\lambda,\mu}^{0}(\Omega)$.
For any $\varepsilon>0$, there exists some $A\in \mathcal{A}^0(\Omega)$
such that
\be\lab{2014-10-22-we1}
\tilde{S}_{\alpha,\beta,\lambda,\mu}^{0}(\Omega)-\varepsilon <S_{\alpha,\beta,\lambda,\mu}(A).
\ee
By the definition of $\mathcal{A}^0(\Omega)$, there exists some $r>0$ such that
$\Omega_r\subset A.$
Then by the monotonicity property, we have
\be\lab{2014-10-22-we2}
S_{\alpha,\beta,\lambda,\mu}(A)\leq S_{\alpha,\beta,\lambda,\mu}(\Omega_r)\leq S_{\alpha,\beta,\lambda,\mu}^{0}(\Omega).
\ee
By \eqref{2014-10-22-we1}, \eqref{2014-10-22-we2} and the arbitrariness  of $\varepsilon$, we obtain that
$\tilde{S}_{\alpha,\beta,\lambda,\mu}^{0}(\Omega)\leq S_{\alpha,\beta,\lambda,\mu}^{0}(\Omega).$
\ep

\br\lab{2014-10-22-wr1}
By Corollary \ref{2014-10-22-cor1}, if we can prove that $$S_{\alpha,\beta,\lambda,\mu}(\Omega)<\min\{S_{\alpha,\beta,\lambda,\mu}^{0}(\Omega), S_{\alpha,\beta,\lambda,\mu}^{\infty}(\Omega)\},$$ then $S_{\alpha,\beta,\lambda,\mu}(\Omega)$ is attained. Since a lot of properties about cones have been studied in Section 7 (subsections 7.1-7.6),
Lemma \ref{2014-10-22-wl1} supplies a useful way to compute $\min\{S_{\alpha,\beta,\lambda,\mu}^{0}(\Omega), S_{\alpha,\beta,\lambda,\mu}^{\infty}(\Omega)\}$. Here, we prefer to give some examples.
\er

\noindent {\bf Example 1:} If  $\Omega$ is bounded with $0\not\in \bar{\Omega}$, then by $S_{\alpha,\beta,\lambda,\mu}(\emptyset)=+\infty$, we see that $S_{\alpha,\beta,\lambda,\mu}^{0}(\Omega)=S_{\alpha,\beta,\lambda,\mu}^{\infty}(\Omega)=+\infty$. Hence, $$ S_{\alpha,\beta,\lambda,\mu}(\Omega)<\min\{S_{\alpha,\beta,\lambda,\mu}^{0}(\Omega), S_{\alpha,\beta,\lambda,\mu}^{\infty}(\Omega)\}$$ and $S_{\alpha,\beta,\lambda,\mu}(\Omega)$ is attained which can also deduce by Theorem \ref{2014-11-26-th1}.

The following examples are also given by Caldiroli, Paolo and Musina, Roberta \cite{CaldiroliMusina.1999}, when they study the case of $p=2, b=0, a\in (-1,0)$. What interesting is that the similar results still hold for our case (a slight modification on Example 3).

\noindent{\bf Example 2:}   Assume $0\in \Omega$ is a cusp point, i.e., there exists a unit versor $\nu$ such that  $\forall\theta\in (0,\pi), \exists\;r_\theta>0$ such that $\Omega_{r_\theta}\subseteq \Omega_\theta$.
Notice that
$$ S_{\alpha,\beta,\lambda,\mu}^{0}(\Omega)=S_{\alpha,\beta,\lambda,\mu}^{0}(\Omega_{r_\theta})\geq S_{\alpha,\beta,\lambda,\mu}(\Omega_{r_\theta}).$$
On the other hand,
$$ S_{\alpha,\beta,\lambda,\mu}(\Omega_{r_\theta})\geq S_{\alpha,\beta,\lambda,\mu}(\Omega_\theta)\rightarrow +\infty\;\hbox{as}\;\theta\rightarrow 0.$$
Hence,
$S_{\alpha,\beta,\lambda,\mu}^{0}(\Omega)=+\infty.$

\noindent{\bf Example 3:}  Let
$\Omega=\Lambda\times \R^{N-k}, 1\leq k<N$,  where $\Lambda$ is an open bounded domain of $\R^k$. Then there exists some $r>0$ such that $\Lambda\subset B_r^k,$ the ball in  $\R^k$ with radial $r$. Now, we let $$ A_n:=\big\{(tx', tx'')\in \R^k\times \R^{N-k}:  t>0,x'\in B_r^k, |x''|_{N-k}\geq n\big\}, $$  then it is easy to see that $\{A_n\}$ are cones such that $A_n\supseteq A_{n+1}, \forall\;n$ and $\displaystyle \cap_{n=1}^{\infty}A_n\subset \{0\}\times \R^{N-k}$. Thus, $meas(\bigcap_{n=1}^{\infty}A_n)=0$.
By (ii) of Lemma \ref{2014-10-16-wl1}, $\displaystyle \lim_{n\rightarrow \infty}S_{\alpha,\beta,\lambda,\mu}(A_n)=+\infty$. Define $\tilde{\Omega}=B_r^k\times \R^{N-k}$, then it is easy to see that $\Omega\subset \tilde{\Omega}$. Moreover, for any $n$, there exists some $r_n>\sqrt{r^2+n^2}>0$ such that $\tilde{\Omega}^{r_n}\subset A_n$,  where $\tilde{\Omega}^{r_n}$ is defined by \eqref{2014-10-23-e1}. Indeed, for any $x=(x_1,x_2)\in \tilde{\Omega}\backslash A_n$, we have $|x_1|_{k}<r$ and $\frac{|x_2|_{N-k}}{|x_1|_k}\leq \frac{n}{r}$, thus $|x_2|_{N-k}\leq n$. Then it follows that $|x|_N\leq \sqrt{r^2+n^2}$. Hence, $\Omega^{r_n}\subset \tilde{\Omega}^{r_n}\subset A_n$.  Then by the monotonicity property we have $S_{\alpha,\beta,\lambda,\mu}(\Omega^{r_n})\geq S_{\alpha,\beta,\lambda,\mu}(\tilde{\Omega}^{r_n})\geq S_{\alpha,\beta,\lambda,\mu}(A_n)\rightarrow +\infty$ as $n\rightarrow +\infty$. Hence $S_{\alpha,\beta,\lambda,\mu}^{+\infty}(\Omega)=\infty$.

\bl\lab{2014-10-23-xl1}
Assume that $s\in (0,2),\kappa>0, \alpha>1,\beta>1,\alpha+\beta=2^*(s)$, and $\Omega=\Lambda\times \R^{N-k}, 1\leq k<N$,  where $\Lambda$ is an open bounded domain of $\R^k$  with $0\not\in \bar{\Lambda}$. Then $S_{\alpha,\beta,\lambda,\mu}(\Omega)$ is attained.
\el
\bp
By Example 3, we see that $S_{\alpha,\beta,\lambda,\mu}^\infty(\Omega)=+\infty$. By $0\not\in \bar{\Lambda}$, we have $S_{\alpha,\beta,\lambda,\mu}^{0}(\Omega)=+\infty$.
Then by Corollary \ref{2014-10-22-cor1}, we obtain the  conclusion.
\ep

Based  on the result of Lemma \ref{2014-10-22-xl1} and the  maximum principle,  we can obtain the following interesting results.
\bc\lab{2014-10-23-cro1}
Assume that $N\geq 3, s\in (0,2),\kappa>0, \alpha>1,\beta>1,\alpha+\beta=2^*(s)$, and let $\Omega$ be a general open domain of $\R^N$.
\begin{itemize}
\item[(i)]If $0\in \Omega$, then $S_{\alpha,\beta,\lambda,\mu}(\Omega)$ is not attained unless $\Omega=\R^N$;
\item[(ii)]If $\Omega$ is an exterior domain, then $S_{\alpha,\beta,\lambda,\mu}(\Omega)$ is not attained unless $\Omega=\R^N$;
\item[(iii)]If $\Omega=A\cup U$,  where $U$ is an open bounded set with $0\not\in \bar{U}$ and $A$ is a cone of $\R^N$, then $$ S_{\alpha,\beta,\lambda,\mu}(\Omega)<S_{\alpha,\beta,\lambda,\mu}(A)=S_{\alpha,\beta,\lambda,\mu}^{0}(\Omega)=S_{\alpha,\beta,\lambda,\mu}^{\infty}(\Omega)$$ and $S_{\alpha,\beta,\lambda,\mu}(\Omega)$ is attained.
\end{itemize}
\ec
\bp
For the case of $(i)$ and (ii), by Lemma \ref{2014-10-22-xl1}, we see that $$ S_{\alpha,\beta,\lambda,\mu}(\Omega)=S_{\alpha,\beta,\lambda,\mu}(\R^N).$$ Then by the maximum principle, it is easy to see that $S_{\alpha,\beta,\lambda,\mu}(\Omega)$ is not attained unless $\Omega=\R^N$. For the case of $(iii)$, since $U$ is bounded and $0\not\in \bar{U}$, there exists $r_1,r_2>0$ such that $\Omega_{r_1}=A_{r_1}, \Omega^{r_2}=A^{r_2}$. Hence $S_{\alpha,\beta,\lambda,\mu}^{0}(\Omega)=S_{\alpha,\beta,\lambda,\mu}^{\infty}(\Omega)=S_{\alpha,\beta,\lambda,\mu}(A)$.
If $S_{\alpha,\beta,\lambda,\mu}(\Omega)<\min\{S_{\alpha,\beta,\lambda,\mu}^{0}(\Omega), \;\; S_{\alpha,\beta,\lambda,\mu}^{\infty}(\Omega)\}$, then by Corollary \ref{2014-10-22-cor1}, $S_{\alpha,\beta,\lambda,\mu}(\Omega)$ is attained and the proof is completed. If not, by Corollary \ref{2014-10-22-cor1} again, $S_{\alpha,\beta,\lambda,\mu}(\Omega)=\min\{S_{\alpha,\beta,\lambda,\mu}^{0}(\Omega), S_{\alpha,\beta,\lambda,\mu}^{\infty}(\Omega)\}$,
Hence, $S_{\alpha,\beta,\lambda,\mu}(\Omega)=S_{\alpha,\beta,\lambda,\mu}(A)$. By maximum principle again, $S_{\alpha,\beta,\lambda,\mu}(A)$ is not attained, a contradiction with Theorem \ref{2014-11-26-th1}.
\ep

\bt\lab{2014-10-24-th1}
Assume that $N\geq 3, s\in (0,2),\kappa>0, \alpha>1,\beta>1,\alpha+\beta=2^*(s)$ and  $\Omega\subset\R^N$ is an open bounded domain.
If $0\in \Omega$, then $S_{\alpha,\beta,\lambda,\mu}(\Omega)$ is not attained.
\et
\bp
Since $\Omega$ is bounded, there exists some $r>0$ such that $\Omega\subset B_r(0)$. When $0\in \Omega$, by Lemma \ref{2014-10-22-xl1}, we have
$$ S_{\alpha,\beta,\lambda,\mu}(\Omega)=S_{\alpha,\beta,\lambda,\mu}(B_r(0))=S_{\alpha,\beta,\lambda,\mu}(\R^N).$$
By way of negation, assume that $S_{\alpha,\beta,\lambda,\mu}(\Omega)$ is attained and let $(u,v) \in D_{0}^{1,2}(\Omega)\times D_{0}^{1,2}(\Omega)$ be an extremal function. We may assume that $u\geq 0, v\geq 0$. By extending $u$ and $v$ outside $\Omega$ by $0$, then we see that $(u,v)$ is also an extremal function of $S_{\alpha,\beta,\lambda,\mu}(B_r(0)$. Hence, $(u,v)\neq (0,0)$ is a nonnegative weak solution to the following problem:
\be\lab{2014-10-24-e1}
\begin{cases}
-\Delta u=S_{\alpha,\beta,a,b}(B_r(0))\Big(\lambda \frac{1}{|x|^s}|u|^{2^*(s)-2}u+\kappa\alpha \frac{1}{|x|^s}|u|^{\alpha-2}u|v|^\beta\Big)\quad &\hbox{in}\;B_r(0),\\
-\Delta v=S_{\alpha,\beta,a,b}(B_r(0))\Big(\mu \frac{1}{|x|^s}|v|^{2^*(s)-2}v+\kappa\beta \frac{1}{|x|^s}|u|^{\alpha}|v|^{\beta-2}v\Big)\quad &\hbox{in}\;B_r(0),\\
u\geq 0, v\geq 0, (u,v)\in \mathscr{D}:=D_{0}^{1,2}(B_r(0))\times D_{0}^{1,2}(B_r(0)).
\end{cases}
\ee
On the other hand, since $B_r(0)$ is a star-shaped domain, by Proposition \ref{2013-10-26-prop1} (see the formula \eqref{2013-10-26-e1}), problem\eqref{2014-10-24-e1} has no nontrivial solution even semi-trivial solution, a contradiction.
Hence, we know  that $S_{\alpha,\beta,\lambda,\mu}(\Omega)$ is not attained.
\ep

\bc\lab{2014-10-24-cro1}
Assume that $N\geq 3, s\in (0,2),\kappa>0, \alpha>1,\beta>1,\alpha+\beta=2^*(s)$, if there exist  some $r_1,r_2>0$ and $\theta\in (0,\pi]$ such that
$$\big(\Omega_\theta \cap B_{r_1}(0)\big)\subsetneq\Omega\subsetneq \big(\Omega_\theta \cap B_{r_2}(0)\big),$$
then $S_{\alpha,\beta,\lambda,\mu}(\Omega)=S_{\alpha,\beta,\lambda,\mu}(\Omega_\theta)$ and is not attained.
\ec
\bp
The proof is similar to that of Theorem \ref{2014-10-24-th1}, we omit the details.
\ep


\section{The case of $s_1\neq s_2\in (0, 2)$}
\renewcommand{\theequation}{8.\arabic{equation}}
\renewcommand{\theremark}{8.\arabic{remark}}
\renewcommand{\thedefinition}{8.\arabic{definition}}
In this section, we study the case of  $s_1\neq s_2\in (0, 2)$.  By constructing a new approximation, the existence of positive ground state solution to the   system  \eqref{P} will be obtained, including   the  regularity and decay estimation.

Assume $\Omega$ is cone.
Define $\displaystyle U_\lambda:=\big(\frac{\mu_{s_1}(\Omega)}{\lambda}\big)^{\frac{1}{2^*(s_1)-2}}U,$
 where $U$ is a ground state solution to  the following problem:
\be\lab{2014-12-10-e1}
\begin{cases}
-\Delta u=\mu_{s_1}(\Omega) \frac{u^{2^*(s_1)-1}}{|x|^{s_1}}\;\hbox{in}\;\Omega,\\
u=0\;\;\hbox{on}\;\partial \Omega.
\end{cases}
\ee
Set
\be\lab{2014-12-9-xxe1}
\eta_{1,0}:=\inf_{v\in \Xi_0} \|v\|^2,\quad\quad \eta_{2,0}:=\inf_{u\in \Theta_0} \|u\|^2,
\ee
where
\be\lab{2014-12-9-xxe2}
\Xi_0:=\Big\{v\in D_{0}^{1,2}(\Omega):\;\int_\Omega \frac{1}{|x|^{s_2}}|U_\lambda|^{2^*(s_2)-2}|v|^2 dx=1\Big\},
\ee
\be\lab{2014-12-9-xxe2+0}
\Theta_0:=\Big\{u\in D_{0}^{1,2}(\Omega):\;\int_\Omega \frac{1}{|x|^{s_2}}|U_\mu|^{2^*(s_2)-2}|u|^2 dx=1\Big\}.
\ee
The corresponding energy functional of the  problem  \eqref{P}  is defined as
\begin{align}\lab{zz=620}
\Phi_0(u,v)&=\frac{1}{2}a(u,v)-\frac{1}{2^*(s_1)}b(u,v)
-\kappa c(u,v)
\end{align}
for all $(u,v)\in\mathscr{D}$, where
\be\lab{2014-12-3-e2-0}
\begin{cases}
a(u,v):=\|(u,v)\|_{\mathscr{D}}^{2},\\
b(u,v):=\lambda \int_{\Omega}\frac{|u|^{2^*(s_1)}}{|x|^{s_1}}dx+\mu \int_{\Omega}\frac{|v|^{2^*(s_1)}}{|x|^{s_1}}dx,\\
c(u,v):=\int_{\Omega}\frac{1}{|x|^{s_2}}|u|^\alpha |v|^\beta dx.
\end{cases}
\ee
Here comes our main result in this section:
{\allowdisplaybreaks  \bt\lab{2014-12-9-mainth1}
Assume that $s_1,s_2\in (0,2), \lambda, \mu\in (0,+\infty), \kappa>0,\alpha>1,\beta>1, \alpha+\beta= 2^*(s_2)$. Suppose  that  one of the following holds:
  \begin{itemize}
\item[$(a_1)$]   $\lambda>\mu$ and either $1<\beta<2$ or $\begin{cases}\beta=2\\ \kappa>\frac{\eta_{1,0}}{2^*(s_2)} \end{cases}$;
\item[$(a_2)$] $\lambda=\mu $ and either $ \min\{\alpha,\beta\}<2$ or $\begin{cases} \min\{\alpha,\beta\}=2,\\ \kappa>\frac{\eta_{1,0}}{2^*(s_2)}= \frac{\eta_{2,0}}{2^*(s_2)}\end{cases}$;
\item[$(a_3)$] $\lambda<\mu$ and either $  1<\alpha<2$ or $\begin{cases}\alpha=2\\ \kappa>\frac{\eta_{2,0}}{2^*(s_2)} \end{cases}$.
\end{itemize}
Then problem \eqref{P} possesses a positive ground state  solution $(u_0,v_0)$ such that
$$\Phi_0(u_0,v_0)< \Big[\frac{1}{2}-\frac{1}{2^*(s_1)}\Big]\big(\mu_{s_1}(\Omega)\big)^{\frac{2^*(s_1)}{2^*(s_1)-2}} \big(\max\{\lambda,\mu\}\big)^{-\frac{2}{2^*(s_1)-2}}.$$
Moreover, we have the following regularity and decay propositions:
\begin{itemize}
\item[$(b_1)$] if $0<\max\{s_1,s_2\}<\frac{N+2}{N}$, then $u_0,v_0\in C^2(\overline{\Omega})$;
\item[$(b_2)$] if $\max\{s_1,s_2\}=\frac{N+2}{N}$,  then $u_0,v_0\in C^{1,\gamma}({\Omega})$ for all $0<\gamma<1$;
\item[$(b_3)$] if $\max\{s_1,s_2\}>\frac{N+2}{N}$,  then $u_0,v_0\in C^{1,\gamma}({\Omega})$ for all $0<\gamma<\frac{N(2-\max\{s_1,s_2\})}{N-2}$.
\end{itemize}
When   $\Omega$ is a cone with $0\in \partial\Omega$ (e.g., $\Omega=\R_+^N$), then there exists a constant $C$ such that
$$ |u_0(x)|, |v_0(x)|\leq C (1+|x|^{-(N-1)}), |\nabla u_0(x)|, |\nabla v_0(x)|\leq C|x|^{-N}.$$
When $\Omega=\R^N$,
$$ |u_0(x)|, |v_0(x)|\leq C (1+|x|^{-N}), \quad \quad |\nabla u_0(x)|, |\nabla  v_0(x)|\leq C|x|^{-N-1}.$$
In particular, if $\Omega=\R_+^N$, then  $(u_0(x), v_0(x))$ is axially symmetric with respect to the $x_N$-axis, i.e., $$ \big(u_0(x), v_0(x)\big)=\big(u_0(x',x_N),
v_0(x',x_N)\big)=\big(u_0(|x'|,x_N), v_0(|x'|,x_N)\big).$$
\et}

\br\lab{main-r1}
The regularity, symmetry results and the decay estimation we  have established in Section 3 of the    present paper. Therefore, in the current section  we only need to focus on the existence of  the positive ground state solution.
\er

\subsection{Approximation}

When $s_1\neq s_2$, the nonlinearities are not homogeneous any more which make the problem much tough.  Here we have to choose a different approximation  to  the original problem in the same domain, i.e., we consider the following problem:
\be\lab{zz-620-1}
\begin{cases}
-\Delta u-\lambda \frac{1}{|x|^{s_1}}|u|^{2^*(s_1)-2}u=\kappa\alpha a_\varepsilon(x)|u|^{\alpha-2}u|v|^\beta\quad &\hbox{in}\;\Omega,\\
-\Delta v-\mu \frac{1}{|x|^{s_1}}|v|^{2^*(s_1)-2}v=\kappa\beta a_\varepsilon(x)|u|^{\alpha}|v|^{\beta-2}v\quad &\hbox{in}\;\Omega,\\
\kappa>0,(u,v)\in \mathscr{D}:=D_{0}^{1,2}(\Omega)\times D_{0}^{1,2}(\Omega),
\end{cases}
\ee
where
\be\lab{2014-12-9-we1}
a_\varepsilon(x):=\begin{cases}
\frac{1}{|x|^{s_2-\varepsilon}}\quad &\hbox{for}\;|x|<1 \\
\frac{1}{|x|^{s_2+\varepsilon}}&\hbox{for}\;|x|\geq 1
\end{cases}
\quad\quad \hbox{for }\;    \varepsilon\in [0,s_2).
\ee
Under some proper assumptions on   $\alpha,\beta,\lambda,\mu $ and $ \kappa>0$, we shall  prove the existence of the positive ground state solution $(u_\varepsilon, v_\varepsilon)$  to  \eqref{zz-620-1}  with a well-dominated  energy (see Theorem \ref{2014-12-6-th1} below).
Finally, we can approach an existence result of  \eqref{P}.

The corresponding energy functional of problem   \eqref{zz-620-1} is defined as
\begin{align}\lab{zz-620-f}
\Phi_\varepsilon(u,v)&=\frac{1}{2}a(u,v)-\frac{1}{2^*(s_1)}b(u,v)
-\kappa c_\varepsilon(u,v)
\end{align}
for all $(u,v)\in\mathscr{D}$, where $ a(u,v)$ and $b(u,v)$ are defined in (\ref{2014-12-3-e2-0}) and
\be\lab{zz-71-1}c_\varepsilon(u,v):=\int_{\Omega}a_\varepsilon(x)|u|^\alpha |v|^\beta dx, \ee which   is decreasing by $\varepsilon$.
Consider the corresponding Nehari manifold
$\displaystyle \mathcal{N}_\varepsilon:=\{(u,v)\in \mathscr{D}\backslash (0,0):\;  J_\varepsilon(u,v)=0\}$,
where
\begin{align}\lab{zz-71-2}
J_\varepsilon(u,v):=&\langle \Phi'_\varepsilon(u,v), (u,v)\rangle
=a(u,v)-b(u,v)-\kappa(\alpha+\beta)c_\varepsilon(u,v).
\end{align}

\bl\lab{2014-12-5-l1}
Assume $s_1,s_2\in (0,2), \lambda, \mu\in (0,+\infty), \kappa>0,\alpha>1,\beta>1$ and $ \alpha+\beta= 2^*(s_2)$. Let $\varepsilon\in [0,s_2)$,
then for any $(u,v)\in \mathscr{D}\backslash \{(0,0)\}$, there exists a unique $t=t_{(\varepsilon,u,v)}>0$ such that $(tu,tv)\in \mathcal{N}_\varepsilon$. Moreover,
$\mathcal{N}_\varepsilon$ is closed and bounded away from $(0,0)$.  Further, $t=t_{(\varepsilon,u,v)}$ is increasing by $\varepsilon$.
\el
\bp
The existence and uniqueness of $t=t_{(\varepsilon,u,v)}$ and that $\mathcal{N}_\varepsilon$ is closed and bounded away from $0$, we refer to Lemma
\ref{2014-11-23-xl1}.  Now, we prove that $t=t_{(\varepsilon,u,v)}$ is increasing by $\varepsilon$.
Assume that $0\leq \varepsilon_1<\varepsilon_2<s_2$, then we see that there exists a  unique $t_1$ and $t_2$ such that
\be\lab{2014-12-5-e1}
J_{\varepsilon_1} (t_1u,t_1v)=J_{\varepsilon_2} (t_2u,t_2v)=0.
\ee
Recalling that $c_\varepsilon(u,v)$ is decreasing by $\varepsilon$, we see that $J_\varepsilon(u,v)$ is increasing by $\varepsilon$. Hence,
\be\lab{2014-12-5-e2}
J_{\varepsilon_2} (t_1u,t_1v)\geq J_{\varepsilon_1} (t_1u,t_1v)=0.
\ee
If $J_{\varepsilon_2} (t_1u,t_1v)=0$, by the uniqueness, we obtain that $t_2=t_1$.
If $J_{\varepsilon_2} (t_1u,t_1v)>0$, noting that $J_{\varepsilon_2} (tu,tv)\rightarrow -\infty$ as $t\rightarrow +\infty$, there exists some $t_*>t_1$ such that $J_{\varepsilon_2} (t_*u,t_*v)=0$.
Then by the uniqueness again, we see that $t_2=t_*>t_1$.
Hence, we always have $t_2\geq t_1$ and we note that $t_2>t_1$ when $uv\not\equiv 0$.
\ep

\vskip 0.26in
Define
\be\lab{2014-12-5-e3}
c_\varepsilon:=\inf_{(u,v)\in \mathcal{N}_\varepsilon}\Phi_\varepsilon(u,v), \quad\quad
\delta_\varepsilon:=\inf_{(u,v)\in \mathcal{N}_\varepsilon}\sqrt{\|u\|^2+\|v\|^2}.
\ee
We have the following results:
\bl\lab{2014-12-5-l4}
$\delta_\varepsilon$ is increasing by $\varepsilon\in [0,s_2)$, i.e., $\delta_0\leq \delta_{\varepsilon_1}\leq \delta_{\varepsilon_2}$ provided $0\leq \varepsilon_1<\varepsilon_2<s_2$.
\el
\bp
For any $(u,v)\neq (0,0)$, set $\phi=\frac{u}{\sqrt{\|u\|^2+\|v\|^2}}, \psi=\frac{v}{\sqrt{\|u\|^2+\|v\|^2}}$. By Lemma \ref{2014-12-5-l1}, there exists $0<t_1\leq t_2$ such that $(t_1\phi, t_1\psi)\in \mathcal{N}_{\varepsilon_1}$ and $(t_2\phi,t_2\psi)\in \mathcal{N}_{\varepsilon_2}$. Hence, we obtain that $\delta_\varepsilon$ is increasing by $\varepsilon\in [0,s_2)$.
\ep

\br\lab{2014-12-6-bur1}
Set $s_{max}:=\max\{s_1,s_2\}$, it is easy to prove that for any $(u,v)\in \mathcal{N}_\varepsilon$, we have
\be\lab{2014-12-6-bue1}
\Phi_\varepsilon(u,v)\geq \Big(\frac{1}{2}-\frac{1}{2^*(s_{max})}\Big) \big(\|u\|^2+\|v\|^2\big)
\ee
and it follows that
\be\lab{2014-12-6-bue2}
c_\varepsilon\geq \Big(\frac{1}{2}-\frac{1}{2^*(s_{max})}\Big) \delta_\varepsilon^2.
\ee

\er

\bl\lab{2014-12-8-l1}
$c_\varepsilon$ is increasing by $\varepsilon$ in $[0,s_2)$.
\el
\bp
Let $(\phi,\psi)\neq (0,0)$ be fixed. By Lemma \ref{2014-12-5-l1}, for any $\varepsilon\in [0,s_2)$, there exists a  unique $t_\varepsilon>0$ such that $t_\varepsilon(\phi,\psi)\in \mathcal{N}_\varepsilon$. In fact,  $t_\varepsilon$ is implicitly  defined by the equation
\be\lab{2014-12-8-e1}
a(\phi,\psi)-b(\phi,\psi)t_{\varepsilon}^{2^*(s_1)-2}-2^*(s_2)\kappa c_\varepsilon(\phi,\psi) t_{\varepsilon}^{2^*(s_2)-2}=0.
\ee
It then follows that
$$\Phi_\varepsilon\big(t(\varepsilon) \phi,t(\varepsilon) \psi\big)$$
\be\lab{2014-12-8-e2}
=\Big[\frac{1}{2}-\frac{1}{2^*(s_2)}\Big]a(\phi,\psi)[t(\varepsilon)]^2
+\Big[\frac{1}{2^*(s_2)}-\frac{1}{2^*(s_1)}\Big]b(\phi,\psi)[t(\varepsilon)]^{2^*(s_1)}.
\ee

\vskip 0.02in
\noindent{\bf Case 1: $s_2>s_1$. } For this case, we see that $\frac{1}{2^*(s_2)}-\frac{1}{2^*(s_1)}>0$. Noting that $a(\phi,\psi)>0, b(\phi,\psi)>0$ and   Lemma \ref{2014-12-5-l1}, we obtain that
 \be\lab{2014-12-8-e3}
 \Phi_\varepsilon\big(t(\varepsilon) \phi,t(\varepsilon) \psi\big)\;\hbox{is increasing by $\varepsilon$ in $[0,s_2)$}.
 \ee
Hence,   we get that $c_\varepsilon$ is increasing by $\varepsilon$ in $[0,s_2)$.

\vskip 0.02in
\noindent{\bf Case 2: $s_2<s_1$.} By the Implicit Function Theorem, we see that $t(\varepsilon)\in C^1(\R)$ and $\frac{d}{d\varepsilon}t(\varepsilon)\geq 0$ by Lemma \ref{2014-12-5-l1}. Hence,
\begin{align}\lab{2014-12-8-e4}
&\frac{d}{d\varepsilon}\Phi_\varepsilon\big(t(\varepsilon) \phi,t(\varepsilon) \psi\big)\nonumber\\
=&2[\frac{1}{2}-\frac{1}{2^*(s_2)}]a(\phi,\psi)t(\varepsilon) t'(\varepsilon)+2^*(s_1)[\frac{1}{2^*(s_2)}-\frac{1}{2^*(s_1)}]b(\phi,\psi)[t(\varepsilon)]^{2^*(s_1)-1}t'(\varepsilon)\nonumber\\
=&\frac{t'(\varepsilon)}{t(\varepsilon)}\Big[[1-\frac{2}{2^*(s_2)}]a(\phi,\psi)[t(\varepsilon)]^{2}+[\frac{2^*(s_1)}{2^*(s_2)}-1]b(\phi,\psi)[t(\varepsilon)]^{2^*(s_1)}\Big]\nonumber\\
=&\frac{t'(\varepsilon)}{t(\varepsilon)}\Big[[1-\frac{2}{2^*(s_2)}]a(\phi,\psi)[t(\varepsilon)]^{2} +[\frac{2^*(s_1)}{2^*(s_2)}-1]\big[a(\phi,\psi)[t(\varepsilon)]^2-2^*(s_2)\kappa c_\varepsilon(\phi,\psi)[t(\varepsilon)]^{2^*(s_2)}\big]\Big]\nonumber\\
=&\frac{t'(\varepsilon)}{t(\varepsilon)}\Big[\frac{2^*(s_1)-2}{2^*(s_2)}a(\phi,\psi)[t(\varepsilon)]^2+[2^*(s_2)-2^*(s_1)]\kappa c_\varepsilon(\phi,\psi)[t(\varepsilon)]^{2^*(s_2)}\Big]\geq0.
\end{align}
Hence, we also obtain the conclusion of \eqref{2014-12-8-e3} for the case of $s_2<s_1$ and the proof is completed.
\ep


\subsection{Estimation on  the least energy of the approximation }
Recall  $U_\lambda:=\big(\frac{\mu_{s_1}(\Omega)}{\lambda}\big)^{\frac{1}{2^*(s_1)-2}}U$, where $U$ is a ground state solution of the following problem:
\be\lab{2014-12-5-e4}
-\Delta u=\mu_{s_1}(\Omega) \frac{u^{2^*(s_1)-1}}{|x|^{s_1}}\;\hbox{in}\;\Omega,\;
u=0\;\;\hbox{on}\;\partial \Omega.
\ee
Define the function
\be\lab{2014-11-23-we3}
\Psi_\lambda(u)=\frac{1}{2}\int_{\Omega}|\nabla u|^2 dx-\frac{\lambda}{2^*(s_1)}\int_{\Omega}\frac{|u|^{2^*(s_1)}}{|x|^{s_1}}dx.
\ee
Then
\be\lab{2014-11-23-we4}
m_\lambda=\Psi_\lambda(U_\lambda)=[\frac{1}{2}-\frac{1}{2^*(s_1)}]\big(\mu_{s_1}(\Omega)\big)^{\frac{2^*(s_1)}{2^*(s_1)-2}} \lambda^{-\frac{2}{2^*(s_1)-2}}
\ee
is the least energy.

\br\lab{2014-12-5-r2} Evidently,
for any $\varepsilon\in [0,s_2)$, we have that   $c_\varepsilon\leq  m_\lambda$ and $c_\varepsilon\leq m_\mu$. Hence,
\be\lab{2014-12-5-e5}
c_\varepsilon\leq [\frac{1}{2}-\frac{1}{2^*(s_1)}]\big(\mu_{s_1}(\Omega)\big)^{\frac{2^*(s_1)}{2^*(s_1)-2}} \big(\max\{\lambda,\mu\}\big)^{-\frac{2}{2^*(s_1)-2}}.
\ee
\er\hfill$\Box$

Define
\be\lab{2014-12-5-e6}
\eta_{1,\varepsilon}:=\inf_{v\in \Xi_\varepsilon} \|v\|^2
\ee
where
\be\lab{2014-12-5-e7}
\Xi_\varepsilon:=\{v\in D_{0}^{1,2}(\Omega):\;\int_\Omega a_\varepsilon(x)|U_\lambda|^{2^*(s_2)-2}|v|^2 dx=1\}.
\ee
Since $a_\varepsilon(x)$ is decreasing by $\varepsilon$, it is easy to see that $\eta_{1,\varepsilon}$ is increasing by $\varepsilon$.
\bl\lab{2014-12-5-l5}
Assume $s_1,s_2\in (0,2), \lambda, \mu\in (0,+\infty), \kappa>0,\alpha>1,\beta>1$ and $ \alpha+\beta= 2^*(s_2)$. Let $\varepsilon\in [0,s_2)$.
\begin{itemize}
\item[$(1)$] If $\beta<2$, then $c_\varepsilon<m_{\lambda}.$
\item[$(2)$] If $\beta>2$, then $(U_\lambda,0)$ is a local minimum point of $\Phi_\varepsilon$ in $\mathcal{N}_\varepsilon$.
\item[$(3)$] If $\beta=2$ and $\kappa>\frac{\eta_{1,\varepsilon}}{2^*(s_2)}$, then
     $c_\varepsilon < m_{\lambda}.$
\item[$(4)$] If $\beta=2$  and $0<\kappa< \frac{\eta_{1,\varepsilon}}{2^*(s_2)}$, then  $(U_\lambda,0)$ is a local minimum point of $\Phi_\varepsilon$ in $\mathcal{N}_\varepsilon$.
\end{itemize}
\el
\bp  The proofs are similar to those  in  Section 6.2.
\ep

\bl\lab{2014-12-6-l1}
$\eta_{1,\varepsilon}$ is continuous with  respect  to $\varepsilon\in [0,s_2)$.
\el
\bp
For any $\varepsilon_0\in [0,s_2)$, we shall prove that $\eta_{1,\varepsilon}$ is continuous at $\varepsilon=\varepsilon_0$.
Apply the argument of Lemma  \ref{2014-11-24-l3},  there exists some $0<v_0\in D_{0}^{1,2}(\Omega)$ such that
\be\lab{2014-12-6-e1}
\|v_0\|^2=\eta_{1,\varepsilon_0}\;\hbox{and}\;\int_\Omega a_{\varepsilon_0}(x)|U_\lambda|^{2^*(s_2)-2}v_0^2dx=1.
\ee
Take a sequence $\{\varepsilon_n\}\subset [0,s_2)$ with  $\varepsilon_n\downarrow \varepsilon_0$ as $n\rightarrow +\infty$.
Recall  that $\eta_{1,\varepsilon}$ is increasing by $\varepsilon$,    then  $\displaystyle\lim_{n\rightarrow +\infty}\eta_{1,\varepsilon_n}$ exists and satisfies
\be\lab{2014-12-6-e2}
\displaystyle\lim_{n\rightarrow +\infty}\eta_{1,\varepsilon_n}\geq \eta_{1,\varepsilon_0}.
\ee
On the other hand, since $a_{\varepsilon_n}(x)\rightarrow a_0(x) $  a.e.  in $\Omega$, recalling the decay property of $U_\lambda$ (see \cite[Theorem 1.2]{GhoussoubRobert.2006}, \cite[Lemma 2.1]{HsiaLinWadade.2010}, \cite[Lemma 2.6]{LinWadadeothers.2012}), it is easy to prove that
\be\lab{2014-12-6-e3}
\lim_{n\rightarrow \infty}\int_\Omega a_{\varepsilon_n}(x)|U_\lambda|^{2^*(s_2)-2}v_0^2dx=\int_\Omega a_{\varepsilon_0}(x)|U_\lambda|^{2^*(s_2)-2}v_0^2dx=1.
\ee
Hence,
\be\lab{2014-12-6-e3}
\lim_{n\rightarrow +\infty}\frac{\|v_0\|^2}{\int_\Omega a_{\varepsilon_n}(x)|U_\lambda|^{2^*(s_2)-2}v_0^2dx}=\eta_{1,\varepsilon_0}.
\ee
Then by the definition of $\eta_{1,\varepsilon}$, we see that
\be\lab{2014-12-6-e4}
\lim_{n\rightarrow +\infty}\eta_{1,\varepsilon_n}\leq \lim_{n\rightarrow +\infty}\frac{\|v_0\|^2}{\int_\Omega a_{\varepsilon_n}(x)|U_\lambda|^{2^*(s_2)-2}v_0^2dx}=\eta_{1,\varepsilon_0}.
\ee
By \eqref{2014-12-6-e2} and \eqref{2014-12-6-e4}, we obtain that $\eta_{1,\varepsilon}$ is right-continuous.

Secondly, we take a sequence $\{\varepsilon_n\}\subset [0,s_2)$ such that $\varepsilon_n\uparrow \varepsilon_0$ as $n\rightarrow +\infty$. By Lemma Lemma  \ref{2014-11-24-l3} again, we may assume that $\{v_n\}\subset D_{0}^{1,2}(\Omega)$ such that
\be\lab{2014-12-6-e5}
\|v_n\|^2=\eta_{1,\varepsilon_n}\;\hbox{and}\; \int_\Omega a_{\varepsilon_n}(x)|U_\lambda|^{2^*(s_2)-2}v_n^2dx\equiv 1.
\ee
Up to a subsequence, we may assume that $v_n\rightharpoonup v_0$ in $D_{0}^{1,2}(\Omega)$ and $v_n\rightarrow v_0$ a.e. in $\Omega$. Similarly, we can prove that
\be\lab{2014-12-6-e6}
\int_\Omega a_{\varepsilon_0}(x)|U_\lambda|^{2^*(s_2)-2}v_0^2dx=\lim_{n\rightarrow +\infty} \int_\Omega a_{\varepsilon_n}(x)|U_\lambda|^{2^*(s_2)-2}v_n^2dx=1.
\ee
It follows that
\be\lab{2014-12-6-e7}
\|v_0\|^2\leq \liminf_{n\rightarrow +\infty}\|v_n\|^2=\lim_{n\rightarrow +\infty}\eta_{1,\varepsilon_n}.
\ee
Therefore,
\be\lab{2014-12-6-e8}
\eta_{1,\varepsilon_0}\leq \frac{\|v_0\|^2}{\int_\Omega a_{\varepsilon_0}(x)|U_\lambda|^{2^*(s_2)-2}v_0^2dx}\leq \lim_{n\rightarrow +\infty}\eta_{1,\varepsilon_n}.
\ee
On the other hand, by the monotonicity, we can obtain that reverse inequality. Hence,
$\displaystyle \eta_{1,\varepsilon_0}\leq \frac{\|v_0\|^2}{\int_\Omega a_{\varepsilon_0}(x)|U_\lambda|^{2^*(s_2)-2}v_0^2dx}= \lim_{n\rightarrow +\infty}\eta_{1,\varepsilon_n},$
i.e., $\eta_{1,\varepsilon}$ is left-continuous. The proof is completed.
\ep

Similarly, we define
\be\lab{2014-12-6-e9}
\eta_{2,\varepsilon}:=\inf_{u\in \Theta_\varepsilon} \|u\|^2
\ee
where
\be\lab{2014-12-6-e10}
\Theta_\varepsilon:=\Big\{u\in D_{0}^{1,2}(\Omega):\;\int_\Omega a_\varepsilon(x)|U_\mu|^{2^*(s_2)-2}|u|^2 dx=1\Big\}.
\ee
We also have that $\eta_{2,\varepsilon}$ is increasing by $\varepsilon\in [0,s_2)$ and continuous. Furthermore, we can propose the following results without proof.

\bl\lab{2014-12-6-l1}
Assume $s_1,s_2\in (0,2), \lambda, \mu\in (0,+\infty), \kappa>0,\alpha>1,\beta>1$ and $ \alpha+\beta= 2^*(s_2)$. Let $\varepsilon\in [0,s_2)$.
\begin{itemize}
\item[$(1)$] If $\alpha<2$, then $c_\varepsilon<m_{\mu}.$
\item[$(2)$] If $\alpha>2$, then $(0,U_\mu)$ is a local minimum point of $\Phi_\varepsilon$ in $\mathcal{N}_\varepsilon$.
\item[$(3)$] If $\alpha=2,\kappa>\frac{\eta_{2,\varepsilon}}{2^*(s_2)}$, then
     $c_\varepsilon <m_{\mu}.$
\item[$(4)$] If $\alpha=2,0<\kappa< \frac{\eta_{2,\varepsilon}}{2^*(s_2)}$, then  $(0,U_\mu)$ is a local minimum point of $\Phi_\varepsilon$ in $\mathcal{N}_\varepsilon$.
\end{itemize}
\el
Now we  can  obtain the following estimation on $c_\varepsilon$:

\bl\lab{2014-12-6-l2}
Assume $s_1,s_2\in (0,2), \lambda, \mu\in (0,+\infty), \kappa>0,\alpha>1,\beta>1$ and $ \alpha+\beta= 2^*(s_2)$. Let $\varepsilon\in [0,s_2)$, then we have
$$c_\varepsilon<\min\{m_\lambda,m_\mu\}= \left[\frac{1}{2}-\frac{1}{2^*(s_1)}\right]\Big(\mu_{s_1}(\Omega)\Big)^{\frac{2^*(s_1)}{2^*(s_1)-2}} \Big(\max\{\lambda,\mu\}\Big)^{-\frac{2}{2^*(s_1)-2}}$$
if one of the following holds:
\begin{itemize}
\item[$(a)$]$\lambda>\mu $ and either $ 1<\beta<2$ or $\begin{cases}\beta=2\\ \kappa>\frac{\eta_{1,\varepsilon}}{2^*(s_2)} \end{cases}$;
\item[$(b)$] $\lambda=\mu  $ and either $   \min\{\alpha,\beta\}<2$ or $\begin{cases} \min\{\alpha,\beta\}=2,\\ \kappa>\frac{\eta_{1,\varepsilon}}{2^*(s_2)}= \frac{\eta_{2,\varepsilon}}{2^*(s_2)}\end{cases}$;
\item[$(c)$] $\lambda<\mu  $ and either $ 1<\alpha<2$ or $\begin{cases}\alpha=2\\ \kappa>\frac{\eta_{2,\varepsilon}}{2^*(s_2)} \end{cases}$.
\end{itemize}
\el
\bp
It is a direct conclusion following by Lemma \ref{2014-12-5-l5} and Lemma \ref{2014-12-6-l1}.
\ep

\subsection{Positive ground state to the approximation  problem  \eqref{zz-620-1}}
In this subsection, we assume that $\varepsilon\in (0,s_2)$ is fixed. Then we can obtain the following result.
\bt\lab{2014-12-6-th1}
Assume $s_1,s_2\in (0,2), \lambda, \mu\in (0,+\infty), \kappa>0,\alpha>1,\beta>1$ and $ \alpha+\beta= 2^*(s_2)$. Then problem $\eqref{zz-620-1}$ possesses a positive ground state solution $(\phi_\varepsilon,\psi_\varepsilon)$ provided further one of the following conditions holds:
\begin{itemize}
\item[$(1)$]$\lambda>\mu$ and either $ 1<\beta<2$ or $\begin{cases}\beta=2\\ \kappa>\frac{\eta_{1,\varepsilon}}{2^*(s)} \end{cases}$;
\item[$(2)$] $\lambda=\mu $ and either $ \min\{\alpha,\beta\}<2$ or $\begin{cases} \min\{\alpha,\beta\}=2,\\ \kappa>\frac{\eta_{1,\varepsilon}}{2^*(s)}= \frac{\eta_{2,\varepsilon}}{2^*(s)}\end{cases}$;
\item[$(3)$] $\lambda<\mu $ and either $ 1<\alpha<2$ or $\begin{cases}\alpha=2\\ \kappa>\frac{\eta_{2,\varepsilon}}{2^*(s)} \end{cases}$.
\end{itemize}
\et

\bo\lab{2014-12-6-prop1}
Assume that $\varepsilon\in (0,s_2)$ and $\{(u_n,v_n)\}$ is a bounded $(PS)_c$ sequence of $\Phi_\varepsilon$.
Up to a subsequence, we assume that $(u_n,v_n)\rightharpoonup (\phi,\psi)$  weakly in $\mathscr{D}$. Set $\tilde{u}_n:=u_n-\phi, \tilde{v}_n:=v_n-\psi$, then we have that
\be\lab{2014-12-6-e12}
\Psi'_\lambda(\tilde{u}_n)\rightarrow 0\;\hbox{and}\;\Psi'_\mu(\tilde{v}_n)\rightarrow 0 \;\hbox{in}\;H^{-1}(\Omega),
\ee
where   $\Psi_\lambda$ is defined in \eqref{2014-11-23-we3}.
\eo
\bp
Under the assumptions, we see that
\be\lab{2014-12-6-e13}
\big\langle\Phi'_\varepsilon(u_n,v_n),(h,0)\big\rangle=o(1)\|h\|
\ee
Since $(u_n,v_n)\rightharpoonup (\phi,\psi)$, it is easy to see that $\Phi'_\varepsilon(\phi,\psi)=0$. Then we have
\be\lab{2014-12-6-e14}
\big\langle \Phi'_\varepsilon(\phi,\psi), (h,0)\big\rangle=0.
\ee
By Lemma \ref{2014-11-26-xl3}  and H\"older inequality, it is easy to see that
\be\lab{2014-12-6-e15}
\int_\Omega a_\varepsilon(x)|u_n|^{\alpha-2}u_n|v_n|^\beta h dx-\int_\Omega a_\varepsilon(x)|\phi|^{\alpha-2}\phi|\psi|^\beta h dx=o(1)\|h\|.
\ee
It follows from \eqref{2014-12-6-e13},\eqref{2014-12-6-e14} and \eqref{2014-12-6-e15} that
\be\lab{2014-12-6-e16}
\int_\Omega \nabla (u_n-\phi)\nabla h dx-\lambda\int_\Omega \Big(\frac{|u_n|^{2^*(s_1)-2}u_n}{|x|^{s_1}}-\frac{|\phi|^{2^*(s_1)-2}\phi}{|x|^{s_1}}\Big)hdx=o(1)\|h\|.
\ee
By \cite[Lemma 3.3]{GhoussoubKang.2004} or \cite[Lemma 3.2]{CeramiZhongZou.2014}, we see that
\be\lab{2014-12-6-e17}
\frac{|u_n|^{2^*(s_1)-2}u_n}{|x|^{s_1}}-\frac{|u_n-\phi|^{2^*(s_1)-2}(u_n-\phi)}{|x|^{s_1}}\rightarrow \frac{|\phi|^{2^*(s_1)-2}\phi}{|x|^{s_1}}\;\hbox{in}\;H^{-1}(\Omega).
\ee
Hence, by \eqref{2014-12-6-e16} and \eqref{2014-12-6-e17}, we obtain that
\be\lab{2014-12-6-e18}
\Psi'_\lambda(\tilde{u}_n)\rightarrow 0\;\hbox{in}\;H^{-1}(\Omega).
\ee
Apply the similar arguments, we can prove that
$
\Psi'_\mu(\tilde{v}_n)\rightarrow 0\;\hbox{in}\;H^{-1}(\Omega).
$
\ep

\bc\lab{2014-12-6-cro1}
Under the assumptions of Proposition \ref{2014-12-6-prop1} and furthermore we assume that
$$c<\min\{m_\lambda,m_\mu\}= \left[\frac{1}{2}-\frac{1}{2^*(s_1)}\right]\big(\mu_{s_1}(\Omega)\big)^{\frac{2^*(s_1)}{2^*(s_1)-2}} \big(\max\{\lambda,\mu\}\big)^{-\frac{2}{2^*(s_1)-2}}.$$
Then up to a subsequence, $(u_n,v_n)\rightarrow (\phi,\psi)$ strongly in $\mathscr{D}$ and $(\phi,\psi)$ satisfies
$\displaystyle \Phi_\varepsilon(\phi,\psi)=c\;\hbox{and}\;\Phi'_\varepsilon(\phi,\psi)=0\;\hbox{in}\;\mathscr{D}^*.$
\ec
\bp
We prove it by way of negation. Assume that $(u_n,v_n)\not\rightarrow (\phi,\psi)$, then at least one of the following holds:\\
(i)$u_n\not\rightarrow \phi$\;in $D_{0}^{1,2}(\Omega)$,\quad
(ii)$v_n\not\rightarrow \psi$\;in $D_{0}^{1,2}(\Omega)$.\\
Without loss of generality, we assume $(i)$.
By Proposition \ref{2014-12-6-prop1}, we see that $\Psi'_\lambda(\tilde{u}_n)\rightarrow 0$ in $H^{-1}(\Omega)$. Since $\tilde{u}_n=u_n-\phi\not\rightarrow 0$ in $D_{0}^{1,2}(\Omega)$, it is easy to see that
\be\lab{2014-12-6-e20}
\liminf_{n\rightarrow +\infty}\Psi_\lambda(\tilde{u}_n)\geq m_\lambda.
\ee
On the other hand, by the Br\'ezis-Lieb type lemma (see \cite[Lemma 3.3]{GhoussoubKang.2004}), we have
\be\lab{2014-12-6-e21}
\Phi_\varepsilon(u_n,v_n)=\Phi_\varepsilon(\tilde{u}_n,\tilde{v}_n)+\Phi_\varepsilon(\phi,\psi)+o(1).
\ee
By Lemma \ref{2014-11-26-xl3} again, we see that
\be\lab{2014-12-6-e22}
\Phi_\varepsilon(\tilde{u}_n,\tilde{v}_n)=\Psi_\lambda(\tilde{u}_n)+\Psi_\mu(\tilde{v}_n)+o(1).
\ee
Since $\Psi'_\mu(\tilde{v}_n)\rightarrow 0$ in $H^{-1}(\Omega)$, it is easy to prove that $\displaystyle \liminf_{n\rightarrow +\infty}\Psi_\mu(\tilde{v}_n)\geq 0$.
We also note that $\Phi_\varepsilon(\phi,\psi)\geq 0$. Then by \eqref{2014-12-6-e21}, \eqref{2014-12-6-e22} and \eqref{2014-12-6-e20}, we have
\be\lab{2014-12-6-e23}
c=\lim_{n\rightarrow +\infty}\Phi_\varepsilon(u_n,v_n)\geq \lim_{n\rightarrow +\infty}\Psi_\lambda(\tilde{u}_n)\geq m_\lambda,
\ee
a contradiction.
\ep

\noindent{\bf Proof of Theorem \ref{2014-12-6-th1}:}
Let $\{(u_n,v_n)\}\subset \mathcal{N}_\varepsilon$ be a minimizing sequence. Then it is easy to see that
$\displaystyle \Phi_\varepsilon(u_n,v_n)\rightarrow c_\varepsilon\;\hbox{and}\;\Phi'_\varepsilon\big|_{\mathcal{N}_\varepsilon}(u_n,v_n)\rightarrow 0\;\hbox{in}\;\mathscr{D}^*.$
It is standard to prove  that $(u_n,v_n)$ is bounded in $\mathscr{D}$ and is also a $(PS)_{c_\varepsilon}$ sequence of $\Phi_\varepsilon$. By Lemma \ref{2014-12-6-l2}, we have
\be\lab{2014-12-6-e24}
c_\varepsilon<\min\{m_\lambda,m_\mu\}= [\frac{1}{2}-\frac{1}{2^*(s_1)}]\big(\mu_{s_1}(\Omega)\big)^{\frac{2^*(s_1)}{2^*(s_1)-2}} \big(\max\{\lambda,\mu\}\big)^{-\frac{2}{2^*(s_1)-2}}.
\ee
Hence, by Corollary \ref{2014-12-6-cro1}, there exists some $(\phi,\psi)\in \mathscr{D}$ and up to a subsequence, $(u_n,v_n)\rightarrow (\phi,\psi)$ strongly in $\mathscr{D}$. Moreover, we have
$
\Phi_\varepsilon(\phi,\psi)=c_\varepsilon\;\hbox{and}\;\Phi'_\varepsilon(\phi,\psi)=0.
$
Thus, $(\phi,\psi)$ is a minimizer of $c_\varepsilon$. It is easy to see that $(|\phi|, |\psi|)$ is also a minimizer. Hence, without loss of generality, we may assume that $\phi\geq 0,\psi\geq 0$  and it follows that $(\phi,\psi)$ is a nonnegative solution of \eqref{zz-620-1}.
Recalling \eqref{2014-12-6-e24}, it is easy to see that $\phi\neq 0, \psi\neq 0$.
 Finally, by the strong maximum principle, we can obtain that $\phi>0,\psi>0$. That is, we obtain that $(\phi,\psi)$ is a positive ground state solution of  \eqref{zz-620-1}.\hfill$\Box$

\subsection{Geometric structure of positive ground state to \eqref{zz-620-1}}
Now, let us define the mountain pass  value
\be\lab{2014-3-10-e2}
\tilde{c}_\varepsilon:=\inf_{\gamma\in \Gamma_\varepsilon}\max_{t\in [0,1]}\Phi_{\varepsilon}(\gamma(t)),
\ee
 where $\Gamma_\varepsilon:=\{\gamma(t)\in C([0,1],  \mathscr{D}):    \gamma(0)=(0,0), \Phi_{\varepsilon}(\gamma(1))<0\}$.
 We have the following result.
 \bt\lab{2014-12-6-th2}
 Assume $s_1,s_2\in (0,2), \lambda, \mu\in (0,+\infty), \kappa>0,\alpha>1,\beta>1$ and $ \alpha+\beta= 2^*(s_2)$. Let $\varepsilon\in (0,s_2)$  and one of the following hold:
  \begin{itemize}
\item[$(i)$]$\lambda>\mu$ and either $ 1<\beta<2$ or $\begin{cases}\beta=2\\ \kappa>\eta_{1,\varepsilon} \end{cases}$;
\item[$(ii)$] $\lambda=\mu$ and either $  \min\{\alpha,\beta\}<2$ or $\begin{cases} \min\{\alpha,\beta\}=2,\\ \kappa>\eta_{1,\varepsilon}= \eta_{2,\varepsilon}\end{cases}$;
\item[$(iii)$] $\lambda<\mu$ and either $  1<\alpha<2$ or $\begin{cases}\alpha=2\\ \kappa>\eta_{2,\varepsilon} \end{cases}$.
\end{itemize}
  Then $c_\varepsilon=\tilde{c}_\varepsilon$ and any positive ground state solution of system   \eqref{zz-620-1}  is a mountain pass solution.
 \et
\bp
It is easy to check that $\Phi_\varepsilon$ satisfies the mountain pass geometric structure.
Recalling the existence result of Theorem \ref{2014-12-6-th1}, let $(\phi,\psi)$ be a positive ground state solution of  \eqref{zz-620-1}. Define $\gamma_0(t):=tT(\phi,\psi)$ for some $T>0$ large enough such that $\Phi_{\varepsilon}(T\phi,T\psi)<0$. Then it is easy to see that $\gamma_0\in \Gamma_\varepsilon$.
By Lemma \ref{2014-12-5-l1}, we have
\be\lab{2014-12-6-e26}
\Phi_\varepsilon(\phi,\psi)=\max_{t>0}\Phi_\varepsilon(t\phi,t\psi).
\ee
Hence,
\be\lab{2014-12-6-e27}
\tilde{c}_\varepsilon\leq \max_{t\in [0,1]}\Phi_\varepsilon(\gamma_0(t))=\Phi_\varepsilon(\phi,\psi)=c_\varepsilon.
\ee
Under the assumptions, it is standard to prove that $\tilde{c}_\varepsilon$ is also a critical value and there exists a solution $(\tilde{\phi},\tilde{\psi})$ such that
$
\Phi_\varepsilon(\tilde{\phi},\tilde{\psi})=\tilde{c}_\varepsilon\;\hbox{and}\;\Phi'_\varepsilon(\tilde{\phi},\tilde{\psi})=0\;\hbox{in}\;\mathscr{D}^*.
$
Then we see that
$
(\tilde{\phi},\tilde{\psi})\in \mathcal{N}_\varepsilon.
$
Hence, by the definition of $c_\varepsilon$ we see that
\be\lab{2014-12-6-e30}
c_\varepsilon:=\inf_{(u,v)\in \mathcal{N}_\varepsilon}\Phi_\varepsilon(u,v)\leq \Phi_\varepsilon(\tilde{\phi},\tilde{\psi})=\tilde{c}_\varepsilon.
\ee
By \eqref{2014-12-6-e27} and \eqref{2014-12-6-e30}, we obtain that
$
c_\varepsilon=\tilde{c}_\varepsilon.
$
For any positive ground state solution, by the arguments as above, we have the  mountain path $\gamma_0\in \Gamma_\varepsilon$ and  thus, the  positive ground state   is indeed a mountain pass solution.
\ep

\br\lab{2014-12-6-xr1}\quad
\begin{itemize}
\item[$(i)$] Recalling that for $\varepsilon\in [0,s_2)$, both $\eta_{1,\varepsilon}$ and $\eta_{2,\varepsilon}$ are increasing by $\varepsilon$ and continuous with  respect  to $\varepsilon$.  When $\kappa>\frac{\eta_{i,0}}{2^*(s_2)}, i\in\{1,2\}$, then by the continuity, we see that  $\kappa>\frac{\eta_{i,\varepsilon}}{2^*(s_2)}$ when $\varepsilon$ is small enough.
\item[$(ii)$] Note that the proof of $\tilde{c}_\varepsilon=c_\varepsilon$ for $\varepsilon\in(0,s_2)$ depends heavily on the existence of  the ground state solution. When $\varepsilon=0$, the existence of ground state solution is still unknown. However, we will prove that the result $\tilde{c}_0=c_0$ is also satisfied (see Corollary \ref{2014-12-6-zcro1} below).
\end{itemize}
\er

\bl\lab{2014-12-6-zl1}
$\tilde{c}_{\varepsilon}\geq \tilde{c}_0$ and $\displaystyle\lim_{\varepsilon\rightarrow 0^+} \tilde{c}_{\varepsilon}=\tilde{c}_0$
\el
\bp
By the monotonicity of $a_\varepsilon(x)$, it is easy to see that $\tilde{c}_{\varepsilon}\geq \tilde{c}_0$. Hence,
\be\lab{2014-12-6-zhe2}
\lim_{\varepsilon\rightarrow 0^+} \tilde{c}_{\varepsilon}\geq\tilde{c}_0.
\ee
Next, we only need to prove the inverse inequality.
For any $\delta>0$, there exists $\gamma_0\in \Gamma_0$ such that
\be\lab{2014-12-6-zhe3}
\max_{t\in [0,1]}\Phi_0(\gamma_0(t))<\tilde{c}_0+\delta.
\ee
Denote $\gamma_0(1):=(\phi,\psi)$, since $\gamma_0\in \Gamma_0$, we have
$\Phi_0(\phi,\psi)<0$.

\vskip 0.02in
\noindent{\bf Case 1:} If $|\phi|^\alpha|\psi|^\beta\equiv 0$, it is easy to see that $\Phi_\varepsilon(\phi,\psi)=\Phi_0(\phi,\psi)<0$. Hence, $\gamma_0\in \Gamma_\varepsilon$ for all $\varepsilon\in [0,s_2)$ for this case.

\vskip 0.02in
\noindent{\bf Case 2:}If $|\phi|^\alpha|\psi|^\beta\not\equiv 0$,
then by the Lebesgue's dominated convergence theorem, we have
\be\lab{2014-12-6-zhe4}
\lim_{\varepsilon\rightarrow 0^+}\int_\Omega a_\varepsilon(x)|\phi|^\alpha |\psi|^\beta dx=\int_\Omega a_0(x)|\phi|^\alpha |\psi|^\beta dx.
\ee
Hence, we have $\Phi_\varepsilon(\phi,\psi)<0$ when $\varepsilon$ is small enough. Thus, we also obtain that $\gamma_0\in \Gamma_\varepsilon$ when $\varepsilon$ is small enough. Now, we take an  arbitrary sequence $\varepsilon_n\downarrow 0$ as $n\rightarrow +\infty$.   Choose  $t_n\in [0,1]$ such that
\be\lab{2014-12-6-zhe5}
\Phi_{\varepsilon_n}(\gamma_0(t_n))=\max_{t\in[0,1]} \Phi_{\varepsilon_n}(\gamma_0(t)).
\ee
Up to a subsequence, we assume that $t_n\rightarrow t^*\in [0,1]$ and denote that
\be\lab{2014-12-6-zhe6}
\gamma_0(t_n):=(u_n,v_n), \gamma_0(t^*):=(u^*,v^*).
\ee
Since $\gamma_0\in C([0,1],  \mathscr{D})$, we obtain that $(u_n,v_n)\rightarrow (u^*,v^*)$ and it follows that
\be\lab{2014-12-6-zhe7}
\Phi_{\varepsilon_n}(u_n,v_n)=\Phi_{\varepsilon_n}(u^*,v^*)+o(1).
\ee
By the Lebesgue's dominated convergence theorem again, we have
\be\lab{2014-12-6-zhe8}
\Phi_{\varepsilon_n}(u^*,v^*)=\Phi_0(u^*,v^*)+o(1).
\ee
Hence, by \eqref{2014-12-6-zhe7} and \eqref{2014-12-6-zhe8}, we have
$
\Phi_{\varepsilon_n}(u_n,v_n)=\Phi_0(u^*,v^*)+o(1).
$
Then
\begin{align}\lab{2014-12-6-zhe10}
\tilde{c}_{\varepsilon_n}\leq &\Phi_{\varepsilon_n}(\gamma_0(t_n))=\Phi_{\varepsilon_n}(u_n,v_n)
=\Phi_0(u^*,v^*)+o(1)=\Phi_0(\gamma_0(t^*))+o(1)\nonumber\\
\leq&\max_{t\in [0,1]}\Phi_0(\gamma_0(t))+o(1)=\Phi_0(\phi,\psi)+o(1)
\leq\tilde{c}_0+\delta+o(1).
\end{align}
Let $n\rightarrow +\infty$, we obtain that
$\displaystyle \lim_{n\rightarrow +\infty}\tilde{c}_{\varepsilon_n}\leq \tilde{c}_0+\delta$. Hence,
$
\lim_{\varepsilon\rightarrow 0^+} \tilde{c}_{\varepsilon}\leq \tilde{c}_0.
$
Insert \eqref{2014-12-6-zhe2},  we complete the proof.
\ep
\bc\lab{2014-12-6-zcro1}
$c_0=\tilde{c}_0$ and $\displaystyle \lim_{\varepsilon\rightarrow 0^+} c_{\varepsilon}=c_0$.
\ec
\bp
For any $(u,v)\neq (0,0)$, define $\gamma(t)=t(u,v)$, then we see that $\gamma\in \Gamma_0$. Hence, it is easy to see that
$
\tilde{c}_0\leq c_0.
$
On the other hand, by Theorem \ref{2014-12-6-th2} and Lemma \ref{2014-12-6-zl1}, we have
$
\tilde{c}_0=\lim_{\varepsilon\rightarrow 0^+}\tilde{c}_\varepsilon=\lim_{\varepsilon\rightarrow 0^+}c_\varepsilon.
$
By Lemma \ref{2014-12-8-l1}, we have
$
\lim_{\varepsilon\rightarrow 0^+}c_\varepsilon\geq c_0.
$
Hence, we obtain that $\tilde{c}_0=c_0$ and $\displaystyle \lim_{\varepsilon\rightarrow 0^+} c_{\varepsilon}=c_0$.
\ep
\subsection{The existence of the positive ground state to the  original system}

Take $\{\varepsilon_n\}\subset (0,s_2)$ such that $\varepsilon_n\downarrow 0$ as $n\rightarrow +\infty$.   By Theorem \ref{2014-12-6-th1},  the system   \eqref{zz-620-1}  possesses a positive ground state solution  $(u_n,v_n)$. By Remark \ref{2014-12-6-bur1}, we have
\be\lab{2014-12-6-e41}
c_{\varepsilon_n}=\Phi_{\varepsilon_n}(u,v)\geq \Big(\frac{1}{2}-\frac{1}{2^*(s_{max})}\Big) \big(\|u_n\|^2+\|v_n\|^2\big).
\ee
By Corollary \ref{2014-12-6-zcro1}, we have $c_{\varepsilon_n}\rightarrow c_0$. Hence,  $\{(u_n,v_n)\}$
is bounded in $\mathscr{D}$. Up to a subsequence, we may assume that $(u_n,v_n)\rightharpoonup (u_0,v_0)$ weakly in $\mathscr{D}$ and $u_n\rightarrow u_0,v_n\rightarrow v_0$ a.e. in $\Omega$.
We shall establish the following results which are useful to prove our main theorem.

\bl\lab{2014-12-6-wl1}
$(u_0,v_0)$ satisfies $\Phi'_0(u_0,v_0)=0$ in $\mathscr{D}^*$.
\el
\bp
We claim that for any $\phi\in D_{0}^{1,2}(\Omega)$, we have
\be\lab{2014-12-6-e42}
\lim_{n\rightarrow+\infty}\int_\Omega a_{\varepsilon_n}(x)|u_n|^{\alpha-2}u_n|v_n|^\beta \phi dx=
\int_\Omega a_{0}(x)|u_0|^{\alpha-2}u_0|v_0|^\beta \phi dx.
\ee
Without loss of generality, we may also assume that $\phi\geq 0$. If not, we view $\phi=\phi_+-\phi_-$, and we discuss  $\phi_+$ and $\phi_-$ respectively.
Firstly, by Fatou's Lemma, we have
\be\lab{2014-12-6-el44}
\int_\Omega a_0(x)|u_0|^{\alpha-2}u_0|v_0|^\beta \phi dx\leq \liminf_{n\rightarrow +\infty} \int_\Omega a_{\varepsilon_n}(x)|u_n|^{\alpha-2}u_n|v_n|^\beta \phi dx.
\ee
On the other hand, since $a_{\varepsilon_n}(x)\leq a_0(x)$, we have
\be\lab{2014-12-6-el45}
\int_\Omega a_{\varepsilon_n}(x)|u_n|^{\alpha-2}u_n|v_n|^\beta \phi dx\leq \int_\Omega a_{0}(x)|u_n|^{\alpha-2}u_n|v_n|^\beta \phi dx.
\ee
Further,  since $(u_n,v_n)\rightharpoonup (u_0,v_0)$ in $\mathscr{D}$, it is easy to see that $$\displaystyle |u_n|^{\alpha-2}u_n|v_n|^\beta \rightharpoonup |u_0|^{\alpha-2}u_0|v_0|^\beta \quad \hbox{in } \;\; L^{\frac{2^*(s_2)}{2^*(s_2)-1}}(\Omega, a_0(x)dx),$$ then we have
\be\lab{2014-12-6-el46}
\lim_{n\rightarrow +\infty}\int_\Omega a_{0}(x)|u_n|^{\alpha-2}u_n|v_n|^\beta \phi dx=\int_\Omega a_{0}(x)|u_0|^{\alpha-2}u_0|v_0|^\beta \phi dx.
\ee
By \eqref{2014-12-6-el45} and \eqref{2014-12-6-el46}, we have
\be\lab{2014-12-6-el47}
\limsup_{n\rightarrow +\infty} \int_\Omega a_{\varepsilon_n}(x)|u_n|^{\alpha-2}u_n|v_n|^\beta \phi dx\leq \int_\Omega a_{0}(x)|u_0|^{\alpha-2}u_0|v_0|^\beta \phi dx.
\ee
Hence, from  \eqref{2014-12-6-el44} and \eqref{2014-12-6-el47}, we prove \eqref{2014-12-6-e42}.
Similarly, we can prove that for any $\psi\in D_{0}^{1,2}(\Omega)$, we have
\be\lab{2014-12-6-e43}
\lim_{n\rightarrow+\infty}\int_\Omega a_{\varepsilon_n}(x)|u_n|^{\alpha}|v_n|^{\beta-2}v_n \psi dx=
\int_\Omega a_{0}(x)|u_0|^{\alpha}|v_0|^{\beta-2}v_0 \psi dx.
\ee
Recalling that $(u_n,v_n)$ are critical point of $\Phi_{\varepsilon_n}$, for any $(\phi,\psi)\in \mathscr{D}$, we have
\be\lab{2014-12-6-e44}
\big\langle\Phi'_{\varepsilon_n}(u_n,v_n), (\phi,\psi)\big\rangle\equiv 0.
\ee
Then by \eqref{2014-12-6-e42}, \eqref{2014-12-6-e43} and $(u_n,v_n)\rightharpoonup (u_0,v_0)$ weakly  in $\mathscr{D}$, we obtain that
\be\lab{2014-12-6-be45}
\big\langle\Phi'_{0}(u_0,v_0), (\phi,\psi)\big\rangle=0.
\ee
Hence,  $\Phi'_{0}(u_0,v_0)=0$ in $\mathscr{D}^*$.
\ep

\vskip0.12in

\bl\lab{2014-12-6-wl2}
 If $(u_0,v_0)\neq (0,0)$, then $\Phi_0(u_0,v_0)=c_0>0$.
\el
\bp
Since $(u_n,v_n)$ is a positive ground state solution of $(P_{\varepsilon_n})$, it is easy to prove that
\be\lab{2014-12-6-be46}
c_{\varepsilon_n}=\Phi_{\varepsilon_n}(u_n,v_n)=\Big[\frac{1}{2}-\frac{1}{2^*(s_1)}\Big]b(u_n,v_n)+\Big[\frac{2^*(s_2)}{2}-1\Big]\kappa c_{\varepsilon_n}(u_n,v_n).
\ee
By Lemma \ref{2014-12-6-wl1}, we also have
\be\lab{2014-12-6-e47}
\Phi_0(u_0,v_0)=\Big[\frac{1}{2}-\frac{1}{2^*(s_1)}\Big]b(u_0,v_0)+\Big[\frac{2^*(s_2)}{2}-1\Big]\kappa c_{0}(u_0,v_0).
\ee
Noting that by Fatou's Lemma, we have
\be\lab{2014-12-6-e48}
b(u_0,v_0)\leq \liminf_{n\rightarrow +\infty}b(u_n,v_n)
\ee
and
\be\lab{2014-12-6-e49}
c_0(u_0,v_0)\leq \liminf_{n\rightarrow +\infty}c_{\varepsilon_n}(u_n,v_n).
\ee
Hence, $\displaystyle\Phi_0(u_0,v_0)\leq \lim_{n\rightarrow +\infty}c_{\varepsilon_n}$ and $\Phi_0(u_0,v_0)\geq 0$ follows by \eqref{2014-12-6-e47}.
By Corollary \ref{2014-12-6-zcro1}, we have $\displaystyle \lim_{n\rightarrow +\infty}c_{\varepsilon_n}=c_0$, hence
\be\lab{2014-12-8-ge1}
\Phi_0(u_0,v_0)\leq c_0.
\ee
On the other hand,
when $(u_0,v_0)\neq (0,0)$, it is trivial that $\Phi_0(u_0,v_0)\geq c_0>0$ by Lemma \ref{2014-12-6-wl1} and the definition of $c_0$. Hence, we obtain that $\Phi_0(u_0,v_0)=c_0$ if $(u_0,v_0)\neq (0,0)$.
\ep

\bl\lab{2014-12-6-wl3}
\begin{align*}\Phi_0(u_0,v_0)<\min\{m_\lambda,m_\mu\}
= \left[\frac{1}{2}-\frac{1}{2^*(s_1)}\right]\Big(\mu_{s_1}(\Omega)\Big)^{\frac{2^*(s_1)}{2^*(s_1)-2}} \Big(\max\{\lambda,\mu\}\Big)^{-\frac{2}{2^*(s_1)-2}}.
\end{align*}
\el
\bp
It is a direct conclusion by Lemma \ref{2014-12-6-l2} and Lemma \ref{2014-12-6-wl2}.
\ep

\bc\lab{2014-12-6-wcro2}
If $(u_0,v_0)\neq (0,0)$, then $(u_0,v_0)$ is a positive ground solution of \eqref{P}.
\ec
\bp
Since $(u_n,v_n)$ are positive and $u_n\rightarrow u_0,v_n\rightarrow v_0$ a.e. in $\Omega$. We see that $u_0\geq 0, v_0\geq 0$. If $v_0=0$, then we see that $\Psi'_\lambda(u_0)=0$ and $u_0\neq 0$. Hence, $\Phi_0(u_0,v_0)=\Psi_\lambda(u_0)\geq m_\lambda$, a contradiction to Lemma \ref{2014-12-6-wl3}. Similarly, if $u_0=0,v_0\neq 0$, we see that $\Phi_0(u_0,v_0)\geq m_\mu$, also a contradiction to Lemma \ref{2014-12-6-wl3}. Hence, $u_0\neq 0, v_0\neq 0$ and $\Phi_0(u_0,v_0)=c_0$ by Lemma \ref{2014-12-6-wl2}. That is, $(u_0,v_0)$ is a nontrivial and nonnegative ground state solution of \eqref{P}. Finally, by the strong maximum principle, we can prove that $(u_0,v_0)$ is a positive solution.
\ep

\bl\lab{2014-12-6-wl4}
Assume that $\displaystyle\liminf_{n\rightarrow +\infty}\int_\Omega a_{\varepsilon_n}|u_n|^\alpha|v_n|^\beta=0$, then $\Psi'_\lambda(u_n)\rightarrow 0, \Psi'_\mu(v_n)\rightarrow 0$ in $H^{-1}(\Omega)$.
\el
\bp
Under the assumptions, we claim that up to a subsequence,
\be\lab{2014-12-6-e53}
\int_\Omega \big|a_{\varepsilon_n}(x)|u_n|^{\alpha-2}u_n|v_n|^\beta h\big|dx=o(1)\|h\|.
\ee
For any $h\in D_{0}^{1,2}(\Omega)$, since $a_{\varepsilon_n}(x)\leq a_0(x)$, we have $\displaystyle \int_\Omega a_{\varepsilon_n}|h|^{2^*(s_2)}dx\leq \int_\Omega a_{0}|h|^{2^*(s_2)}dx.$ Then by the Hardy-Sobolev inequality, we obtain that there exists some $C>0$ independent of $n$ such that
\be\lab{2014-12-6-e54}
\Big(\int_\Omega a_{\varepsilon_n}|h|^{2^*(s_2)}dx\Big)^{\frac{1}{2^*(s_2)}}\leq C\|h\|.
\ee
Noting that $\frac{\alpha-1}{\alpha}+\frac{\beta}{2^*(s_2)\alpha}+\frac{1}{2^*(s_2)}=1$, by H\"older inequality and \eqref{2014-12-6-e54}, we have
\begin{align}\lab{2014-12-6-e55}
&\int_\Omega \big|a_{\varepsilon_n}(x)|u_n|^{\alpha-2}u_n|v_n|^\beta h\big|dx\nonumber\\
&\leq \Big(\int_\Omega a_{\varepsilon_n}|u_n|^\alpha|v_n|^\beta dx\Big)^{\frac{\alpha-1}{\alpha}}
  \Big(\int_\Omega a_{\varepsilon_n}(x)|v_n|^{2^*(s_2)}dx\Big)^{\frac{\beta}{2^*(s_2)\alpha}}\;\Big(\int_\Omega a_{\varepsilon_n}|h|^{2^*(s_2)}dx\Big)^{\frac{1}{2^*(s_2)}}\nonumber\\
&= o(1)\|h\|,
\end{align}
which means that \eqref{2014-12-6-e53} is proved.
Since $(u_n,v_n)$ is a positive ground state solution of  the system \eqref{zz-620-1}, we have
\be\lab{2014-12-6-e56}
\big\langle\Phi'_{\varepsilon_n}(u_n,v_n), (h,0)\big\rangle\equiv 0.
\ee
That is,
\be\lab{2014-12-6-e57}
\int_\Omega \Big[\nabla u_n\nabla h-\lambda \frac{|u_n|^{2^*(s_1)-1}u_n}{|x|^{s_1}}h-\kappa \alpha \int_\Omega a_{\varepsilon_n}(x)|u_n|^{\alpha-2}u_n|v_n|^\beta h \Big]dx\equiv 0\ee
for all $n$ and $h\in D_{0}^{1,2}(\Omega)$.  Then by \eqref{2014-12-6-e53} and \eqref{2014-12-6-e57}, we obtain that
\be\lab{2014-12-6-e58}
\int_\Omega \Big[\nabla u_n\nabla h-\lambda \frac{|u_n|^{2^*(s_1)-1}u_n}{|x|^{s_1}}h\Big]dx=o(1)\|h\|.
\ee
 Hence, $\Psi'_\lambda(u_n)\rightarrow 0$ in $H^{-1}(\Omega)$.
 Similarly, we can prove that $\Psi'_\mu(v_n)\rightarrow 0$ in $H^{-1}(\Omega)$.
\ep
\bc\lab{2014-12-6-pcro1}
$\displaystyle\liminf_{n\rightarrow +\infty}\int_\Omega a_{\varepsilon_n}|u_n|^\alpha|v_n|^\beta>0$.
\ec
\bp
By Lemma \ref{2014-12-5-l4}, we see that
\be\lab{2014-12-6-e51}
\|u_n\|^2+\|v_n\|^2\geq \delta_{\varepsilon_n}^{2}\geq \delta_0^2>0.
\ee
Hence, we obtain that
\be\lab{2014-12-6-e52}
(u_n,v_n)\not\rightarrow (0,0)\;\hbox{in}\;\mathscr{D}.
\ee
If $\displaystyle\liminf_{n\rightarrow +\infty}\int_\Omega a_{\varepsilon_n}|u_n|^\alpha|v_n|^\beta=0$, by Lemma \ref{2014-12-6-wl4}, we obtain that $\Psi'_\lambda(u_n)\rightarrow 0, \Psi'_\mu(v_n)\rightarrow 0$ in $H^{-1}(\Omega)$. Since either $u_n\not\rightarrow 0$ or $v_n\not\rightarrow 0$, it is easy to see that either
$\displaystyle \lim_{n\rightarrow +\infty}\Psi_\lambda(u_n)\geq m_\lambda$
or
$\displaystyle \lim_{n\rightarrow +\infty}\Psi_\mu(v_n)\geq m_\mu.$
By the assumption of $\displaystyle\liminf_{n\rightarrow +\infty}\int_\Omega a_{\varepsilon_n}|u_n|^\alpha|v_n|^\beta=0$ again, we have
\be\lab{2014-12-6-mle1}
\lim_{n\rightarrow +\infty}\Phi_{\varepsilon_n}(u_n,v_n)=\lim_{n\rightarrow +\infty} \Psi_\lambda(u_n)+\Psi_\mu(v_n)\geq \min\{m_\lambda,m_\mu\},
\ee
a contradiction to Lemma \ref{2014-12-6-wl3}. Thereby this corollary is proved.
\ep


\bl\lab{2014-12-9-l1}
Assume that $\{(\phi_n,\psi_n)\}$ is a bounded sequence of $\mathscr{D}$ such that
\be\lab{2014-12-9-e1}
\lim_{n\rightarrow +\infty}J_{\varepsilon_n}(\phi_n,\psi_n)=0
\ee
and
\be\lab{2014-12-9-e2}
\liminf_{n\rightarrow +\infty}c_{\varepsilon_n}(\phi_n,\psi_n)>0,
\ee
where   the functionals $ J_{\varepsilon}(u, v)$ and  $c_{\varepsilon}(u, v)$ are defined in  \eqref{zz-71-2} and \eqref{zz-71-1}, respectively.
Then
\be\lab{2014-12-9-e3}
\liminf_{n\rightarrow +\infty}\Phi_{\varepsilon_n}(\phi_n,\psi_n)\geq c_0.
\ee
\el
\bp
Since $\displaystyle \liminf_{n\rightarrow +\infty}c_{\varepsilon_n}(\phi_n,\psi_n)>0$, we see that $(\phi_n,\psi_n)\not\rightarrow (0,0)$ in $\mathscr{D}$. Without loss of generality, we
may assume that $(\phi_n,\psi_n)\neq (0,0)$ for all $n$. Combining with the boundedness of $\{(\phi_n,\psi_n)\}$, we obtain that there exists some $d_0, d_1>0$ such that
\be\lab{2014-12-9-e4}
0< d_0\leq a(\phi_n,\psi_n):=\|\phi_n\|^2+\|\psi_n\|^2\leq d_1\;\hbox{for all}\;n.
\ee
We also claim that $\displaystyle b(\phi_n,\psi_n)$ is bounded  and away from $0$,   i.e., there exists $d_3,d_4>0$ such that
\be\lab{2014-12-9-e5}
0<d_3\leq b(\phi_n,\psi_n):=\lambda |\phi_n|_{2^*(s_1),s_1}^{2^*(s_1)}+\mu |\psi_n|_{2^*(s_1),s_1}^{2^*(s_1)}\leq d_4.
\ee
The  right-hand inequality in   \eqref{2014-12-9-e5} is trivial due to the Hardy-Sobolev inequality. Now, we only need to prove the existence of $d_3$. We proceed by contradiction. If $b(\phi_n,\psi_n)\rightarrow 0$ up to a subsequence, then $\phi_n\rightarrow 0, \psi_n\rightarrow 0$ strongly in $L^{2^*(s_1)}(\Omega, \frac{dx}{|x|^{s_1}})$. Recalling the boundedness of $\{(\phi_n,\psi_n)\}$ again, by Proposition \ref{2014-5-5-interpolation-corollary} and Proposition \ref{2014-4-22-interpolation-corollary}, we obtain that
$\phi_n\rightarrow 0, \psi_n\rightarrow 0$ strongly in $L^{2^*(s_2)}(\Omega, a_0(x)dx)$. Noting that $a_{\varepsilon_n}(x)\leq a_0(x)$, then by the H\"older inequality, it is easy to prove that
\be\lab{2014-12-9-e6}
c_{\varepsilon_n}(\phi_n,\psi_n)\leq c_0(\phi_n,\psi_n)\rightarrow 0,
\ee
a contradiction to \eqref{2014-12-9-e2} and thereby \eqref{2014-12-9-e5} is proved.
We also note that $c_{\varepsilon_n}(\phi_n,\psi_n)$ is bounded. Hence, up to a subsequence, we may assume that
\be\lab{2014-12-9-e7}
a(\phi_n,\psi_n)\rightarrow a^*>0, b(\phi_n,\psi_n)\rightarrow b^*>0, c_{\varepsilon_n}(\phi_n,\psi_n)\rightarrow c^*>0.
\ee
Then by \eqref{2014-12-9-e1}, we see that
\be\lab{2014-12-9-e8}
a^*-b^*-2^*(s_2)\kappa c^*=0.
\ee
On the other hand, for any $n$, by Lemma \ref{2014-12-5-l1}, there exists a  unique $t_n>0$ such that $J_{\varepsilon_n}(t_n\phi_n,t_n\psi_n)=0$ and $t_n$ is implicity given by the following equation
\be\lab{2014-12-9-e9}
a(\phi_n,\psi_n)-b(\phi_n,\psi_n)t_{n}^{2^*(s_1)-2}-2^*(s_2)\kappa c_{\varepsilon_n}(\phi_n,\psi_n)t_{n}^{2^*(s_2)-2}=0.
\ee
Since $a(\phi_n,\psi_n)$ is bounded and $\displaystyle \liminf_{n\rightarrow +\infty}c_{\varepsilon_n}(\phi_n,\psi_n)>0$,
it is easy to see that $t_n$ is bounded.
On the other hand, by the Hardy-Sobolev inequality again, we obtain that
\be\lab{2014-12-9-e10}
a(\phi_n,\psi_n)\leq C_1 \big(a(\phi_n,\psi_n)\big)^{\frac{2^*(s_1)}{2}}t_{n}^{2^*(s_1)-2}+C_2 \big(a(\phi_n,\psi_n)\big)^{\frac{2^*(s_2)}{2}}t_{n}^{2^*(s_2)-2},
\ee
for some positive constants $C_1,C_2$ independent of $n$. Then it is easy to see that at least one of the following holds:
\begin{itemize}
\item[$(i)$]$\frac{1}{2}a(\phi_n,\psi_n)\leq C_1 \big(a(\phi_n,\psi_n)\big)^{\frac{2^*(s_1)}{2}}t_{n}^{2^*(s_1)-2}$;
\item[$(ii)$]$\frac{1}{2}a(\phi_n,\psi_n)\leq C_2 \big(a(\phi_n,\psi_n)\big)^{\frac{2^*(s_2)}{2}}t_{n}^{2^*(s_2)-2}$.
\end{itemize}
Hence, we see that $t_n$ is bounded away from $0$. Up to a subsequence if necessary, we assume that $t_n\rightarrow t^*>0$.
Then we have that

$$\left.
\begin{array}{lll}
&J_{\varepsilon_n}(t_n\phi_n,t_n\psi_n)\equiv 0,\\
&\{(\phi_n,\psi_n)\}\;\hbox{is bounded in}\;\mathscr{D},\\
&t_n\rightarrow t^*>0,
\end{array}\right\}\Rightarrow \lim_{n\rightarrow +\infty}J_{\varepsilon_n}(t^*\phi_n,t^*\psi_n)=0.
$$
By \eqref{2014-12-9-e7}, we obtain that
\be\lab{2014-12-9-e11}
a^*(t^*)^2-b^*(t^*)^{2^*(s_1)}-2^*(s_2)\kappa c^*(t^*)^{2^*(s_2)}=0,\quad (t^*>0).
\ee
It is easy to see that  the algebraic equation  $a^*=b^* t^{2^*(s_1)-2}+c^* t^{2^*(s_2)-2}$ has a  unique  positive solution. Hence, by \eqref{2014-12-9-e8} and \eqref{2014-12-9-e11}, we obtain that $t^*=1$. Then recalling the boundedness of $\{(\phi_n,\psi_n)\}$ again, we see that
\be\lab{2014-12-9-e12}
\lim_{n\rightarrow +\infty}\Phi_{\varepsilon_n}(\phi_n,\psi_n)=\lim_{n\rightarrow +\infty}\Phi_{\varepsilon_n}(t_n\phi_n,t_n\psi_n).
\ee
By the definition of $t_n$, we see that $(t_n\phi_n,t_n\psi_n)\in \mathcal{N}_{\varepsilon_n}$. Hence, $\Phi_{\varepsilon_n}(t_n\phi_n,t_n\psi_n)\geq c_{\varepsilon_n}$. Then by \eqref{2014-12-9-e12} and Corollary \ref{2014-12-6-zcro1}, we obtain that
\be\lab{2014-12-9-e13}
\lim_{n\rightarrow +\infty}\Phi_{\varepsilon_n}(\phi_n,\psi_n)\geq \lim_{n\rightarrow +\infty}c_{\varepsilon_n}=c_0.
\ee
\ep

\bl\lab{2014-12-6-wcro1}
Assume $s_1,s_2\in (0,2), \lambda, \mu\in (0,+\infty), \kappa>0,\alpha>1,\beta>1$ and $ \alpha+\beta= 2^*(s_2)$. Let $\varepsilon\in (0,s_2)$, then any  solution $(u,v)$ of  the  system \eqref{zz-620-1}  satisfies
\be\lab{2014-12-6-e40}
\int_{\Omega\cap \mathbb{B}_1}\kappa a_\varepsilon(x)|u|^\alpha|v|^\beta dx=\int_{\Omega\cap \mathbb{B}_1^c}\kappa a_\varepsilon(x)|u|^\alpha|v|^\beta dx,
\ee
where $\mathbb{B}_1$ is the unit ball of  $\, \R^N$ entered at zero.
\el
\bp
Let
\be\lab{2014-12-6-e36}
G(x,u,v)=\frac{\lambda}{2^*(s_1)}\frac{|u|^{2^*(s_1)}}{|x|^{s_1}}
+\frac{\mu}{2^*(s_1)}\frac{|v|^{2^*(s_1)}}{|x|^{s_1}}
+\kappa a_\varepsilon(x)|u|^\alpha|v|^\beta.
\ee
Noting that
\be\lab{2014-11-27-e2}
\frac{\partial}{\partial x_i}a_\varepsilon(x)=\begin{cases}
-(s_2-\varepsilon)\frac{1}{|x|^{s_2+2-\varepsilon}}x_i\quad &\hbox{for}\;|x|<1,\\
-(s_2+\varepsilon)\frac{1}{|x|^{s_2+2+\varepsilon}}x_i\quad&\hbox{for}\;|x|>1,
\end{cases}
\ee
we have
\begin{align}\lab{2014-12-6-e37}
&x_i\cdot G_{x_i}(x,u,v)\\
=&\begin{cases}-s_1\Big[\frac{\lambda}{2^*(s_1)}\frac{|u|^{2^*(s_1)}}{|x|^{s_1}}
+\frac{\mu}{2^*(s_1)}\frac{|v|^{2^*(s_1)}}{|x|^{s_1}}\Big]\frac{x_i^2}{|x|^2}-(s_2-\varepsilon)\frac{\kappa|u|^\alpha|v|^\beta}{|x|^{s_2+2-\varepsilon}}x_i^2\;\;\hbox{ if}\;|x|<1,\\
-s_1\Big[\frac{\lambda}{2^*(s_1)}\frac{|u|^{2^*(s_1)}}{|x|^{s_1}}
+\frac{\mu}{2^*(s_1)}\frac{|v|^{2^*(s_1)}}{|x|^{s_1}}\Big]\frac{x_i^2}{|x|^2}-(s_2+\varepsilon)\frac{\kappa|u|^\alpha|v|^\beta}{|x|^{s_2+2+\varepsilon}}x_i^2\;\;\hbox{ if}\;|x|>1.
\end{cases}\nonumber
\end{align}
Hence, by \eqref{2013-10-26-e2}, we have
\begin{align}\lab{2014-12-6-e38}
&-2(N-s_1)\int_\Omega \big[\frac{\lambda}{2^*(s_1)}\frac{|u|^{2^*(s_1)}}{|x|^{s_1}}+\frac{\mu}{2^*(s_1)}\frac{|v|^{2^*(s_1)}}{|x|^{s_1}}\big]dx\nonumber\\
&-2(N-s_2)\int_\Omega \kappa a_\varepsilon(x)|u|^\alpha|v|^\beta dx
+(N-2)\int_\Omega (|\nabla u|^2+|\nabla v|^2)dx\\
=&2\varepsilon\int_{\Omega\cap \mathbb{B}_1}\kappa a_\varepsilon(x)|u|^\alpha|v|^\beta dx-2\varepsilon \int_{\Omega\cap \mathbb{B}_1^c}\kappa a_\varepsilon(x)|u|^\alpha|v|^\beta dx.\nonumber
\end{align}
On the other hand, since $(u,v)$ is a solution of the  system \eqref{zz-620-1}, we have
\begin{align}\lab{2014-12-6-e39}
&\int_\Omega \big(|\nabla u|^2+|\nabla v|^2\big)dx
=\int_\Omega \big(\lambda\frac{|u|^{2^*(s_1)}}{|x|^{s_1}}+\mu\frac{|v|^{2^*(s_1)}}{|x|^{s_1}}+2^*(s_2)\kappa a_\varepsilon(x)|u|^\alpha|v|^\beta\big)dx.
\end{align}
Since $\varepsilon>0$, by \eqref{2014-12-6-e38} and \eqref{2014-12-6-e39}, we obtain the result of \eqref{2014-12-6-e40}.
\ep

\noindent{\bf Proof of Theorem \ref{2014-12-9-mainth1}:} By Corollary \ref{2014-12-6-wcro2}, we only need to prove that $(u_0,v_0)\neq (0,0)$.
Now,we proceed by contradiction. We assume that $(u_0,v_0)=(0,0)$.
By Corollary \ref{2014-12-6-pcro1}, we have, up to a subsequence if necessary, that
\be\lab{2014-12-9-xe2}
\lim_{n\rightarrow +\infty}\int_\Omega a_{\varepsilon_n}(x)|u_n|^\alpha|v_n|^\beta dx:=\tau>0.
\ee
On the other hand, by Corollary \ref{2014-12-6-wcro1}, we have
\be\lab{2014-12-9-xe3}
\int_{\Omega\cap {\mathbb{B}}_1} a_{\varepsilon_n}(x)|u_n|^\alpha|v_n|^\beta dx=\int_{\Omega\cap {\mathbb{B}}_1^c} a_{\varepsilon_n}(x)|u_n|^\alpha|v_n|^\beta dx\;\;\; \hbox{ for all }\;n.
\ee
Hence,
\be\lab{2014-12-9-xe4}
\lim_{n\rightarrow +\infty}\int_{\Omega\cap {\mathbb{B}}_1} a_{\varepsilon_n}(x)|u_n|^\alpha|v_n|^\beta dx=\lim_{n\rightarrow +\infty}\int_{\Omega\cap {\mathbb{B}}_1^c} a_{\varepsilon_n}(x)|u_n|^\alpha|v_n|^\beta dx=\frac{\tau}{2}>0.
\ee

Let $\chi(x)\in C_c^\infty(\R^N)$ be a cut-off function such that $\chi(x)\equiv 1$ in ${\mathbb{B}}_{\frac{1}{2}}$, $\chi(x)\equiv 0$ in $\R^N\backslash {\mathbb{B}}_{1}$ and take $\tilde{\chi}(x)\in C^\infty(\R^N)$ such that $\tilde{\chi}(x)\equiv 0$ in ${\mathbb{B}}_1$ and $\tilde{\chi}\equiv 1$ in $\R^N\backslash {\mathbb{B}}_{2}$.
Denote
\be\lab{2014-12-9-xe5}
\begin{cases}
\phi_{1,n}(x):=\chi(x)u_n(x),\; \phi_{2,n}(x):=\tilde{\chi}(x)u_n(x),\\
\psi_{1,n}(x):=\chi(x)v_n(x),\;  \psi_{2,n}(x):=\tilde{\chi}(x)v_n(x).
\end{cases}
\ee
Recalling that $(u_n,v_n)$ is a positive ground state solution of the   system \eqref{zz-620-1} with $\varepsilon=\varepsilon_n,$ then
\be\lab{2014-12-9-bue1}
\big\langle \Phi'_{\varepsilon_n}(u_n,v_n), (u_n-\phi_{1,n}-\phi_{2,n}, v_n-\psi_{1,n}-\psi_{2,n}) \big\rangle\equiv 0\;\hbox{for all}\;n.
\ee
when $(u_0,v_0)=(0,0)$, by Rellich-Kondrachov compactness theorem, if $0\not\in \overline{\tilde{\Omega}}$,  we have that $u_n\rightarrow 0, v_n\rightarrow 0$ strongly in $L^{2^*(s_1)}(\tilde{\Omega}, \frac{dx}{|x|^{s_1}})$ and $L^{2^*(s_2)}(\tilde{\Omega}, a_{\varepsilon_n}dx)$ uniformly for all $n$.
Hence, it is easy to prove that
\be\lab{2014-12-9-bue2}
(u_n-\phi_{1,n}-\phi_{2,n}, v_n-\psi_{1,n}-\psi_{2,n})\rightarrow (0,0)\;\hbox{strongly in}\;\mathscr{D}.
\ee
We also have that
\be\lab{2014-12-9-xe6}
\big\langle \Phi'_{\varepsilon_n}(u_n,v_n), (\phi_{1,n}, \psi_{1,n}) \big\rangle\equiv 0\;\hbox{for all}\;n.
\ee
Then by $(u_0,v_0)=(0,0)$ and Rellich-Kondrachov compactness theorem again,  it is easy to see that
\be\lab{2014-12-9-xe7}
\lim_{n\rightarrow +\infty}J_{\varepsilon_n}(\phi_{1,n},\psi_{1,n})=0.
\ee
We also obtain that
\be\lab{2014-12-9-xe8}
\liminf_{n\rightarrow +\infty}c_{\varepsilon_n}(\phi_{1,n}, \psi_{1,n})=\lim_{n\rightarrow +\infty}\int_{\Omega\cap {\mathbb{B}}_1} a_{\varepsilon_n}(x)|u_n|^\alpha|v_n|^\beta dx=\frac{\tau}{2}>0.
\ee
Hence, by Lemma \ref{2014-12-9-l1}, we have
\be\lab{2014-12-9-xe9}
\lim_{n\rightarrow +\infty}\Phi_{\varepsilon_n}(\phi_{1,n}, \psi_{1,n})\geq c_0.
\ee
Similarly, we can prove that
\be\lab{2014-12-9-xe10}
\lim_{n\rightarrow +\infty}\Phi_{\varepsilon_n}(\phi_{2,n}, \psi_{2,n})\geq c_0.
\ee
By $(u_0,v_0)=(0,0)$ and the Rellich-Kondrachov compact theorem again, we have that
\be\lab{2014-12-9-xe11}
c_{\varepsilon_n}=\lim_{n\rightarrow +\infty}\Phi_{\varepsilon_n}(u_n,v_n)=\lim_{n\rightarrow +\infty}\big[\Phi_{\varepsilon_n}(\phi_{1,n},\psi_{1,n})+\Phi_{\varepsilon_n}(\phi_{2,n},\psi_{2,n})\big].
\ee
Then by \eqref{2014-12-9-xe9},\eqref{2014-12-9-xe10} and Corollary \ref{2014-12-6-zcro1}, we obtain that
\be\lab{2014-12-9-xe12}
c_0\geq 2c_0,
\ee
a contradict to the fact of that $c_0>0$.
Hence, $(u_0,v_0)\neq (0,0)$ and the proof is completed.

\vskip0.26in


\end{document}